\documentclass{siamart190516}
\usepackage{graphicx}
\usepackage{booktabs}
\usepackage{amsmath}
\usepackage[english]{babel}
\usepackage[utf8]{inputenc}
\usepackage{amssymb}
\usepackage{multicol}
\usepackage[utf8]{inputenc}
\usepackage[english]{babel}
\usepackage{multirow}
\usepackage{hyperref,enumerate}
\usepackage{cancel} 
\usepackage{caption}
\usepackage{subcaption}
\usepackage{bbm}
\newcommand{\R}{\mathbb{R}}

\newcommand*\samethanks[1][\value{footnote}]{\footnotemark[#1]}

\graphicspath{{Figures/}}
\usepackage{ifthen}
\usepackage{float}
\floatstyle{ruled}
\newfloat{algorithm}{htb}{alg}%[section]
\floatname{algorithm}{Algorithm}
\newcounter{algo@row}
\newcounter{algo@rowindent}
\newcommand{\algofont}[1]{\textbf{#1}}% S1
\newcommand{\algonumbersize}[1]{\scriptsize{#1}}% S2
\newcommand{\algopreitem}[1][\arabic{algo@row}]{\texttt{\algonumbersize{#1}}}
\newcommand{\algoitemskip}{\hspace{\value{algo@rowindent}cc}}
\newcommand{\algonewnestedopen}[2]{
	\newcommand{#1}[1][]{%
		\ifthenelse{\equal{##1}{}}{\item}{\item[{\algopreitem[##1]}]}
		\algoitemskip\algofont{#2}%
		\addtocounter{algo@rowindent}{1}%
		\ignorespaces
	}
}
\newcommand{\algonewnestedaux}[2]{
	\newcommand{#1}[1][]{
		\addtocounter{algo@rowindent}{-1}
		\ifthenelse{\equal{##1}{}}{\item}{\item[{\algopreitem[##1]}]}
		\algoitemskip\algofont{#2}%
		\addtocounter{algo@rowindent}{+1}%
		\ignorespaces
	}
}
\newcommand{\algonewnestedclose}[2]{
	\newcommand{#1}[1][]{
		\addtocounter{algo@rowindent}{-1}
		\ifthenelse{\equal{##1}{}}{\item}{\item[{\algopreitem[##1]}]}
		\algoitemskip\algofont{#2}%
		\ignorespaces
	}
}
\newcommand{\algonewcommand}[2]{
	\newcommand{#1}[1][default]{
		\ifthenelse{\equal{##1}{default}}{\item}{\item[{\algopreitem[##1]}]}%
		\algoitemskip\algofont{#2}%
		\ignorespaces
	}%
}
\newcommand{\algonewkeyword}[2]{\newcommand{#1}{\algofont{#2}}}
\algonewcommand{\STATE}{\ignorespaces}
\algonewcommand{\INPUT}{Input: }
\algonewcommand{\pINPUT}{\phantom{Input: }}
\algonewcommand{\COMPUTE}{Compute: }
\algonewcommand{\OUTPUT}{Output: }
\algonewcommand{\pOUTPUT}{\phantom{Output: }}
\DeclareMathOperator*{\argmin}{arg\,min}				   % arg min
\newcommand{\bGamma}{{\boldsymbol{\Gamma}}}

\newcommand{\bA}{{\bf A}}
\newcommand{\bB}{{\bf B}}

\newcommand{\bD}{{\bf D}}

\newcommand{\bF}{{\bf F}}

\newcommand{\bI}{{\bf I}}

\newcommand{\bL}{{\bf L}}
\newcommand{\bM}{{\bf M}}

\newcommand{\bQ}{{\bf Q}}
\newcommand{\bR}{{\bf R}}

\newcommand{\bU}{{\bf U}}
\newcommand{\bV}{{\bf V}}
\newcommand{\bW}{{\bf W}}
\newcommand{\bX}{{\bf X}}
\newcommand{\bY}{{\bf Y}}
\newcommand{\bZ}{{\bf Z}}

\newcommand{\bd}{{\bf d}}
\newcommand{\be}{{\bf e}}

\newcommand{\br}{{\bf r}}

\newcommand{\bu}{{\bf u}}

\newcommand{\bw}{{\bf w}}
\newcommand{\bx}{{\bf x}}
\newcommand{\by}{{\bf y}}
\newcommand{\bz}{{\bf z}}
\newcommand{\bzero}{{\bf0}}

\newcommand{\bbR}{\mathbb{R}}

\newcommand{\T}[1]{\boldsymbol{\mathcal{#1}}}
\newcommand{\bmat}[1]{\begin{bmatrix}#1 \end{bmatrix} }

\definecolor{edscol}{rgb}{0, 0.4, 0 }
\definecolor{sgcol}{rgb}{0.54, 0.17, 0.89}

\algonewnestedopen{\IF}{if }
\algonewnestedaux{\ELSEIF}{else if }
\algonewnestedaux{\ELSE}{else }
\algonewnestedclose{\ENDIF}{end if }
\algonewnestedopen{\FOR}{for }
\algonewnestedclose{\ENDFOR}{end for }
\algonewnestedopen{\WHILE}{while }
\algonewnestedclose{\ENDWHILE}{end while }
\algonewcommand{\BREAK}{break}%

\algonewkeyword{\For}{for }%
\algonewkeyword{\To}{to }%
\algonewkeyword{\Do}{do }%
\algonewkeyword{\If}{if }%
\algonewkeyword{\Then}{then }%
\algonewkeyword{\Else}{else }%
\algonewkeyword{\End}{end }%
\algonewkeyword{\AND}{and }%
\algonewkeyword{\True}{true }%
\algonewkeyword{\False}{false }%
\algonewkeyword{\Call}{call }%
\algonewkeyword{\irbleigs}{irbleigs }%
\algonewkeyword{\tridiag}{tridiag}%
\algonewkeyword{\reorth}{reorth}%
\newcommand{\TheTitle}{  \MakeLowercase{$\boldsymbol{\ell_p}$-$\boldsymbol{\ell_q}$}} 
\newcommand{\ShortTitle}{Edge preservation in dynamic inverse problems}

\newcommand{\TheAuthors}{M. Pasha, A. K. Saibaba, S. Gazzola, M. I. Espa\~nol, and E. de Sturler}

\headers{\ShortTitle}{\TheAuthors}
\title{{\TheTitle}\thanks{Edge-preserving MM Krylov methods}
	}

\title{Efficient edge-preserving methods for dynamic inverse problems}

\author{Mirjeta Pasha\thanks{School of Mathematical and Statistical Sciences, Arizona State University, United States of America}\and Arvind K. Saibaba\thanks{Department of Mathematics, North Carolina State University, United States of America}\and Silvia Gazzola\thanks{Department of Mathematical  Sciences, University of Bath, United Kingdom }\and Malena I. Espa\~nol\samethanks[1]\and Eric de Sturler\thanks{Department of Mathematics, Virginia Tech, United States of America}}

\begin{document}
\maketitle
\begin{abstract} We consider efficient methods for computing solutions to dynamic inverse problems, where both the quantities of interest and the forward operator (measurement process) may change at different time instances but we want to solve for all the images simultaneously. 
We are interested in large-scale ill-posed problems that are made more challenging by their dynamic nature and, possibly, by the limited amount of available data per measurement step.
To remedy these difficulties, we apply regularization methods that enforce simultaneous regularization in space and time (such as edge enhancement at each time instant and proximity at consecutive time instants) and achieve this with low computational cost and enhanced accuracy. More precisely, we develop iterative methods based on a majorization-minimization (MM) strategy with quadratic tangent majorant, which allows the resulting least squares problem to be solved with a generalized Krylov subspace (GKS) method; the regularization parameter can be defined automatically and efficiently at each iteration. 

Numerical examples from a wide range of applications, such as limited-angle computerized tomography (CT), space-time image deblurring, and photoacoustic tomography (PAT), illustrate the effectiveness of the described approaches.  
\end{abstract}


\section{Introduction}\label{sec:intro}
In the classical setting, inverse problems are commonly formulated as static, where the underlying parameters that define the problem do not change during the measurement process. There exists a very rich literature and many numerical methods for this setting; see \cite{engl1996regularization, hansen1998rank, kaipio2006statistical, mueller2012linear, vogel2002computational} and the references therein. Motivated by new developments in science and engineering applications, dynamic inverse problems have recently obtained considerable attention, shifting the focus of the research community to the latter, where time-dependent information
needs to be recovered from time-dependent data.
Such applications include dynamical impedance tomography \cite{schmitt2002efficient, schmitt2002efficient2}, process tomography \cite{beck2012process}, undersampled dynamic x-ray tomography \cite{burger2017variational}, and passive seismic tomography \cite{westman2012passive, zhang2009passive}, to mention a few.
A common question of interest is the reconstruction of non-stationary objects from time-dependent projection measurements. For instance, moving objects during a CT scan lead to artifacts in the stationary reconstruction even if the change in time is small. More specifically, in the imaging of organs like heart and lungs, small changes from the heart beat or changes in the lungs during breathing can significantly affect the quality of the reconstructed solution. In \cite{achenbach2001detection, blondel20043d, li2005motion, trabold2003estimation}, approaches to reconstruct a static image from dynamic data are discussed. In~\cite{burger2017variational}, the authors discuss the reconstruction of dynamic data in space and time. More recent work on computationally feasible methods in the Bayesian framework for dynamic inverse problems is presented in \cite{chung2018efficient} and the quantification of the uncertainties was discussed in~\cite{saibaba2020efficient}. In this work, we are interested in similar scenarios where the target of interest changes in space and time, and we are not limited to any specific motion of the objects during the measurement process. Furthermore, we seek to preserve the edges of the desired solution. Edge preserving reconstruction is a technique to smooth images while preserving edges, which has been employed in many fundamental applications in image processing such as artifact removal \cite{yang2005adaptive}, denoising \cite{guo2018edge, rudin1992nonlinear,  vogel1998fast}, image segmentation \cite{cui2017classification, he2016image}, and feature selection \cite{yin2007total}. Despite its important role and wide use, edge-preserving filtering is still an active area of research. 
The methods that we propose rely on total variation (TV)-type regularization. 
There has been considerable work on edge preserving methods, but there are only a few contributions on edge preserving methods for dynamic inverse problems. The latter have been mostly developed in the recent years highlighting the need to develop methods for dynamic inverse problems in parallel with recent advancements in science and technology. We first mention  \cite{semerci2014tensor} where an iterative reconstruction method is presented for solving the multi-energy CT problem where the multi-spectral unknown is modeled as a low rank 3-way tensor that is combined
with total variation regularization to enhance the regularization
capabilities especially at low energy images where the effects of
noise are most notable. More recently, a tensor low-rank and sparse
representation model for moving object detection was proposed in \cite{hu2016moving}, where the nuclear norm constraint is used to exploit the spatio-temporal redundancy of
the background. In \cite{gao2019iterated}, a framework for solving state estimation problems with an additional sparsity-promoting $\ell_1$-regularization term is presented.

\subsection{Background on dynamic inverse problems}\label{ssec:background1}
\begin{figure}[h!]
     \centering
     \includegraphics[scale=0.30]{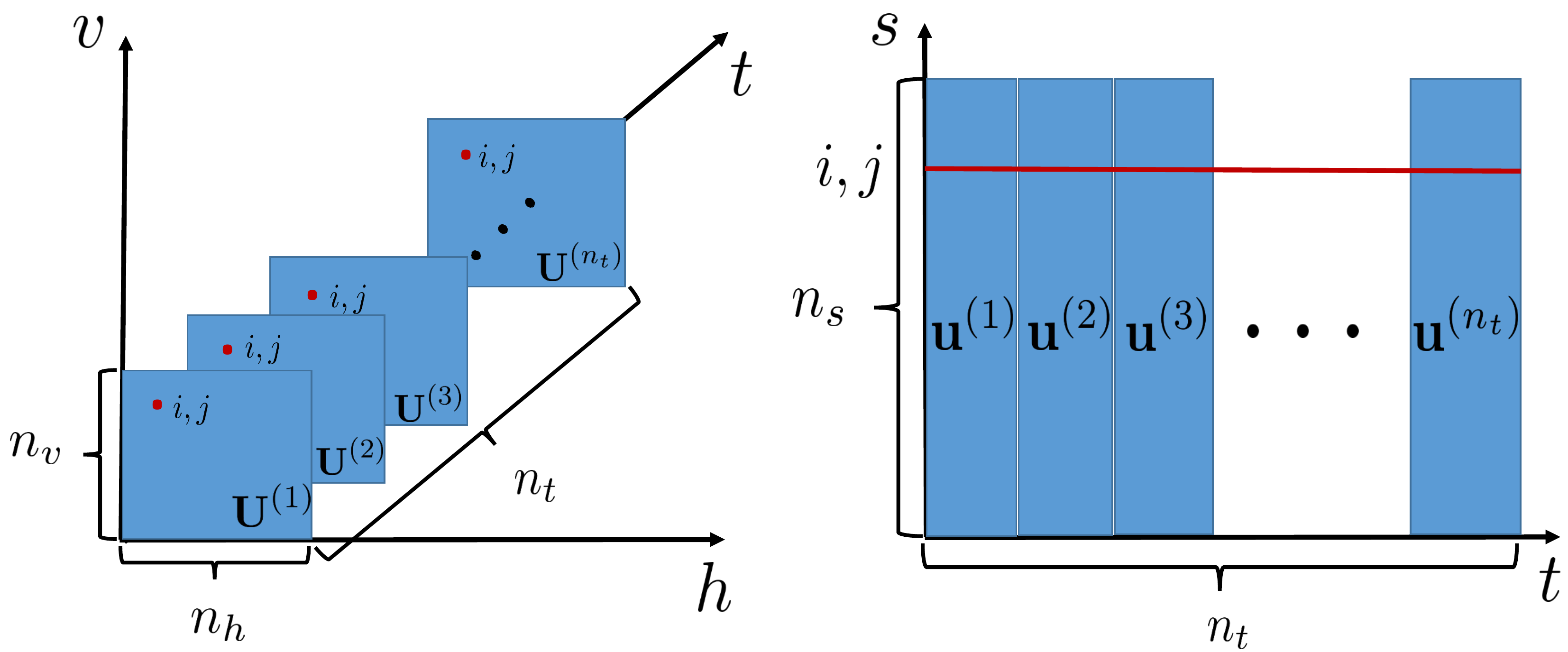}
\caption{Images $\bU^{(t)}$ to be reconstructed with pixels $i,j$ in red (left), and their corresponding vectorization $\bu^{(t)}$, which are the columns of the matrix $\bU$ where the pixels $i,j$ are now in the same row (right).}
\label{fig:urepr}
\end{figure}

First we set some notation that we use throughout the paper. 
Let $\bU^{(t)} \in \mathbb{R}^{n_v\times n_h}$ be the 2D (matrix) representation of an image with $n_v$ rows and $n_h$ columns at time instance $t = 1,2,\dots, n_t$. Define $n_s=n_vn_h$ and let $\bu^{(t)}$ be the column vector
by a lexicographical ordering of the two-dimensional $\bU^{(t)}$, that is, $\bu^{(t)}={\rm{vec}}(\bU^{(t)})\in \mathbb{R}^{n_s}$, with ${\rm vec}$ being the operation that vectorizes a matrix by stacking its columns. Then, let
$
\bU =\left[\bu^{(1)},\dots,\bu^{(n_t)}\right]\in \mathbb{R}^{n_s\times n_t}
$
be such that $\bu={\rm{vec}}(\bU)\in\mathbb{R}^{n}$ and $n=n_sn_t$. A pictorial representation of these quantities is displayed in Figure \ref{fig:urepr}.

We are interested in solving inverse problems with an unknown target of interest in space and time, where the goal is to recover the parameters $\bu^{(t)} \in \mathbb{R}^{n_s}$ that, in our case, represent pixels in the image to be reconstructed from the available measurements $\bd^{(t)} \in \mathbb{R}^{m_t}$ for $t = 1, 2,\dots, n_t$. Since we focus on imaging applications, we use the term `pixels' (rather than `parameters') throughout the paper. The value $n_t$ represents the number of time points and, therefore, the value $m = \sum_{t=1}^{n_t} m_t$ is the total number of available measurements. We consider the number of pixels, $n_s$, to be fixed for all time points. 
Dynamic problems may also involve reconstructing
a sequence of 
images with varying numbers of pixels (e.g., in image registration), but we do not consider that setting in this paper. 
For completeness, we define static and dynamic inverse problems in the context of this paper.

\paragraph{1. Dynamic inverse problems}    In a dynamic inverse problem, both the object of interest and the measurement process are known to change in time and, therefore, combining prior information at different time instances is found to enhance the reconstruction and recover dynamic information about the objects of interest. More specifically, we have the measurement equation
\begin{equation}\label{eq: minEq}
\bd = \bF \bu  + \be \mbox,
\end{equation}
where we consider two cases for the parameter-to-observation map $\bF \in \R^{m \times n}$:
\begin{enumerate}
\item[a.]  The forward operator is time-dependent, that is, it changes during the data acquisition process and results in the operator $\bF$ being a block diagonal matrix
\begin{equation}\label{eq: blockF} \bF= \begin{bmatrix} \bA^{(1)} \\ &  \ddots\\ & & \bA^{(n_t)}  \end{bmatrix}.\end{equation}
\item[b.] The forward operator is time-independent, i.e., $\bA^{(t)} = \bA$ for $t=1,\dots,n_t$, where $\bF$ simplifies to $\bF = \bI_{n_t} \otimes \bA$, with $\otimes$ being the Kronecker product. 
\end{enumerate}
The vector $\bd ={\rm{vec}}([\bd^{(1)},\dots,\bd^{(n_t)}])\in \R^{m}$ represents measured data that are contaminated by an (unknown) error $\be \in \R^{m}$ 
that may stem from measurement errors. 
We assume that the noise vector follows a normal Gaussian distribution with mean zero and covariance ${\boldsymbol{\Gamma}}$, i.e., $\be \sim \mathcal{N}(\bf{0}, {\boldsymbol{\Gamma}})$\footnote{Throughout this paper we use $\mathcal{N}({\boldsymbol{\mu}}, {\boldsymbol{\Gamma}})$ to denote the multivariate normal distribution with mean ${\boldsymbol {\mu}}$ and covariance matrix ${\boldsymbol {\Gamma}}$.}. 
The inverse problem involves recovering the pixels $\bu$ from the data $\bd$. That is, we seek to solve the general regularized problem
\begin{equation}\label{eq: l2-lq}
\bu_\text{dynamic}=\argmin_{\bu \in \R^{n}}{\mathcal{J}(\bu)} := \mathcal{F}(\bu)+\lambda\mathcal{R}(\bu),
\end{equation}
where the functional $\mathcal{F}(\bu)$ is a data-misfit term that takes the form $\frac{1}{2}\|\bF\bu-\bd\|_{\boldsymbol{\Gamma}^{-1}}^2$ and the term $\mathcal{R}(\bu)$ is a regularization term that can take different forms; several new forms for $\mathcal{R}(\bu)$ will be proposed in Section~\ref{sec: newmethods}.
\paragraph{2. Static inverse problems}By contrast, in a static inverse problem, the information from the current time step is used to reconstruct the unknown pixels $\bu^{(t)}$, $t = 1,2,\dots, n_t$. We assume that the measurement noise at each time step is independent of other time steps, so that the overall noise covariance matrix, defined above, \linebreak[4]$\bGamma = \text{BlockDiag}(\bGamma_1,\dots,\bGamma_{n_t})$ is a block-diagonal matrix, where $\bGamma_t$ is the noise covariance matrix at step $t$.  We then solve the sequence of optimization problems
\begin{equation}\label{eq: static}
\bu^{(t)}_\text{static} = \argmin_{\bu \in \R^{n_s}} \frac12 \| \bA^{(t)}\bu -\bd^{(t)}\|_{\bGamma_t^{-1}}^2 + \lambda \mathcal{R}(\bu), \qquad t = 1,2,\dots, n_t 
\end{equation}
independently, to obtain the solution to the static inverse problem. Throughout this paper, $\lambda > 0$ is an appropriate regularization parameter that determines a balance between the data-misfit and the regularization terms.

\paragraph{Challenges  and need for regularization} The considered inverse problems are typically ill-posed since a solution may not exist, may not be unique, or may not depend continuously on the data. A clear first challenge in dynamic inverse problems stems from the limited information available per time instance during the measurement process.
Moreover, when solving dynamic inverse problems, the unknown has $n = n_sn_t$ pixels, which can be as high as $\mathcal{O}(10^6)$: a clear second challenge is therefore the large-scale of the considered problems. 

In order to obtain meaningful solutions of  ill-posed inverse problems one must resort to regularization. In this paper, we focus on developing efficient regularization approaches for dynamic inverse problems that promote edge-preservation in the resulting images by incorporating specific representations of the prior information. Namely, we propose a combination of spatial and temporal prior information representations that allow for the recovery of piecewise constant solutions and the use of efficient numerical methods that can exploit these representations.

\subsection{Overview of the main contributions}
Our main goal in this paper is to provide a suite of techniques for edge-preserving reconstructions in dynamic inverse problems. 
We summarize the main contributions as follows:
\begin{enumerate}
\item We propose six different regularization techniques, involving anisotropic and isotropic total variation, and group sparsity, which promote edge-preserving reconstructions in dynamic inverse problems, {where the images to be reconstructed are changing in time}. 
We provide motivation for each technique, which combines spatiotemporal information in different ways. For each regularization technique, we also provide an interpretation using tensor notation.

\item For each regularization technique, we write down the corresponding optimization problem for reconstructing the desired solution, whose objective functions are convex but non-differentiable. To remedy the non-differentiability, we consider a smoothed functional instead, and we derive an iterative \break reweighted least squares (IRLS) approach for each optimization problem using the  majorization-minimization (MM) technique.

\item To efficiently solve the sequence of least squares problems and define the regularization parameter, we use a generalized Krylov subspace (GKS) method, resulting in a so-called MM-GKS method. For instance, typically we are able to reconstruct over $1.9$ million pixels in less than $100$ MM-GKS iterations.

\item We illustrate the performance and demonstrate the range of applicability and effectiveness of
the described approaches on a variety of test problems with simulated and real data arising from space-time image deblurring, photoacoustic tomography (PAT), and limited angle computerized tomography (CT). 

\end{enumerate}
In summary, we present a unified framework for edge-preserving reconstructions in dynamic inverse problems. The definition of the specific regularization terms allows to use the same MM-GKS technique to solve the corresponding optimization problems. The presented approaches are generic and extend to other problem settings such as
multichannel imaging \cite{kazantsev2018joint, li2010multichannel, thompson1985dynamic},
electroencephalographic current density
reconstruction \cite{giraldo2013weighted}, and anatomical image analysis to study changes in organ anatomy \cite{rohde2003comprehensive}; in all these applications, the solution techniques combine limited information from different sources to improve the quality of the resulting reconstruction and {recover}  dynamic information from different channels.

\emph{Overview of the paper.} This paper is organized as follows. In Section~\ref{sec: problem}, we present some background material, including additional notation, a survey of well-established regularization terms, and an iterative method used to solve the inverse problem by the aid of an MM strategy. In Section~\ref{sec: newmethods}, we propose six different methods for edge-preserving regularization in dynamic inverse problems, write a unifying framework and derive, by using an MM approach, an iteratively reweighted least squares method for solving the resulting optimization problem.  Some alternative approaches and extensions that fit within our framework are presented in Section~\ref{sec: alternative}. In Section \ref{sec: iterativemethods}, we describe iterative methods based on generalized Krylov subspaces to efficiently solve the resulting optimization problem and define the regularization parameter at each iteration. In Section~\ref{sec: nresults}, we present numerical examples that demonstrate the performance of the proposed regularization terms and the MM solvers. Finally, some conclusions, remarks, and future directions are presented in Section~\ref{sec: conclusions}.

\section{Background}\label{sec: problem}

In this section, we review known facts about tensors, regularization terms such as (discrete) isotropic and anisotropic total variation, and the majorization-minimization approach for solving optimization problems.

\subsection{Tensor notation} \label{ssec:tensor}
The use of tensor notation is very convenient for describing dynamic images. 
A \textit{tensor} $\T{X}$ is a multi-dimensional array (also called \textit{n-way} or \textit{n-mode} array), whose entries are scalars. The \emph{order} of a tensor refers to the number of ways or modes. For instance, vectors are tensors of order one, and matrices are tensors of order two. More details on tensors can be found in~\cite{kolda2009tensor}. 

In this work, we primarily focus on 3rd-order tensors $\T{X}  \in \mathbb{R}^{n_1 \times n_2 \times n_3}$ with entries $x_{i,j,k}$. 
 
Fibers are higher-order analogues of matrix rows and columns. A \emph{fiber} of a third order tensor is a vector that is obtained by fixing two of the indices of the tensor $\T{X}$. We define $\T{X}_{:, j, k}$, $\T{X}_{i,:, k}$, and $\T{X}_{i, j, :}$ to be mode-1, mode-2, and mode-3 fibers, respectively. We implicitly assume that once a mode fiber has been extracted, it is reshaped as a column vector. 
\emph{Slices} are  two dimensional sections of a tensor that are obtained by fixing one of the indices. We define $\T{X}_{i,,:}, \T{X}_{:,j,:},$ and $\T{X}_{:,:,k}$ to be horizontal, lateral, and frontal slices, respectively. As before, when a slice is extracted, we implicitly assume that it is a matrix. The \emph{mode-$j$ unfolding} or \emph{matricization} of a tensor $\T{X}$ is obtained by arranging the mode-$j$ fibers to be the columns of a resulting matrix. We denote these by $\bX_{(1)} \in \R^{n_1\times (n_2n_3)}, \bX_{(2)} \in \R^{n_2\times (n_1n_3)}$, and $\bX_{(3)} \in \R^{n_3\times (n_1n_2)}$.

Another important concept here is the \emph{mode-$j$ product} that defines the operation of multiplying a tensor $\T{X} \in \R^{n_1\times n_2 \times n_3}$ by a matrix $\bL_j \in \R^{r\times n_j}$ for $j=1,2,3$ given in the following definition. We write $\T{Y} = \T{X} \times_j\bL_j$ in terms of the mode unfoldings as $\bY_{(j)} = \bL_j\bX_{(j)}$. For distinct modes
in a series of multiplications, the order of the multiplication is irrelevant.

We will also need to use \emph{norms} for tensors, which we define entrywise. That is, for $q\in [1,\infty)$, we define
\begin{equation}\label{eq:tensornorm}
    \|\T{X}\|_q = \left(\sum_{i=1}^{n_1}\sum_{j=1}^{n_2}\sum_{k=1}^{n_3}|x_{i,j,k}|^q\right)^{1/q}.
\end{equation}

A tensor representation of the dynamic inverse problem solution described in Section \ref{ssec:background1} is obtained by defining the multidimensional array $\T{U} \in \R^{n_v\times n_h\times n_t},$ with its frontal slices taken to be 2D representations of the image $\bu^{(t)}$. That is, we let 
\begin{equation}\label{tensorU}
\T{U}_{:,:,t} = \text{mat}(\bu^{(t)}) \in \R^{n_v\times n_h} \qquad t=1,\dots,n_t\,.
\end{equation} 
  Furthermore, $\bu^{(t)}$ are the mode-3 fibers and $\bU = \bU_{(3)}^T$ is the transposed mode-$3$ unfolding.
\subsection{Regularization terms based on the first derivative operator}

When the desired solution is known to be piecewise constant with regular and sharp edges, total variation (TV) regularization is a popular choice, as it allows the solution to have discontinuities by preserving edges and discouraging noisy oscillations. TV regularization essentially enforces sparse gradient representations for the solution. 

The idea of (isotropic) TV regularization was first proposed in \cite{rudin1992nonlinear} for denoising, and it has then been used in compressive sensing \cite{candes2006robust}, segmentation \cite{sun2011image}, image deblurring \cite{rudin1994total}, electrical impedance tomography \cite{gonzalez2017isotropic},
and computerized tomography \cite{sidky2008image}, just to mention a few;  see \cite{caselles2015total} for a review.  The anisotropic TV formulation was addressed in \cite{choksi2010anisotropic, esedoglu2004decomposition, lou2015weighted}. 

Let 
\begin{equation}\label{D2}
	\bL_{d}= \alpha_d
	\begin{bmatrix} 
	1 &-1 & & &  \\
	&1 &-1 & &  \\
	& &\ddots  &\ddots  & \\
	 &  & &1  &-1 
	\end{bmatrix}\in \R^{(n_d-1)\times n_d} %\in\mathbb{R}.
	\quad
	\mbox{and}\quad
	\bI_{n_d}\in\mathbb{R}^{n_d\times n_d}
\end{equation}
be a rescaled finite difference discretization of the first derivative operator with $\alpha_d>0$ and the identity matrix of order $n_d$, respectively. In defining some of the operators below, we will augment 
the matrix $\bL_{d}$ with one zero row (at the bottom) and denote it by $\bL_{d}^0$. The matrices $\bL_d$ and $\bL_d^0$ are used to obtain discretizations of the first derivatives in the $d$-direction, with $d = v, h, t$ respectively (that is, multiplications by these matrices produce vectors containing approximations of vertical derivatives ($v$-direction), horizontal derivatives ($h$-direction), and derivatives with respect to time ($t$-direction)). {For simplicity, in the following, we let  $\alpha_d=1$, but different values can be used in practice: an $\alpha_d\neq 1$ can be treated as a regularization parameter that must be estimated as part of the inversion process.}

Considering only the spatial derivatives for now, these have the form  

\begin{equation}\label{eq: Deriv}
\begin{array}{lcl}
     {\rm vec}(\bL_{v}\bU^{(t)}) &=& (\bI_{n_h} \otimes \bL_{v}) \bu^{(t)} \in \mathbb{R}^{(n_{v} -1) n_h}\\
    {\rm vec}(\bU^{(t)}\bL^T_h) &=&(\bL_{h} \otimes \bI_{n_v}) \bu^{(t)} \in \mathbb{R}^{(n_{h} -1) n_v}
\end{array},
\quad t=1,\dots,n_t\,.
\end{equation}

Define also the spatial derivative matrix 
\begin{equation}\label{eq: Ls}
\bL_s = \bmat{ \bI_{n_h} \otimes \bL_{v} \\ \bL_{h} \otimes \bI_{n_v}}. 
\end{equation}
When time is considered, we have
\[
{\rm vec}(\bU\bL_{t}^T) =(\bL_{t}\otimes \bI_{n_s} ) \bu \in \mathbb{R}^{(n_{t} -1)n_s}.
\]
By letting $n_t = 1$ (i.e., $n = n_s$) for now, so that $\bu=\bu^{(1)}={\rm vec}(\bU^{(1)})$, we define the anisotropic total variation (TV$_{\rm aniso}$) as \\

\begin{align}\label{eq: anisoTV}
\mbox{TV}_{\rm aniso} (\bu) & = \sum_{k = 1}^{(n_v-1)}\sum_{\ell = 1}^{n_h} \left\vert \left(\bL_{v}\bU^{(1)}\right)_{k,\ell}\right\vert  + \sum_{k = 1}^{(n_h-1)}\sum_{\ell = 1}^{n_v}\left\vert \left(\bU^{(1)} \bL_{h}^T\right)_{k,\ell}\right\vert \nonumber \\
& = \|(\bI_{n_h} \otimes \bL_{v})\bu\|_1+ \|(\bL_{h} \otimes \bI_{n_v})\bu\|_1 = \|\bL_s\bu\|_1\,.
\end{align}
Assuming, for simplicity, that $n_h=n_v$, we define the isotropic total variation (TV$_{\rm iso}$) as
\begin{align}\label{eq: isoTV}
\mbox{TV}_{\rm iso} (\bu) & = \sum_{k=1}^{n_v}  \sum_{\ell=1}^{n_h} \sqrt{	(\bL_{v}^0\bU^{(1)})^2_{k,\ell} + (\bU^{(1)}(\bL_{h}^0)^T)^2_{k,\ell}} \nonumber \\
&=   \sum_{\ell=1}^{n_vn_h}\sqrt{	((\bI_{n_h} \otimes \bL^0_{v})\bu)^2_\ell + ((\bL^0_{h} \otimes \bI_{n_v})\bu)^2_\ell} \nonumber \\ 
& =\left\|\left[(\bI_{n_h} \otimes \bL^0_{v})\bu,\, (\bL^0_{h} \otimes \bI_{n_v})\bu\right]\right\|_{2,1},
\end{align} 
where $\|\cdot\|_{2,1}$ denotes the functional defined, for a matrix  $\bY
\in\mathbb{R}^{m_y\times n_y}$, as 
\[
\|\bY\|_{2,1}=\sum_{i=1}^{m_y}
\left(\sum_{j=1}^{n_y}(\bY)_{i,j}^2\right)^{1/2}=\sum_{i=1}^{m_y}\;
\|\bY_{i,:}\|_2\,.
\]
The $\|\cdot\|_{2,1}$ functional is typically used to enforce some kind of group sparsity. {The minimization of TV$_{\rm aniso}$ favors horizontal and vertical structures, since the edges not aligned with the coordinate axes are penalized heavily, so  TV$_{\rm iso}$ is used instead; however, this also has issues since  rotating an image by 90$^\circ$, or its multiples, changes the value of TV$_{\rm iso}$, which is undesirable~\cite{condat2017discrete}.}
\subsection{A majorization-minimization method}\label{sec: MM}
In this section, we provide an overview of the majorization-minimization technique for approximating the solution of \eqref{eq: l2-lq} by solving a sequence of optimization problems; 

see \cite{hunter2004tutorial, lange2016mm} for more details on the MM methods used. 

Suppose we want to minimize an objective function $\mathcal{J}(\bu)$. We shall need the following definition of a quadratic tangent majorant.
\begin{definition}[\cite{huang2017majorization}]\label{def: 1} Let $\by \in \R^n$ be fixed. The functional $\mathcal{Q} (\cdot; \by) \colon \R^{n} \rightarrow \R$ is said to be a 
quadratic tangent majorant for $\mathcal{J}(\bx)$ at $\bx=\by\in\R^n$ if it satisfies the following conditions:
\begin{enumerate}
\item $\mathcal{Q}(\bx;\by)$ is quadratic in $\bx$,
\item $\mathcal{Q}(\by;\by)=\mathcal{J}(\by)$,
\item $\bigtriangledown_{\bx}\mathcal{Q}(\by;\by)=\bigtriangledown_{\bx}\mathcal{J}(\by)$, 
%\AKS{Should the gradient be wrt $\by$}
\item $\mathcal{Q}(\bx;\by)\geq \mathcal{J}(\bx)\quad\forall \bx\in \R^{n}$.
\end{enumerate}
\end{definition}

The MM methods considered in this paper establish an iterative scheme whereby, starting from a given approximation of $\bu_{\rm true}$, a quadratic tangent majorant functional for $\mathcal{J}(\bu)$ at the previous iteration is defined and minimized to get the next approximation of $\bu_{\rm true}$. In other words, after the approximation $\bu_{(k)}$ has been computed at the $k$th iteration of the MM scheme, the $(k+1)$th approximate solution is computed as 
\[ \bu_{(k+1)} = \argmin_{\bu\in\R^n} \mathcal{Q}(\bu;\bu_{(k)}) \qquad k=0,1,\dots.  \]
The process is iterated. At the first iteration, one may take $\bu_{(0)}=\mathbf{0}$. 

The convergence of the MM approach with quadratic tangent majorants was established in~\cite{huang2017majorization}, which we use in this paper as well.  

\section{Dynamic edge-preserving regularization}\label{sec: newmethods}

We propose a unified framework with six main methods for edge preserving reconstruction applied to dynamic inverse problems with a spatial and time component. Furthermore, the framework is generic and can be extended to many other applications. For each technique, we motivate the kind of regularization, and using an MM approach we derive an iteratively reweighted least squares method for solving the resulting optimization problem. Finally, we also provide an interpretation for the regularization term using tensor notation. 
\subsection{Anisotropic space-time total variation (AnisoTV)} \label{subsec: R1a}
In this first technique we use the
summation of the TV of the images at each time step as a regularizer as well as regularization for temporal information.
Let $\bL_{s}$ be a matrix that represents the discretized finite difference operator corresponding to the spatial first derivative defined as in \eqref{eq: Ls}. The anisotropic TV terms $\|\bL_s\bu^{(t)}\|_1$,  $t=1,\dots, n_t$ ensure that the spatial discrete gradients of the images are sparse at each time instant. In addition, to incorporate temporal information, assuming that the images do not change considerably from one time instant to the next, we also want to penalize the  difference between any two consecutive images; we do so by considering the 1-norm differences  $\|\bu^{(t+1)}-\bu^{(t)}\|_1$ for $t=1,\dots,n_t-1$.   
These two requirements can be imposed by considering the following regularization term
\begin{equation}\label{eqn:R1}
\begin{aligned}
\mathcal{R}_{1}(\bu) &= \>  \sum_{t=1}^{n_t} \|\bL_s\bu^{(t)}\|_1 +  \sum_{t=1}^{n_t-1}\|\bu^{(t+1)}-\bu^{(t)}\|_1\\
&= \> \|(\bI_{n_t}\otimes\bL_s)\bu\|_1 + \|(\bL_t\otimes \bI_{n_s})\bu\|_1.
\end{aligned} 
\end{equation}

Notice that we can write the regularization terms that impose sparsity in time and in space under the same $\ell_1$ norm by letting 
\begin{equation}\label{eq: D1}
\bD_1 = \bmat{ \bI_{n_t} \otimes  \bL_s \\ \bL_t \otimes \bI_{n_s}} = \begin{bmatrix} \bI_{n_t}\otimes\bI_{n_h}\otimes\bL_{v} \\  
\bI_{n_t}\otimes \bL_{h}\otimes\bI_{n_v} \\
\bL_{t}\otimes\bI_{n_h}\otimes\bI_{n_v} 
\end{bmatrix}, \end{equation}
so that $\mathcal{R}_{1}(\bu) = \|\bD_1\bu\|_1$. 

To enable comparisons with the other methods proposed in this paper,  we provide an alternative representation for the regularization term in tensor notation. Recall the tensor representation $\T{U}$ of $\bu$ as in (\ref{tensorU}).  

Then we can write
\[\mathcal{R}_{\rm 1}(\bu) = \|\T{U} \times_1 \bL_{v}\|_1 + \|\T{U} \times_2 \bL_{h}\|_1 + \|\T{U} \times_3 \bL_{t}\|_1.\]

\paragraph{Optimization problem and MM approach} With the regularization term defined as in~\eqref{eqn:R1}, the optimization problem that we seek to solve takes the form
\begin{equation}\label{eq:R1}
\min_{\bu \in \R^{n}}\mathcal{J}_{1}(\bu):=  \mathcal{F}(\bu) +  \lambda \mathcal{R}_{1}(\bu), 
\end{equation}
where $\lambda > 0$.

We now derive an MM approach for solving this optimization problem by solving a sequence of simpler optimization problems whose closed form solutions exist. We do this in detail here, since the other regularization terms we propose have similar derivations. At the $k$th iteration of the MM method, let $\bz_{(k)}=\bD_1\bu_{(k)}$, where $\bu_{(k)}$ is the current iterate. Since the regularization term is nondifferentiable, we first majorize as 
\begin{equation}\label{R1epsilon}
    \mathcal{R}_{1}(\bu) \leq \sum_{\ell} \sqrt{(\bD_1\bu)_\ell^2 + \epsilon^2} =: R_{1\epsilon}(\bu),
\end{equation}
where $R_{1\epsilon}$ is the smoothed regularization term with $\epsilon > 0$. Similarly, we define the smoothed objective function $J_{1\epsilon}$, by replacing $\mathcal{R}_1(\bu)$ with $\mathcal{R}_{1\epsilon}(\bu)$ in~\eqref{eq:R1}.

To obtain a quadratic tangent majorant, we use the elementary inequality~\cite[Equation (1.5)]{lange2016mm} for $u,v > 0$ 
\begin{equation}\label{inequv}
    \sqrt{u} \leq \sqrt{v} + \frac{1}{2\sqrt{v}} (u-v);
\end{equation}
this is an equality if $u = v$. 
By applying (\ref{inequv}) to each term in the sum (\ref{R1epsilon}), with $u = (\bD_1\bu)_\ell^2 + \epsilon^2$ and $v = (\bD_1\bu_{(k)})_\ell^2 + \epsilon^2$, we obtain that
\[ \begin{aligned}\mathcal{R}_{1}(\bu) \leq & \>  \sum_{\ell}  \frac{1}{2\sqrt{(\bD_1\bu_{(k)})_\ell^2 + \epsilon^2 }} (\bD_1\bu)_\ell^2 + \tilde{c}_1 \\
= & \> \frac{1}{2} \|\bM_{\rm 1}^{(k)}\bu\|^{2}_{2}+\tilde{c}_1, 
\end{aligned}\] 
where $\tilde{c}_1$ is a constant independent of $\bu$ (but dependent on $\bu_{(k)}, \bD_1$, and $\epsilon$) and $\bM_1^{(k)}$ is the weighting matrix
\begin{equation}\label{eq:M1}
\bM_{1}^{(k)} := \bW_{\rm 1}^{(k)}\bD_1 \quad \mbox{ with } \quad \bW_{1}^{(k)} = \text{diag}((\bz_{(k)})^2+\epsilon^2)^{-1/4}).
\end{equation}  
Note that all operations in the expressions on the right-hand sides, including squaring, are performed element-wise.

We can now define the quadratic tangent majorant $\mathcal{Q}_{\rm 1}(\bu; \bu_{(k)})$ for the objective function $\mathcal{J}_{1\epsilon}(\bu)$ as
\begin{equation}\label{eq:Q1}
\mathcal{Q}_{\rm 1}(\bu; \bu_{(k)}):=\mathcal{F}(\bu)
+\displaystyle{\frac{\lambda}{2}} \|\bM_{\rm 1}^{(k)}\bu\|^{2}_{2}+c_1,\nonumber
\end{equation}
where $c_1 = \lambda\tilde{c}_1$.

Thus, as described in Section \ref{sec: MM}, we state the IRLS approach for solving the optimization problem~\eqref{eq:R1}: given an initial guess $\bu_{(0)}$, we solve the sequence of optimization problems 
\begin{equation}\label{eq:MMQ1}
    \bu_{(k+1)} = \argmin_{\bu\in \R^{n}} \mathcal{Q}_{1}(\bu;\bu_{(k)}), \qquad k = 0,1,2,\dots,
\end{equation}
to obtain the next iterate $\bu_{(k+1)}$. This can be interpreted as an iteratively reweighted least squares approach since, at each iteration, it replaces the regularization term $\mathcal{R}_{1\epsilon}(\bu)$ by an iteratively reweighted  $\ell_2$ regularization term.

\subsection{Total variation in space and Tikhonov in time (TVplusTikhonov)}\label{subsec: R1b}
In this technique, we consider anisotropic total variation in space and assume that the target of interest has small changes in time. Then, we define a new regularization term as 
\begin{equation}\label{eqn:R2}
\begin{aligned} 
\mathcal{R}_{2}(\bu) & := \>  \sum_{t=1}^{n_t} \|\bL_s\bu^{(t)}\|_1 + \sum_{t=1}^{n_t-1}\|\bu^{(t+1)}-\bu^{(t)}\|_2^2 \\
&= \>  \|(\bI_{n_t}\otimes\bL_s)\bu\|_1 +\|(\bL_t\otimes \bI_{n_s})\bu\|_2^2.
\end{aligned}
\end{equation}
In tensor notation, similar to $\mathcal{R}_1(\bu)$, we can succinctly write
\[\mathcal{R}_{\rm 2}(\bu) =\|\T{U} \times_1 \bL_{v}\|_1 + \|\T{U} \times_2 \bL_{h}\|_1 +  \|\T{U} \times_3 \bL_{t}\|_2^2.\]
Note that, when compared with $\mathcal{R}_1(\bu)$, $\mathcal{R}_{2}(\bu)$ requires the difference between the images at consecutive time steps to be small, while $\mathcal{R}_{1}(\bu)$ additionally promotes the sparsity of the difference.

\paragraph{Optimization problem and MM approach}{
We solve the inverse problem (\ref{eq: minEq}) by solving the optimization problem:
\begin{equation}\label{eq:8}
\min_{\bu \in \R^{n}}\mathcal{J}_{2}(\bu):=  \mathcal{F}(\bu)+{\lambda} \mathcal{R}_{2}(\bu), 
\end{equation}
where $\lambda > 0$. To achieve this, we can apply the MM approach similar to Section~\ref{subsec: R1a}. In particular, we consider the smoothed version $\mathcal{R}_{2\epsilon}(\bu)$ of $\mathcal{R}_2(\bu)$, where the smoothing is applied only to the first term in (\ref{eqn:R2}); the corresponding smoothed objective function is denoted by $\mathcal{J}_{2\epsilon}(\bu)$. To derive a majorant, we only need to majorize the first term, since the second term is already expressed in the squared $\ell_2$ norm. We skip the details and only provide a summary of the approach. We have the quadratic tangent majorant
\begin{equation}\label{eq:Q2}
\mathcal{Q}_{\rm 2}(\bu; \bu_{(k)}) :=\mathcal{F}(\bu) 
+\displaystyle{\frac{\lambda}{2}} \|\bM_{\rm 2}^{(k)}\bu\|^{2}_{2}+c_2,
\end{equation}}
where $c_2$ is a constant independent of $\bu$, and the matrix $\bM_{\rm 2}^{(k)}$ is defined as
\begin{equation}\label{eq:M2} \bM_{\rm 2}^{(k)} := \begin{bmatrix} \bW_{\rm 2}^{(k)} \\  & \bI  \end{bmatrix}\bD_1,
\end{equation} 
with $\bD_1$ as in \eqref{eq: D1}.

By letting $\bz_{(k)}=\bD_1\bu_{(k)}$, the weighting matrix $\bW_{\rm 2}^{(k)}$ is defined as 
\[\bW_{\rm 2}^{(k)} = \text{diag}\left((\bz_{(k)})^2+\epsilon^2\right)^{-1/4}.\]
%\MHP{Need to change notation for $\bw$. I think we need to use $(\bw_1)_{(k)}$ to denote the vector at iteration $k$ to be consistent with the notation above.}
As in (\ref{eq:MMQ1}), to solve the optimization problem~\eqref{eq:8}, we solve a sequence of reweighted least squares problems with the objective function $\mathcal{Q}_2$ defined in~\eqref{eq:Q2}. 

\subsection{Anisotropic space-time total variation (Aniso3DTV)}
To explain this approach, it is easier to directly consider the tensor notation. 
We define the tensor $\T{Y}$ in which the finite difference tensor is applied simultaneously across all three modes
 \begin{equation}\label{eq: tensorY}
\T{Y} = \T{U}\times_1 \bL_v\times_2\bL_h\times_3 \bL_t.
\end{equation}
We can write the 3D Anisotropic TV norm as the vectorized 1-norm of this tensor. That is
\[ \mathcal{R}_3(\bu) = \|\T{Y}\|_1 = \sum_{v=1}^{n_v}\sum_{h=1}^{n_h}\sum_{t=1}^{n_t} |\by_{v,h,t}|.\]
This is in contrast to $\mathcal{R}_1(\bu)$, which computes the sum of the tensor 1-norms in which only one derivative is applied per summand.

To derive an equivalent representation using matrix notation, consider the mode-$1$ unfolding of the tensor $\T{Y}$, $\bY_{(1)} = \bL_v \bU_{(1)}(\bL_t^T\otimes\bL^T_{h})$. Let $\by = \text{vec}(\bY_{(1)})$ and $\bu = \text{vec}(\bU_{(1)})$ denote the vectorizations of the mode-$1$ unfoldings of $\T{Y}$ and $\T{U}$ respectively, which are related through the formula
\[ \by = \bD_3 \bu \quad \mbox{ with }\quad \bD_3=  (\bL_t \otimes \bL_h \otimes \bL_v).\]
Therefore, we have $\mathcal{R}_3(\bu) := \|\bD_3\bu\|_1$.

\paragraph{Optimization problem and MM approach}{The problem that we want to solve can be formulated as
\begin{equation}\label{eq: Functional3}
\min_{\bu \in \R^{n}}\mathcal{J}_{3}(\bu):= \mathcal{F}(\bu)+\lambda \mathcal{R}_{3}(\bu), \nonumber
\end{equation}
which can be tackled with the MM approach similar to the one described in Section \ref{subsec: R1a}. Again, we consider the smoothed version $\mathcal{R}_{3\epsilon}(\bu)$ of $\mathcal{R}_3(\bu)$;  
the corresponding smoothed objective function is denoted by  $\mathcal{J}_{3\epsilon}(\bu)$. We majorize $\mathcal{J}_{3\epsilon}(\bu)$ by the quadratic tangent majorant
\begin{equation}\label{eq:Q3}
\mathcal{Q}_{\rm 3}(\bu; \bu_{(k)}) :=\mathcal{F}(\bu)
+\displaystyle{\frac{\lambda}{2}} \|\bM_{\rm 3}^{(k)}\bu\|^{2}_{2}+c_3, \nonumber
\end{equation}
 where $c_3$ is a constant independent of $\bu$ and
 \begin{equation}\label{eq:M3}
    \bM_3^{(k)} = \bW_3^{(k)} \bD_3.
\end{equation}
  The weighting matrix $\bW_3^{(k)}$ at iteration $k$ is defined as \[\bW_3^{(k)} = \text{diag}(((\bz_{(k)})^2 + \epsilon^2)^{-1/4}) \quad \mbox{ with } \quad \bz_{(k)}= \bD_3\bu_{(k)} .\] 
}

\subsection{3D space-time isotropic total variation (Iso3DTV)}\label{ssec:Iso3DTV}

In this next approach, we apply  isotropic total variation in all three directions, i.e., two spatial and one temporal direction. We first introduce the variables
\begin{equation}\label{eq:z3DisoTV}
%\begin{align*}
\begin{array}{ccl}
\bz_v^0(\bu):=&\>(\bI_{n_t}\otimes\bI_{n_h}\otimes\bL_{v}^0)\bu\,,\\
\bz_h^0(\bu):=&\>(\bI_{n_t}\otimes\bL_{h}^0\otimes\bI_{n_v})\bu\,,\\
\bz_t^0(\bu):=& \>(\bL_{t}^0\otimes\bI_{n_h}\otimes\bI_{n_v})\bu\,.
\end{array}
%\end{align*}
\end{equation}
%and
%\[
%\bZ=[\bz_v, \bz_h, \bz_t]\in\mathbb{R}^{n_vn_hn_t\times 3},
%\]
Recall that $\bL_d^0$, $d = v, h, t$ is obtained by augmenting $\bL_d$ with a row of zeros.
Then, we can compactly write 
%consider 
the following regularization term
\begin{equation}\label{eqn:R4}
%\mathcal{R}_4(\bu) =
%\sum_{k=1}^{n_vn_hn_t}\sqrt{
%((\bI_{n_t}\otimes\bI_{n_h}\otimes\bL_{v}^0)\bu)_k^2 + ((\bI_{n_t}\otimes\bL_{h}^0\otimes\bI_{n_v})\bu)_k^2 + ((\bL_{t}^0\otimes\bI_{n_h}\otimes\bI_{n_v})\bu)_k}^2
\begin{aligned}
\mathcal{R}_4(\bu) & := \>
\sum_{\ell=1}^{n_vn_hn_t}\sqrt{
(\bz_v^0(\bu))_{\ell}^2 + (\bz_h^0(\bu))_{\ell}^2 + (\bz_t^0(\bu))_{\ell}^2}\\ 
& = \> \|\,[\bz_v^0(\bu), \bz_h^0(\bu), \bz_t^0(\bu)]\,\|_{2,1}.
\end{aligned}
\end{equation}
To devise a tensor formulation for $\mathcal{R}_4(\bu)$, first  consider the following tensors 
\[
\T{Z}_v = \T{U}\times_1\bL_v^0,\quad
\T{Z}_h = \T{U}\times_2\bL_h^0,\quad
\T{Z}_t = \T{U}\times_3\bL_t^0,
\]
and their mode-3 unfoldings $(\bZ_v)_{(3)}$, $(\bZ_h)_{(3)}$, $(\bZ_t)_{(3)}$, respectively. Define a new tensor $\T{Y}\in\mathbb{R}^{n_s\times n_t\times 3}$ such that
\[
\T{Y}_{:,:,1} = (\bZ_v)_{(3)}^T,\quad
\T{Y}_{:,:,2} = (\bZ_h)_{(3)}^T,\quad
\T{Y}_{:,:,3} = (\bZ_t)_{(3)}^T.
\]
Then $\mathcal{R}_4(\bu)$ is the sum of the 2-norms of the mode-3 fibers of $\T{Y}$, that is
\[ 
\mathcal{R}_4(\bu) = \sum_{i=1}^{n_s} \sum_{j=1}^{n_t} \|\T{Y}_{i,j,:}\|_2\,.
\] 
To interpret this representation, the frontal slices of the tensor $\T{Y}$ are the collection of gradient images at all time instances, and the derivatives are taken one direction at a time. The regularization operator $\mathcal{R}_4(\bu)$ is the sum of two norms of its tubal fibers.

\paragraph{Optimization problem and MM approach} We have the following problem 
\begin{equation}\label{eq: Functional4}
\min_{\bu \in \R^n}\mathcal{J}_{4}(\bu):= \mathcal{F}(\bu)+\lambda \mathcal{R}_{4}(\bu).
\end{equation}
We first consider, instead of $\mathcal{R}_{4}(\bu)$, the smoothed regularization term
\[ \mathcal{R}_{4\epsilon}(\bu) := \>
\sum_{\ell=1}^{n_vn_hn_t}\sqrt{
(\bz_v^0(\bu))_{\ell}^2 + (\bz_h^0(\bu))_{\ell}^2 + (\bz_t^0(\by))_{\ell}^2 + \epsilon^2} \]
and the corresponding objective function $\mathcal{J}_{4\epsilon}(\bu)$. 
Following the derivation in \cite{wohlberg2007tv}, we devise weights to be used in an MM approach to Iso3DTV. We can define the quadratic tangent majorant $\mathcal{Q}_{\rm 4}(\bu; \bu_{(k)})$ for the objective function $\mathcal{J}_{4\epsilon}(\bu)$ as
\begin{equation}\label{eq:Q4}
\mathcal{Q}_{\rm 4}(\bu; \bu_{(k)}):=\mathcal{F}(\bu)
+\displaystyle{\frac{\lambda}{2}} \|\bM_{\rm 4}^{(k)}\bu\|^{2}_{2}+c_4,\nonumber
\end{equation}

where $c_4$ is a constant independent of $\bu$, and $\bM_4^{(k)}$ is the weighted matrix
\begin{equation}\label{eq:M4}
\bM_{4}^{(k)} := \bW_{\rm 4}^{(k)}\bD_4 \qquad \mbox{ with } \quad \bD_4 := \begin{bmatrix} \bI_{n_t}\otimes\bI_{n_h}\otimes\bL^0_{v} \\  
\bI_{n_t}\otimes \bL^0_{h}\otimes\bI_{n_v} \\
\bL_{t}^0\otimes\bI_{n_h}\otimes\bI_{n_v} 
\end{bmatrix},
\end{equation}  
and
\begin{equation}\label{eq:w4}
\bW_{4}^{(k)} = \bI_3\otimes \text{diag}\left(\left(
%(\bz_{(k)})^2+\epsilon^2)^{-1/4}
(\bz_v^0(\bu_{(k)}))^2+(\bz_h^0(\bu_{(k)}))^2+(\bz_t^0(\bu_{(k)}))^2+\epsilon^2
\right)^{-1/4}\right),\nonumber
\end{equation} 
where $(\bz_d^0(\bu_{(k)}))$ are the vectors $\bz_d^0$ in \eqref{eq:z3DisoTV}, $d=v,h,t$, evaluated at $\bu=\bu_{(k)}$, i.e., at the $k$th iteration. Finally, the matrix $\bD_4$ is similar to $\bD_1$ defined in \eqref{eq: D1}, with the augmented derivative matrices $\bL_d^0$ instead of $\bL_d$.

\subsection{Isotropic in space, anisotropic in time total variation (IsoTV)} 

This method can be considered as a variation of the AnisoTV method presented in Section \ref{subsec: R1a}, where only the spatial anisotropic total variation is replaced by spatial isotropic total variation. Namely, using the notation in \eqref{eq:z3DisoTV}, we consider the regularization term
\begin{equation}\label{eqn:R5}
\begin{aligned}
\mathcal{R}_{5}(\bu) &= \>  \sum_{\ell=1}^{n_vn_hn_t}\sqrt{
(\bz_v^0(\bu))_{\ell}^2 + (\bz_h^0(\bu))_{\ell}^2} +  \sum_{t=1}^{n_t-1}\|\bu^{(t+1)}-\bu^{(t)}\|_1\\
&= \> \|\,[\bz_v^0(\bu), \bz_h^0(\bu)]\,\|_{2,1} + \|(\bL_t\otimes \bI_{n_s})\bu\|_1.
\end{aligned} 
\end{equation}

The associated tensor formulation reads similar to the ones presented in subsections \ref{subsec: R1a} and \ref{ssec:Iso3DTV}, namely,
\[
\mathcal{R}_{5}(\bu)=\sum_{i=1}^{n_s} \sum_{j=1}^{n_t} \|\T{Y}_{i,j,:}\|_2 +\|\T{U}\times_3\bL_t\|_1,
\]
where $\T{Y}\in\mathbb{R}^{n_s\times n_t\times 2}$ is such that
\[
\T{Y}_{:,:,1} = (\bZ_v)_{(3)}^T,\quad
\T{Y}_{:,:,2} = (\bZ_h)_{(3)}^T.
\]

\paragraph{Optimization problem and MM approach} We have the following problem 
\begin{equation}\label{eq: Functional5}
\min_{\bu\in \R^n}\mathcal{J}_{5}(\bu)=\min_{\bu \in \R^n} \mathcal{F}(\bu)+\lambda \mathcal{R}_{5}(\bu).
\end{equation}
We define a smoothed version of $\mathcal{R}_5(\bu)$, denoted $\mathcal{R}_{5\epsilon}(\bu)$ where the smoothing is applied separately to the first and second terms; the corresponding smoothed objective function is denoted $\mathcal{J}_{5\epsilon}(\bu)$. We can then define the quadratic tangent majorant $\mathcal{Q}_{\rm 5}(\bu; \bu_{(k)})$ for the objective function $\mathcal{J}_{5\epsilon}(\bu)$ as
\begin{equation}\label{eq:Q5}
\mathcal{Q}_{\rm 5}(\bu; \bu_{(k)}):=\mathcal{F}(\bu)
+\displaystyle{\frac{\lambda}{2}} \|\bM_{\rm 5}^{(k)}\bu\|^{2}_{2}+c_5,\nonumber
\end{equation}

where $c_5$ is a constant independent of $\bu$, and $\bM_5^{(k)}$ is the weighted matrix
\begin{equation}\label{eq:M5}
\bM_{5}^{(k)} := \bW_{\rm 5}^{(k)}\bD_5 
\end{equation}  
with  
\begin{equation}\label{eq:w5}
\bD_5 := 
\begin{bmatrix} \bI_{n_t}\otimes\bI_{n_h}\otimes\bL^0_{v} \\  
\bI_{n_t}\otimes \bL^0_{h}\otimes\bI_{n_v} \\
\bL_{t}\otimes\bI_{n_h}\otimes\bI_{n_v} 
\end{bmatrix} \mbox{ and }\bW_{5}^{(k)} = 
\left[
\begin{array}{cc}
\bI_2\otimes \text{diag}\left(\bw_{(k)}^s\right) & \\
 & \text{diag}\left(\bw_{(k)}^t\right)
\end{array}
\right],
 \nonumber 
\end{equation}
where
\begin{equation}\bw_{(k)}^s=\left((
\bz_v^0(\bu_{(k)}))^2+(\bz_h^0(\bu_{(k)}))^2+\epsilon^2\right)^{-1/4} \mbox{ and }
\bw_{(k)}^t=\left(( \bz_t(\bu_{(k)}))^2+\epsilon^2\right)^{-1/4}. \nonumber\end{equation}
Here $\bz_d(\bu_{(k)})$ are again the vectors $\bz_d(\bu)$ in \eqref{eq:z3DisoTV}, $d=v,h$, evaluated at $\bu=\bu_{(k)}$, i.e., at the $k$th iteration.

\subsection{Group sparsity (GS)}\label{subsec: R2}

Group sparsity allows to promote sparsity when reconstructing a vector of unknown pixels that are naturally partitioned in subsets; see~\cite{bach2012optimization,herzog2012directional}. In our applications, there are several possible ways to define groups. For example, we can naturally group the variables corresponding to pixels at each time instant, i.e., $\{\bu^{(t)}\}_{t=1}^{n_t}$. To enforce piecewise constant structure in space and time, we adopt the following approach. Let $n_s' = (n_v-1)n_h + (n_h-1)n_v$ be the total number of pixels in the gradient images. 
Consider the groups defined by the vectors
\[ \bz_{\ell} = \begin{bmatrix} (\bL_s\bu^{(1)})_{\ell},\dots,(\bL_s\bu^{(n_t)})_{\ell} \end{bmatrix} = \left(\bI_{n_t} \otimes \be_\ell^T\bL_s\right)\bu \in \mathbb{R}^{n_t}, \qquad \ell=1,\dots,n_s'.\]
\vspace*{-4mm}
\begin{figure}[h!]
     \centering
     \includegraphics[scale=0.30]{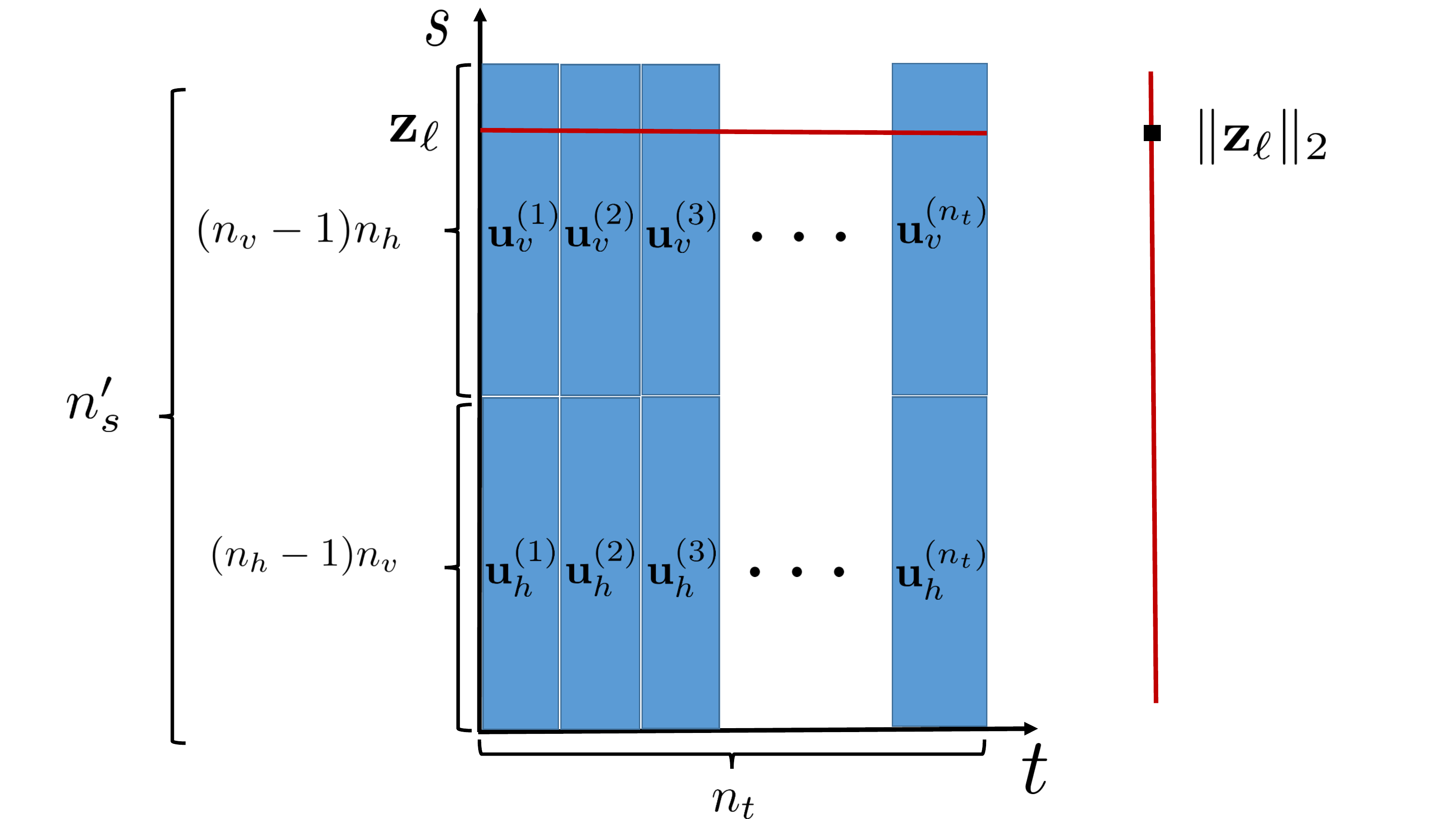}
     \caption{The vector of spatial derivatives $\bL_s\bu^{(t)}$ contains the partial derivatives with respect to the vertical ($\bu_v^{(t)}$)  and horizontal ($\bu_h^{(t)}$) directions for each image. These vectors are the columns of the matrix $\bZ$. We compute the 2-norm of each row $\bz_\ell$ of $\bZ$ and add them.}\label{fig:groupS}
\end{figure}

Alternatively, define the matrix $\bZ$ whose columns represents the vectorized gradient images at different time $t$. 
\[
\begin{aligned}
\bZ &=[\bL_s\bu^{(1)},\dots,\bL_s\bu^{(n_t)}]=\bL_s\bU\in\R^{n_s'\times v_t},\\
\bz &={\rm{vec}}(\bZ)=(\bI_{n_t}\otimes \bL_s)\bu\,.
\end{aligned}
\]
Note that $\bz_\ell$ are the rows of $\bZ$. These are also illustrated in Figure \ref{fig:groupS}. 
The regularization term corresponding to group sparsity can then be expressed as a mixture of norms 
\begin{equation}\label{eq: SparsityReg} 
\mathcal{R}_6(\bu) := \sum_{\ell=1}^{n_s'}\|\bz_{\ell}\|_2 = \sum_{\ell=1}^{n_s'} \left( \sum_{t=1}^{n_t} (\bL_s\bu^{(t)})_{\ell}^2 \right)^{1/2} = \|\bL_s\bU\|_{2,1}. \nonumber 
\end{equation}
This was already used in Equations \eqref{eq: isoTV} and \eqref{eqn:R4}. 
In other words, the regularization term behaves like a 1-norm on the vector $$\begin{bmatrix}\|\bz_1\|_2& \dots& \|\bz_{n_s'}\|_2\end{bmatrix}.$$ By inducing sparsity on the vector of 2-norms of $\bz_{\ell}$, $\ell=1,\dots,n_s'$, we encourage entries of the vector $\|\bz_{\ell}\|_2$ (and, in turn, each vector $\bz_{\ell}$) to be zero. On one hand, using this regularization, we are ensuring that the sparsity in the gradient images is being shared across time instances. On the other hand, this regularization formulation does not enforce sparsity across the groups, i.e., across vectors $\bz_{\ell}$. 

Let $\T{U}$ be the tensor of images, and let $\T{X} = \T{U}\times_1\bL_v$ and $\T{Y} = \T{U}\times_2\bL_h$ be the tensors obtained by taking the gradient in the vertical and horizontal directions. Then $\mathcal{R}_6(\bu)$ is the sum of 2-norm of the mode-3 fibers of $\T{X}$ and $\T{Y}$. That is
\[ \mathcal{R}_6(\bu) = \sum_{i=1}^{(n_v-1)} \sum_{j=1}^{n_h} \|\T{X}_{i,j,:}\|_2 + \sum_{i=1}^{(n_h-1)}\sum_{j=1}^{n_v}\|\T{Y}_{i,j,:}\|_2. \] 
Note also that $\bZ = \bmat{\bX_{(3)}, & \bY_{(3)}}^T$. 
\paragraph{Optimization problem and MM approach} 
Corresponding to the regularization operator $\mathcal{R}_6$, we can define the optimization problem:
\begin{equation}\label{eq:J6}
\min_{\bu \in \R^{n}}\mathcal{J}_{6}(\bu) :=  \mathcal{F}(\bu)+{\lambda} \mathcal{R}_{6}(\bu), \nonumber
\end{equation}
where $\lambda > 0$. We can apply the MM approach similar to Section~\ref{subsec: R1a}.
We now seek a quadratic tangent majorant for a smoothed version of $\mathcal{R}_{6}(\bu)$. To this end, let $\bu_{(k)}$ be the current iterate (similarly, define $\bz_{(k)} = (\bI_{n_t}\otimes \bL_s)\bu_{(k)}$). Then, we have that
\[
\begin{aligned}
\mathcal{R}_6(\bu) \leq & \>   {\sum_{\ell=1}^{n_s'} \sqrt{\|\bz_{\ell}\|_2^2 + \epsilon^2}= : \mathcal{R}_{6\epsilon}(\bu)},  \\
\mathcal{R}_{6\epsilon}(\bu) \leq & \> \sum_{\ell=1}^{n_s'}\frac{\|\bz_{\ell}\|_2^2}{2\sqrt{\|(\bI_{n_t}\otimes \be_{\ell}^T\bL_s)\bu_{(k)}\|_2^2+\epsilon^2}} + \tilde{c}_6, 
\end{aligned}
\]
where $\tilde{c}_6$ is a constant independent of $\bz_{\ell}$ and $\bu$. The corresponding smoothed optimization function is defined as $\mathcal{J}_{6\epsilon}(\bu)$. Let us define the weighting matrix $\bW_6^{(k)}$ of size $n_s'\times n_s'$ as 
\[ \bW_6^{(k)} := \text{diag}\left(\frac{1}{\sqrt{\| (\bI_{n_t}\otimes \be_1^T\bL_s)\bu_{(k)}\|_2^2+\epsilon^2}}, \dots,  \frac{1}{\sqrt{\| (\bI_{n_t}\otimes \be_{n_s'}^T\bL_s)\bu_{(k)}\|_2^2+\epsilon^2}}\right)^{1/2} .\]
We can use this weighting matrix to define the quadratic tangent majorant
\[ \mathcal{Q}_6(\bu;\bu_{(k)}) := \mathcal{F}(\bu)  + \frac{\lambda}{2} \| \bM_6^{(k)} \bu\|_2^2 + c_6,  \]
where $c_6 = \lambda \tilde{c}_6$ and the matrix $\bM_6^{(k)}$ takes the form
\begin{equation}\label{eq:M6}
    \bM_6^{(k)}  := (\bI_{n_t} \otimes \bW_6^{(k)}) \bD_6 \quad \mbox{ with } \quad \bD_6 :=  (\bI_{n_t} \otimes \bL_s).
\end{equation}

\subsection{Summary of proposed approaches} In this section, we have presented six different regularization terms for promoting edge-preserving reconstructions in dynamic inverse problems. Here we show that they can be treated in a unified fashion, thereby providing a succinct summary of all the proposed methods. For each regularization term, we solve an optimization problem of the form
\begin{equation}\label{eqn:dynamic}
    \min_{\bu \in \mathbb{R}^{n}} \mathcal{J}_{j\epsilon}(\bu) := \mathcal{F}(\bu) + \lambda \mathcal{R}_{j\epsilon}(\bu), \qquad j = 1,\dots,6,
\end{equation}
where $\mathcal{R}_{j\epsilon}(\bu)$ is a smoothed regularization term depending on the method used, and $\mathcal{F}(\bu)$ is a term that measures the data-misfit. For each optimization problem, we have derived an MM approach that solves a sequence of iteratively reweighted least squares problem. That is, at step $k$ given an initial guess $\bu_{(0)}$, we solve the sequence of optimization problems
\begin{equation}\label{eqn:genMM}  
\bu_{(k+1)} =  \argmin_{\bu \in \mathbb{R}^{n}}\displaystyle{\frac{1}{2}} \|\bF\bu-\bd\|_{\bGamma^{-1}}^2
+\displaystyle{\frac{\lambda}{2}} \|\bM_j^{(k)}\bu\|^{2}_{2} \qquad k=0,1\dots\,.
\end{equation}
The matrix $\bM_j^{(k)}$ takes different forms depending on the regularization technique used. 

Table~\ref{Table: Summary} summarizes some details about the proposed regularization terms, and points to the formulas defining the reweighting matrices appearing within $\bM_j^{(k)}$ in the MM step. In Section~\ref{sec: iterativemethods}, we discuss iterative methods to efficiently solve the sequence of least squares problems (\ref{eqn:genMM}) and select the regularization parameter $\lambda$. 
\begin{table}[h!]

\begin{center}
\caption{The proposed methods, together with their associated regularization terms and the weighting matrices for the MM step. The vectors $\bz_d^0(\bu)$, $d=v,h,t$ are defined in (\ref{eq:z3DisoTV}).} 
	\label{Table: Summary}
\begin{tabular}{ c| c| c }
\textbf{Method} & $\mathcal{R}_i(\bu)$ & \textbf{MM weights} \\ \hline
 AnisoTV & $\left\|(\bI_{n_t} \otimes  \bL_s)\bu\right\|_1 + \|(\bL_t\otimes \bI_{n_s})\bu\|_1$ & \eqref{eq:M1} \\ \hline 
TVplusTikhonov & $ \>  \|(\bI_{n_t}\otimes\bL_s)\bu\|_1 + \|(\bL_t\otimes \bI_{n_s})\bu\|_2^2$ & \eqref{eq:M2} 
\\ \hline
 Aniso3DTV & $\|(\bL_t \otimes \bL_h \otimes \bL_v)\bu\|_1$ & \eqref{eq:M3} \\ \hline
 Iso3DTV & $\|\,[\bz_v^0(\bu), \bz_h^0(\bu), \bz_t^0(\bu)]\,\|_{2,1}$ & \eqref{eq:M4}
  \\\hline
 IsoTV & $\|\,[\bz_v^0(\bu), \bz_h^0(\bu)]\,\|_{2,1} + \|(\bL_t\otimes \bI_{n_s})\bu\|_1$ & \eqref{eq:M5} \\\hline
 GS & $\|\bL_s\bU\|_{2,1}$ & \eqref{eq:M6}  
\end{tabular}
\end{center}
\end{table}

\section{Alternative approaches}\label{sec: alternative}
In Section \ref{sec: newmethods} we presented a variety of regularization methods that use different forms of TV regularization to obtain solutions methods that enhance edge representation. In this section, we summarize some alternative approaches that can be used, still within the MM framework, to enforce edge-preserving reconstructions.

\paragraph{Beyond the $\ell_1$ and $\ell_2$ norms}
One way to interpret the anisotropic TV is that it enforces sparsity in the gradient images. A natural measure of sparsity of a vector is the $\ell_0$-``norm'', which counts the number of nonzero entries. However, solving minimization problems that involve the $\ell_0$ term is known to be NP-hard, hence to remedy this difficulty one approximates the $\ell_0$-``norm'' by $\ell_1$ convex relaxation. Several nonconvex penalties with $0<q<1$ have been used alternatively to $\ell_1$; see \cite{chartrand2007exact, krishnan2009fast, xu2012l}. 
 
The $\ell_0$``norm'' of the gradient can be regarded as the length of the partition boundaries as in Potts model \cite{potts1952some} or piecewise constant Mumford–Shah model \cite{mumford1985boundary}.  
With this analogy, the methods that we discuss in Section \ref{sec: newmethods} can be generalized by using $\ell_q$ regularization. For example, the regularization term (\ref{eqn:R1}) in Section~\ref{subsec: R1a} can be generalized by choosing 
\[ \mathcal{R}_{1}^q(\bu) = \frac{1}{q}\|\bD_1\bu\|_q^q, \qquad 0<q\leq 2. \] 
Similarly, the GS method presented in Section \ref{subsec: R2} can be espressed using general mixed $p$-$q$ ``norms'' instead of 2-1 ``norms''. 
\paragraph{Beyond the gradient operator}
One can build appropriate sparsity transforms using and combining operators other than the first order finite difference operator $\bL_d$ defined in (\ref{D2}), where $d= v,h,t$. A first simple extension replaces $\bL_d$ by a discretization of the \emph{second order derivative operator}, which can still assist in preserving edges \cite{acton2009diffusion, pal2015brief, tan2018digital}. 
Moreover, one can incorporate a \emph{wavelet} transform: since being introduced in the 1910s by Haar \cite{stankovic2003haar} as a family of piecewise constant functions from which one can generate orthonormal bases for the square integrable functions, an extensive literature on wavelets was lead by discoveries from Str\"{o}mberg, Meyer, and culminated with the celebrated Daubechies wavelets of compact support; see for instance \cite{daubechies1992ten, dgm1986painless, meyer1992wavelets, stromberg1980weights} and references therein for more details.
Similar to wavelets, \emph{framelet} representations of images are orthogonal basis transformations that form a dictionary of minimum size that initially decomposes the images into transformed coefficients.
A variety of framelets can be used; for instance, 
one can use tight frames as in \cite{buccini2,COS09};
they are determined by linear B-splines that are formed by a low-pass filter and two high-pass filters that define the corresponding masks of the analysis operator. Finally, a number of variations are also possible when specifically considering the GS regularizer proposed in Section \ref{subsec: R2}. For instance, one can consider `overlapping groups' and also replace $\bL_s$ by other operators, such as the ones mentioned above. 
%the second order derivative, wavelets and framelets.  
It is well-known that, beyond dynamic inverse problems, sparse representations can improve pattern recognition, feature extraction, compression, multi-task regression and noise reduction; see, for example, \cite{ashwini2017sparse, bach2008exploring, chambolle2005total, kim2010tree, patel2013sparse, stephane1999wavelet}. 
\paragraph{Beyond one single regularization parameter}
Specifically for dynamic inverse \linebreak[4]problems, it may be meaningful to adapt the  regularization parameters based on the dynamics. For instance, one can define dedicated regularization parameters for different channels or domains (spatial or temporal). Within the framework presented in Section \ref{sec: newmethods}, this can be achieved by setting, in addition or as an alternative to $\lambda$, appropriate valued for the parameters $\alpha_d$ in (\ref{D2}). For instance, \cite{4193460} considers a  scenario where the regularization parameters are different for the spatial and the temporal domains.

\section{Iterative methods and parameter selection techniques}\label{sec: iterativemethods}
In this section we describe a numerical method to solve the optimization problems arising from the approaches described in Section \ref{sec: newmethods}, focusing on the minimization step. 
%\MHP{Should we leave the minimization step here or move it in the background?}
Considering problem \eqref{eqn:dynamic} and still assuming, for now, a fixed regularization parameter $\lambda > 0$, to determine $\bu_{(k+1)}$ as in \eqref{eqn:genMM} 
%\ }
we solve for the zero gradient of %$\mathcal{J}_{j\epsilon}(\bu)$,
%\MP{Here we solve for the zero gradient of $\mathcal{Q}_j(\bu,\bu_{(k)})$}
$\mathcal{Q}_j(\bu;\bu_{(k)})$
which leads to the regularized normal equations
(or general Tikhonov problem),
\begin{equation}\label{eq: normaleqQuadMajorant}
\left(
  \bF^T \bGamma^{-1}\bF + \lambda(\bM_j^{(k)})^{T}\bM_j^{(k)}
\right)\bu_{(k+1)} = \bF^T\bGamma^{-1}\bd.
\end{equation}
The system \eqref{eq: normaleqQuadMajorant} has a unique solution, if
\begin{equation}\label{eq: NANL}
\mathcal{N}(\bF^T \bGamma^{-1}\bF)\cap \mathcal{N}((\bM_j^{(k)})^T \bM_j^{(k)})=\{\mathbf{0}\}.
\end{equation}

This condition is equivalent to 
\begin{equation}\label{eq:nullspcondition}
    \mathcal{N}(\bF) \cap \mathcal{N}(\bD_j) = \{\mathbf{0}\} \qquad j=1,\dots,6,
\end{equation}
where we have used the properties that $\mathcal{N}(\bA^T\bA) = \mathcal{N}(\bA)$ and $\mathcal{N}(\bA\bB) = \mathcal{N}(\bB)$ if $\bA$ has full column rank, and 
where the matrices $\bD_j$ were defined in~\eqref{eq:M1},~\eqref{eq:M3},~\eqref{eq:M4}, \eqref{eq:M5}, and \eqref{eq:M6} (for convenience, we have defined $\bD_2 = \bD_1$). If \eqref{eq:nullspcondition} is satisfied, the solution to
\eqref{eq: normaleqQuadMajorant} is the unique
minimizer of $\mathcal{Q}_j(\bu; \bu_{(k)})$. Therefore, for methods 1-3 (AnisoTV, TVplusTikhonov, and Aniso3DTV), this fits the assumptions of~\cite[Theorem 5]{huang2017majorization}, and as a consequence the sequence $\{\bu_{(k)}\}$ converges to a stationary point of $\mathcal{J}_{j\epsilon}(\bu)$ for each method. For methods 4-6 (Iso3DTV, IsoTV, and GS), it may be possible to extend the analysis from that paper; however, we do not pursue it here. 

Since solving \eqref{eq: normaleqQuadMajorant} for large-scale matrices $\bF$ and $\bM^{(k)}$ may be computationally demanding or even prohibitive, we project \eqref{eq: normaleqQuadMajorant} unto a low dimensional subspace (namely, a generalized Krylov space) and solve a much smaller projected problem. If the approximate solution is not satisfactory, we extend the search space with the (normalized) residual. This leads to the generalized Golub-Kahan (GKS) process \cite{lampe2012large}; our description follows \cite{huang2017majorization}, adapted for the problems presented here. 

Given a $d$-dimensional ($d \ll n$) search space $\mathcal{V}_{d} = \mathcal{R}(\bV_{d})$, with $\bV_d \in \mathbb{R}^{n \times d}$ where, for numerical reasons, we keep the columns of $\bV_{d}$ orthonormal, we compute an approximate solution to \eqref{eq: normaleqQuadMajorant} as follows. Given the thin QR-decompositions
$\bGamma^{- 1/2}\bF\bV_{d} = \bQ_{\bF}^{(k)}\bR_{\bF}^{(k)}$ and 
$\bM_j^{(k)}\bV_{d} = \bQ^{(k)}_{\bM}\bR^{(k)}_{\bM}$, substituting  $\bu = \bV_d \by$ in \eqref{eqn:genMM} (or, equivalently, \eqref{eq: normaleqQuadMajorant}) 
leads to the small minimization problem
\begin{equation} \label{eq: SmallQuadraticMajorant1}
\by_{(k+1)} = \arg \min_{\by \in \R^d}
\left\| 
\left[ \begin{array}{c}
  \bR_{\bF}^{(k)} \\ \lambda^{1/2}\bR_{\bM}^{(k)}
\end{array} \right] \by - 
\left[ \begin{array}{c}
     (\bQ_{\bF}^{(k)})^T \bd \\ \bzero
\end{array} \right]
\right\|_2^2 ,
\end{equation}
and the corresponding $d \times d$ regularized normal equations for 
$\by_{(k+1)}$,
\begin{equation}\label{eq: SmallnormaleqQuadMajorant}
\left(
  (\bR_{\bF}^{(k)})^T\bR_{\bF}^{(k)} + \lambda(\bR_{\bM}^{(k)})^{T}\bR_{\bM}^{(k)}
\right)\by_{(k+1)} = \bR_{\bF}^T \bQ_{\bF}^T\bd .
\end{equation}
At each iteration, we use a GCV-like condition to determine the regularization parameter $\lambda$ (see below), after which the system 
\eqref{eq: SmallnormaleqQuadMajorant} can be solved at low cost. 
This gives the approximate solution $\bu_{(k+1)} = \bV_d \by_{(k+1)}$. The residual for 
\eqref{eq: normaleqQuadMajorant} can be computed
as
\begin{equation} \label{eq: ResidualNormalEqs}
  \br^{(k+1)} = 
  \bF^T\bGamma^{-1}(\bF\bV_d \by_{(k+1)} - \bd)
  + \lambda(\bM^{(k)})^T\bM^{(k)}\bV_d \by_{(k+1)} .
\end{equation}
The iteration is halted if the relative change in the solution drops below a given tolerance. 
Otherwise, we use the normalized residual to expand the search space, 
$\bV_{d+1} = [\bV_d \;\; \br^{(k+1)}/ \|\br^{(k+1)}\|_2 ]$. While in exact arithmetic 
$\br^{(k+1)} \perp \bV_d$, in practice, for numerical stability, we first explicitly orthogonalize the new residual against $\bV_d$. 
Next, we compute $\bW_j^{(k+1)}$ and $\bM_j^{(k+1)}$ as discussed in Section \ref{sec: newmethods} (for each method) and continue the iteration, solving for $\bu_{(k+2)}$. As $\bGamma$ is fixed, the thin QR-decomposition $\bGamma^{- 1/2}\bF\bV_{d} = \bQ_{\bF}^{(k)}\bR_{\bF}^{(k)}$
can be updated efficiently for the new column. To compute a small initial search space, the GKS algorithm is generally started by a few steps, say $\ell$, of the Golub-Kahan bidiagonalization for $\bF$ and $\bF^T \bd$. So, for $k=0$ we have $d = \ell$ and more generally, at step $k$, $d = k+\ell$. We emphasise again that an approximation of $\bu_{(k+1)}$ in \eqref{eqn:dynamic} is obtained by solving a single projected problem of dimension $k+\ell$, and that this hold for every $k$.

We now discuss briefly the choice of the regularization parameter $\lambda > 0$,
which balances the misfit term and the regularization term. 
Different techniques can be used to 
determine the regularization parameter, such as the L-curve, the discrepancy principle (DP), the unbiased predictive risk estimator (UPRE), and generalized cross validation (GCV) \cite{vogel2002computational, engl1996regularization, hanke1993regularization, hansen1998rank}. 
Here, at each iteration, we 
determine the regularization parameter by generalized cross validation applied to the 
projected problem  
\eqref{eq: SmallQuadraticMajorant1}
or \eqref{eq: SmallnormaleqQuadMajorant}.
To compute the GCV functional we use the generalized singular value decomposition of 
$[(\bR_{\bF}^{(k)})^T \; (\bR_{\bM}^{(k)})^T]^T$, 
which can be computed efficiently as 
$\bR_{\bF}^{(k)}, \bR_{\bM}^{(k)} \in \bbR^{d \times d}$ ($d = \ell + k$).
For further details, see \cite{buccini2021generalized}.
\section{Numerical experiments}\label{sec: nresults}
In this section we provide numerical examples from three different dynamic inverse problems: image deblurring, x-ray CT with simulated data, and x-ray CT with real data (for which the true solution is not available). Our goal is two-fold: to show that using dynamic information can be advantageous and to compare the different methods that we propose in this paper. 
\paragraph{Discussion on selecting the numerical examples} The first example that we consider concerns a synthetic space-time image deblurring where images change in time, but the blurring operator is fixed for all the time instances. Even though this is not an inverse problem with limited measurements,  we use it as an example to compare all the proposed methods since the true solution is available.
The second example is a problem from dynamic photoacoustic tomography (PAT), where we have a small number of measurements per time step (since information is collected from limited angles), but we have a lot of time steps yielding a large number of measurements overall.  This is our largest test problem in which we have over $1.9$ million unknowns, and the forward operator $\bA^{(t)}$ changes at each time step. In this example, we compare a few of the proposed methods for dynamic inverse problems with the results from the static inverse problem. The last example concerns real data arising from limited angle CT where the target of interest is a sequence of ``emoji images''. For this example, the true solution is not available and we can only provide qualitative assessment, but this example clearly illustrates the impact of incorporating temporal information in the reconstruction process; we also illustrate the effect of heuristically incorporating nonnegativity constraints. 

\paragraph{Quality measures and stopping criteria} 
To assess the quality of the reconstructed solution, we compute the Relative Reconstruction Errors (RREs) which are obtained using the $\ell_2$ error norms. That is, for some recovered $\bu_{(k)}$ at the $k$-th iteration, the RRE is defined as follows
\begin{equation} \label{eq: RRE} 
{\rm RRE}(\bu_{(k)},\bu_\mathrm{true}) = \frac{||\bu_{(k)}-\bu_\mathrm{true}||_2}{||\bu_\mathrm{true}||_2}. \nonumber
\end{equation}
In addition to RRE, in one of the examples, we use the Structural SIMilarity index (SSIM) to measure the quality of the computed approximate solutions. The definition of the SSIM is involved and we refer to \cite{SSIM} for details. Here we just recall that the SSIM measures how well the overall structure of the image is recovered; the higher the index, the better the reconstruction. The highest value achievable is $1$. 

We stop the iterations when the maximum number of 150 iterations is reached or if the discrepancy principle (DP) is satisfied, that is, when the following condition is satisfied
\begin{equation}\label{eq: DP}
\|\bF\bu-\bd\|_{\bGamma^{-1}}\leq\eta\delta,
\end{equation}
where $\eta>1$ is a user defined constant (we pick $\eta$ to be $1.01$) and $\delta$ is an estimate of the noise magnitude ($\|\be\|_{\bGamma^{-1}}$).
For fair comparison, in all the numerical examples, we set the smoothing parameter $\epsilon = 10^{-3}$ and $\ell = 5$, that is, we run 5 iterations of the Golub-Kahan bidiagonalization algorithm to generate an initial subspace. 
In all the examples we perturb the measurements with white Gaussian noise that is obtained when the entries of the vector $\be$ in the data vector $\bd$ are uncorrelated realizations of a Gaussian random variable with $0$ mean. In this case we refer to the ratio 
$\sigma=\|\be\|_{\bGamma^{-1}}/\|\bF\bu\|_{\bGamma^{-1}}$ as the noise level. 
\subsection{Example 1: Space-time image deblurring}
\begin{figure}[!ht]
\centering
    \begin{minipage}{0.45\textwidth}
		\includegraphics[width=\textwidth]{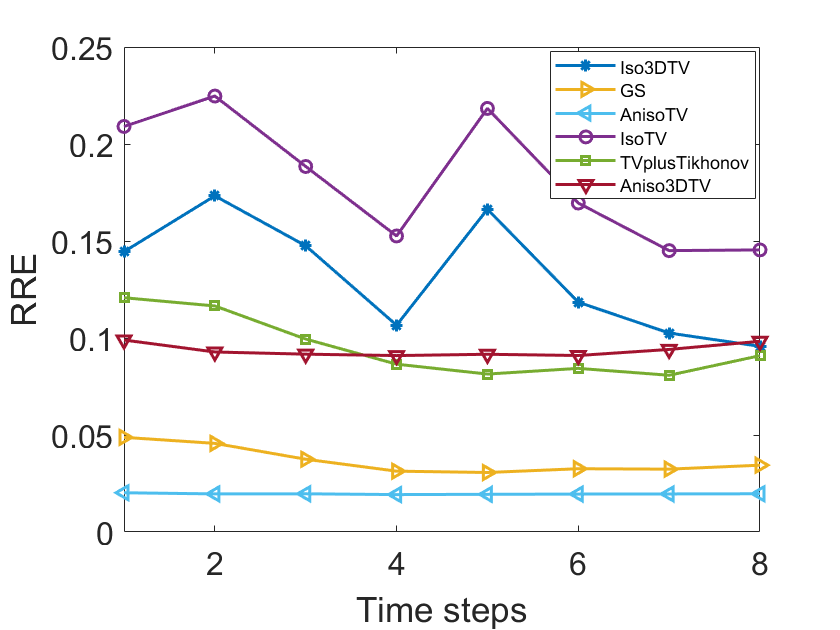}\\(a)
	\end{minipage}
	 \begin{minipage}{0.45\textwidth}
		\includegraphics[width=\textwidth]{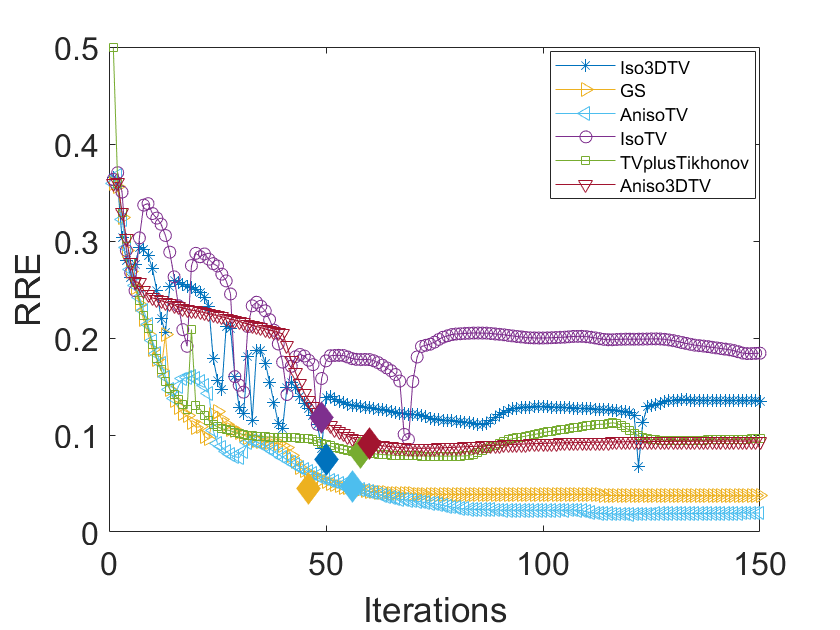}\\(b)
	\end{minipage}
	\caption{Space-time image deblurring test problem: a) RRE computed at the iteration when the DP \eqref{eq: DP} is first satisfied,  for each time step. b) Convergence history of the RRE (all of the time steps together). The methods considered here are AnisoTV, TVplusTikhonov, IsoTV, Aniso3DTV, Iso3DTV, and GS. The solid diamond markers highlight the iteration satisfying the DP.}
   \label{Fig: ErrorSL} 
\end{figure}

\begin{figure}[ht!]
\begin{center}
    \begin{minipage}{0.18\textwidth}
		\includegraphics[width=\textwidth]{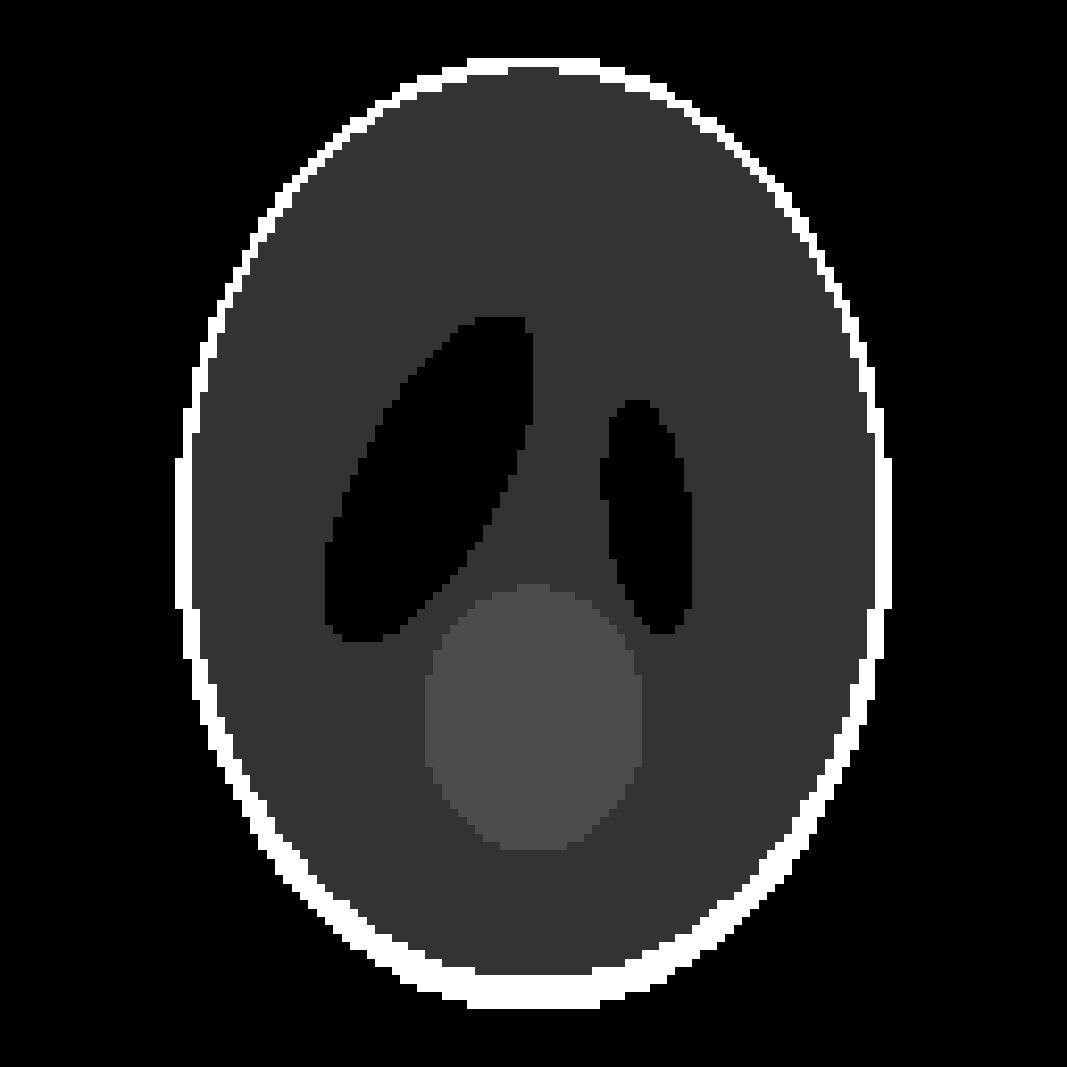}\\%(a)
	\end{minipage}
	\begin{minipage}{0.18\textwidth}
		\includegraphics[width=\textwidth, height = \textwidth]{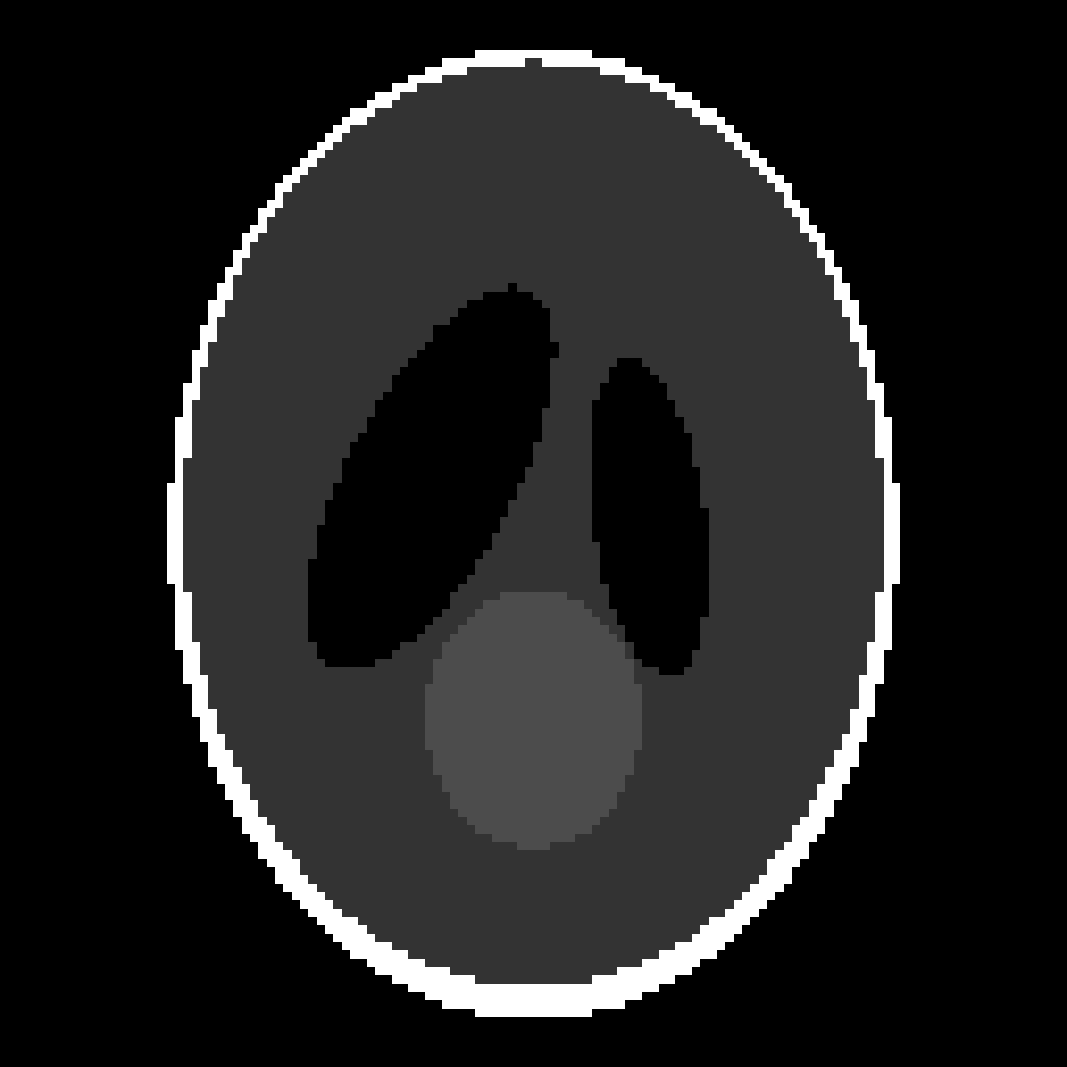}\\%(b)
	\end{minipage}
	\begin{minipage}{0.18\textwidth}
		\includegraphics[width=\textwidth]{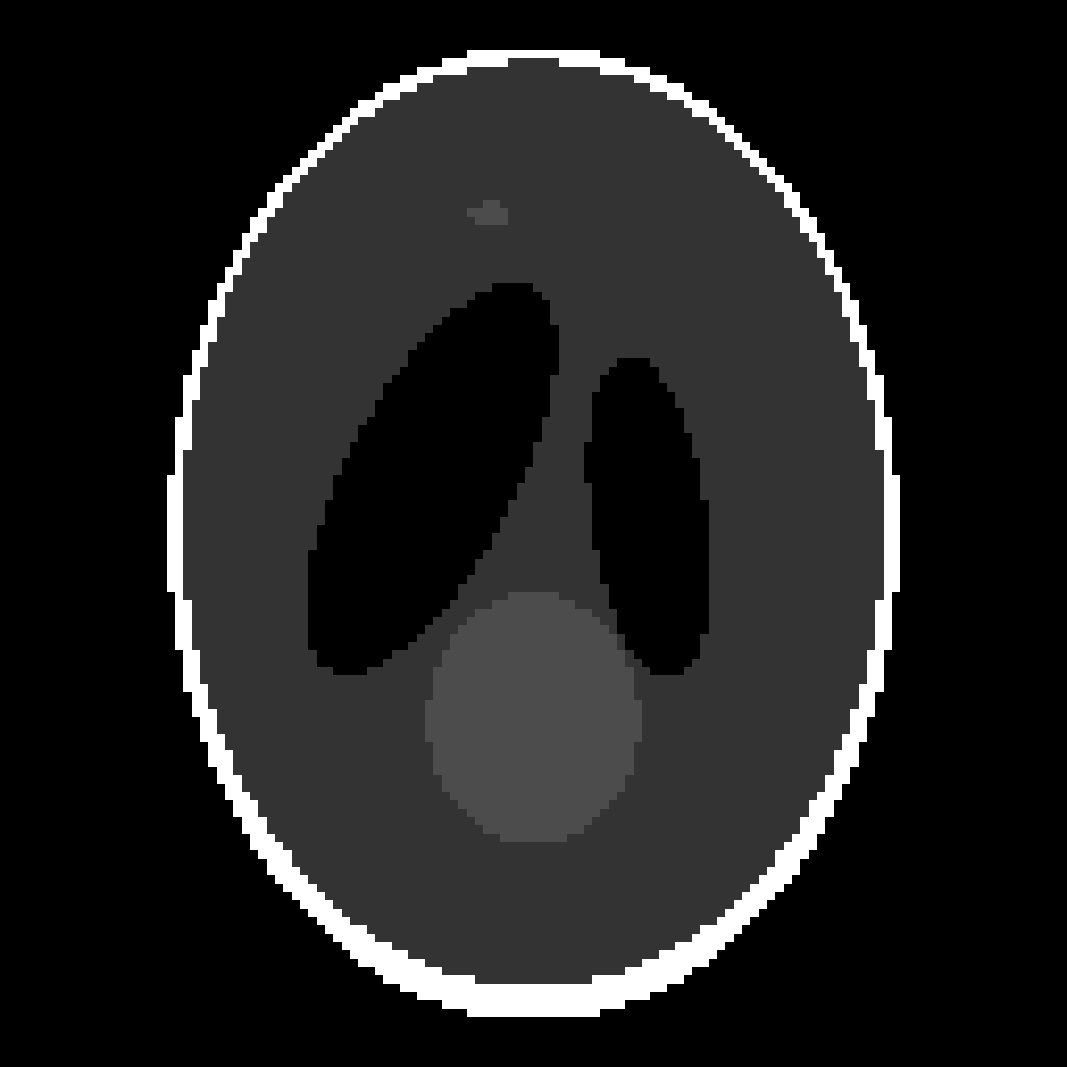}\\%(c)
	\end{minipage}
	\begin{minipage}{0.18\textwidth}
		\includegraphics[width=\textwidth]{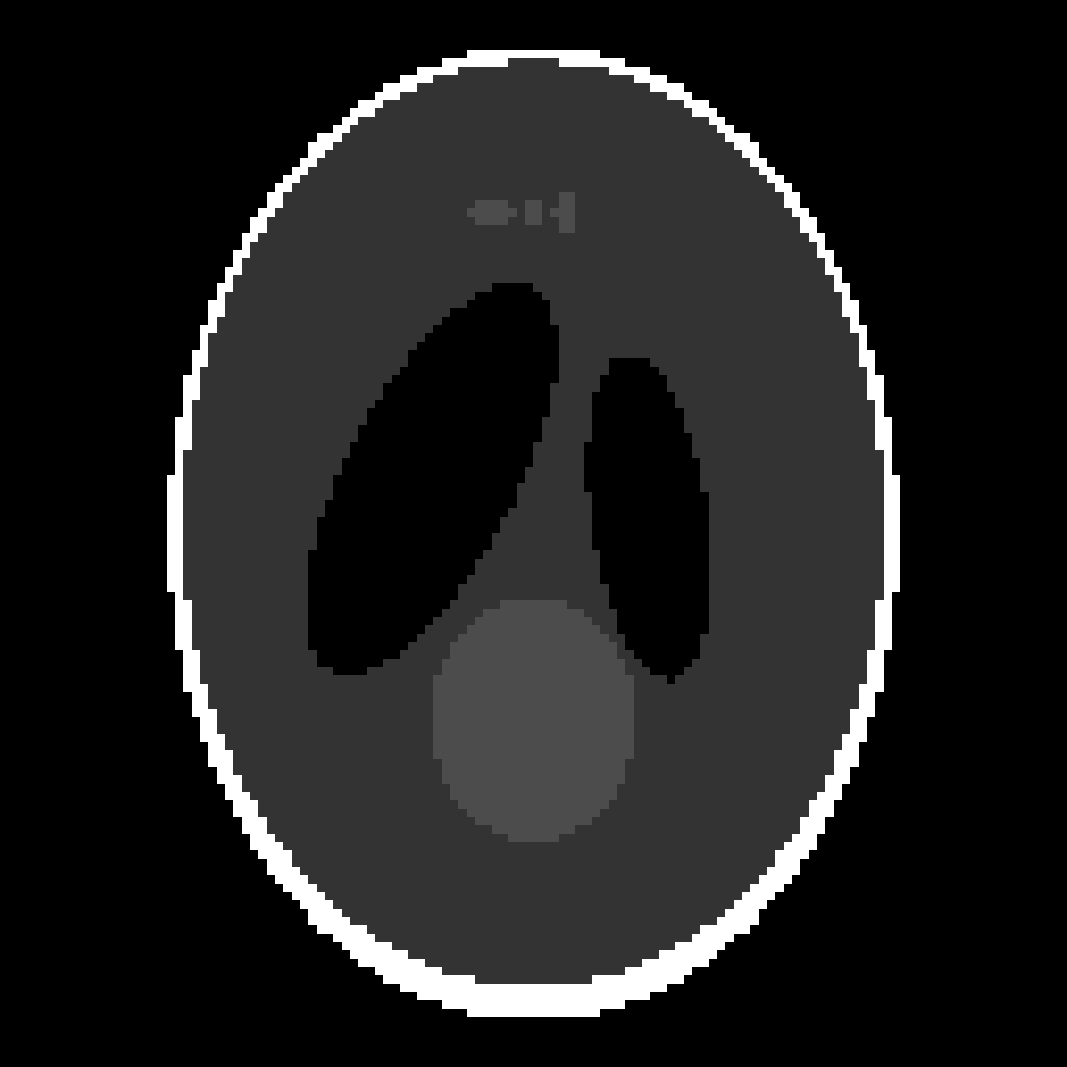}\\%(d)
	\end{minipage}
	\begin{minipage}{0.18\textwidth}
		\includegraphics[width=\textwidth]{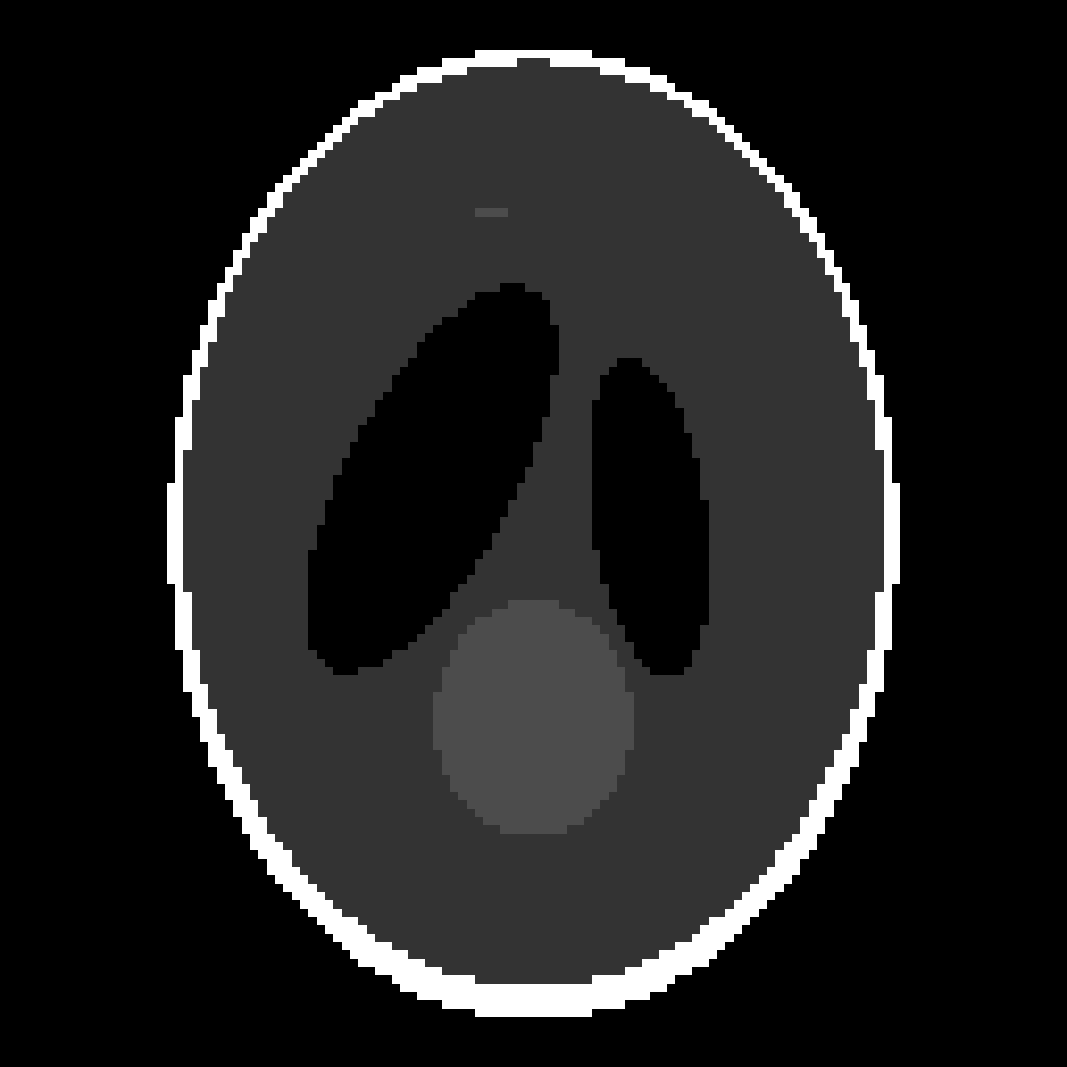}\\%(e)
	\end{minipage}

	\begin{minipage}{0.18\textwidth}
		\includegraphics[width=\textwidth]{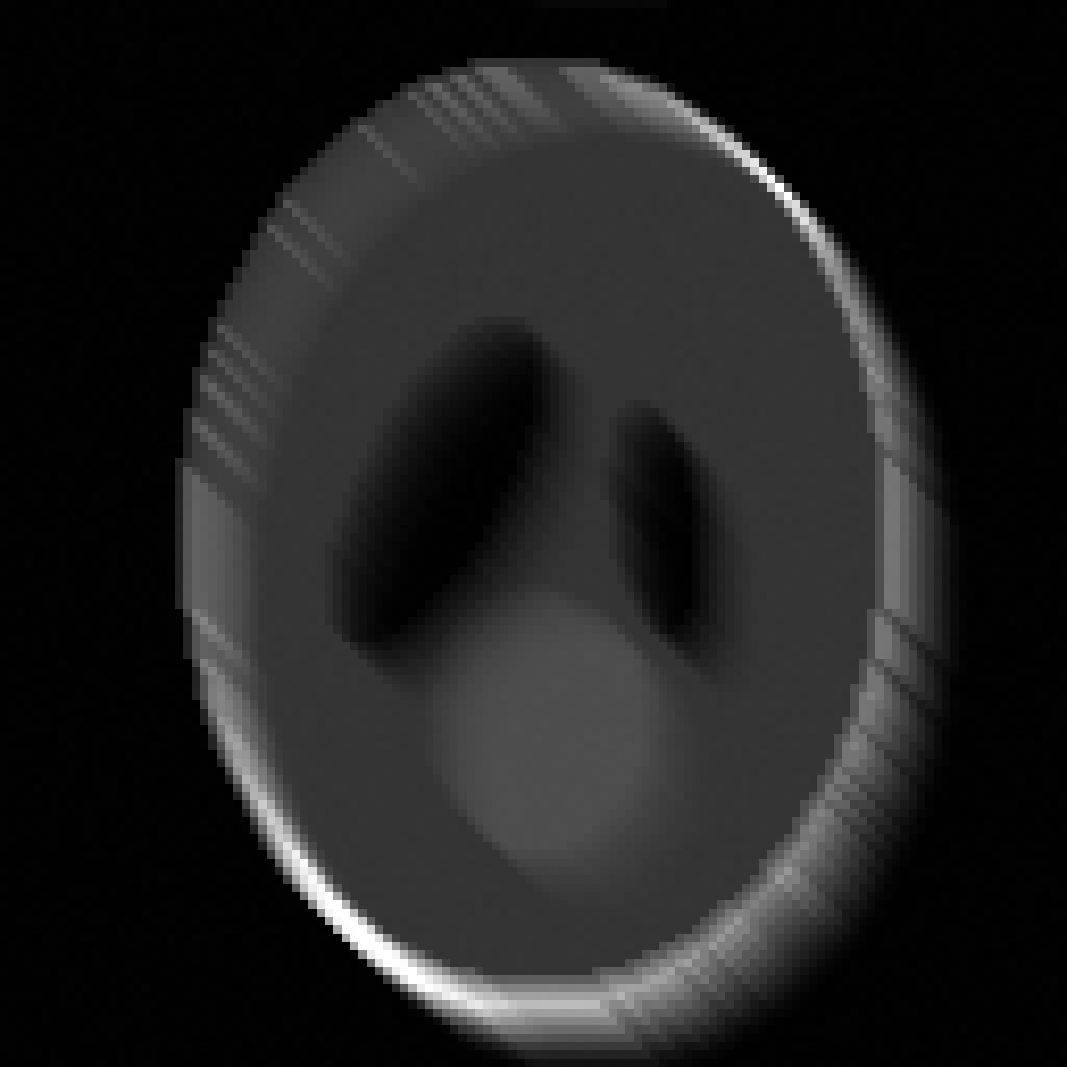}\\%(i)
	\end{minipage}
	\begin{minipage}{0.18\textwidth}
		\includegraphics[width=\textwidth]{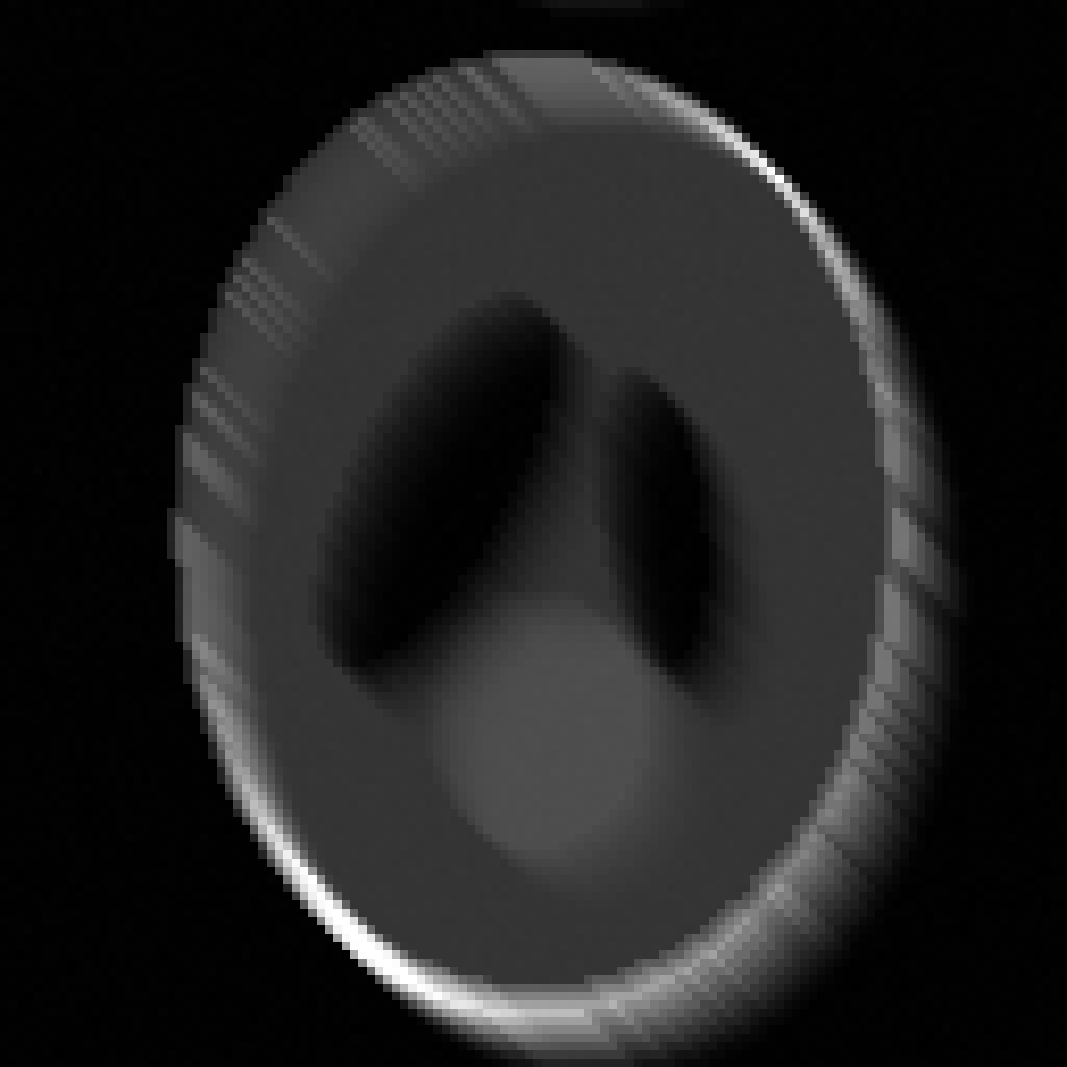}\\%(j)
	\end{minipage}
		\begin{minipage}{0.18\textwidth}
		\includegraphics[width=\textwidth]{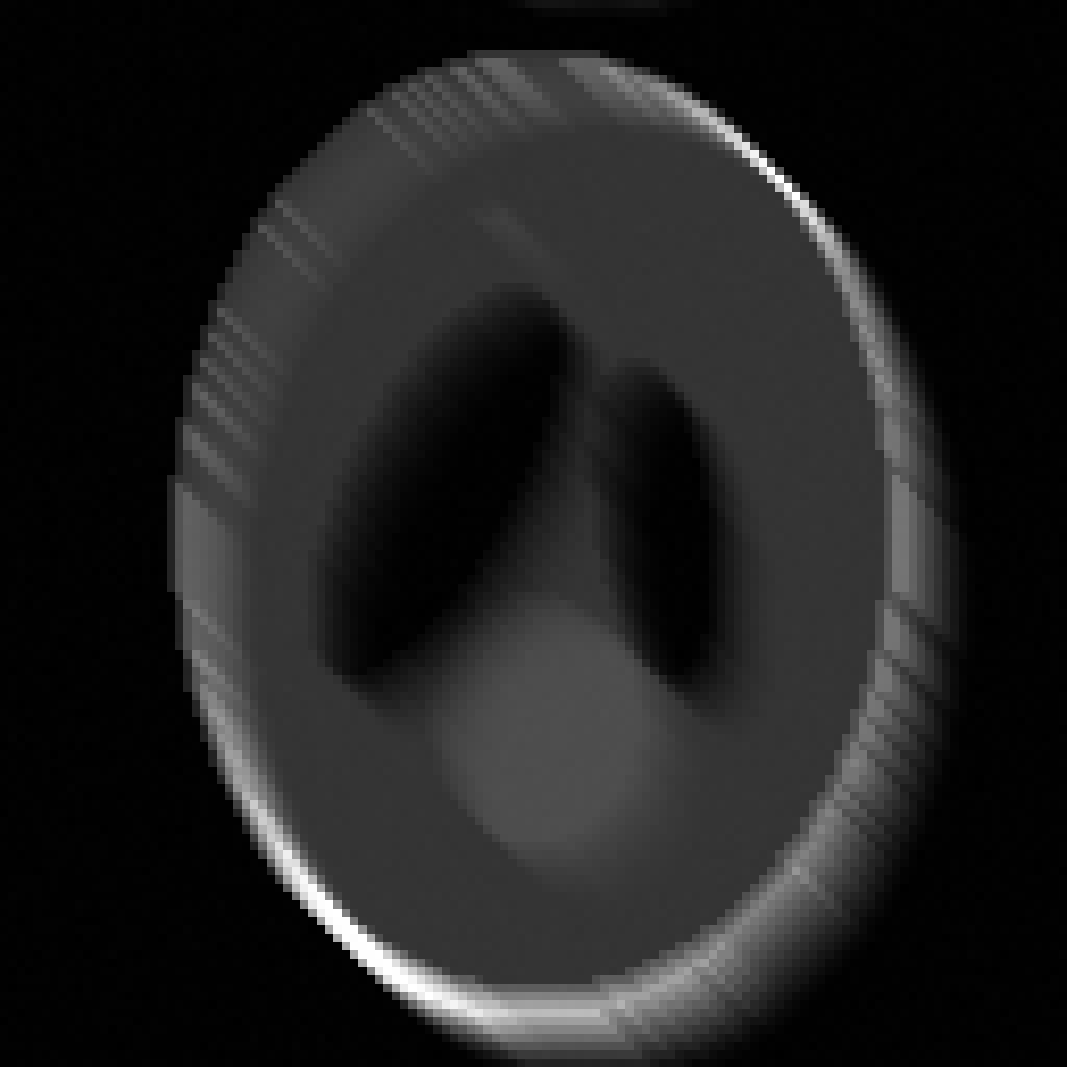}\\%(k)
	\end{minipage}
		\begin{minipage}{0.18\textwidth}
		\includegraphics[width=\textwidth]{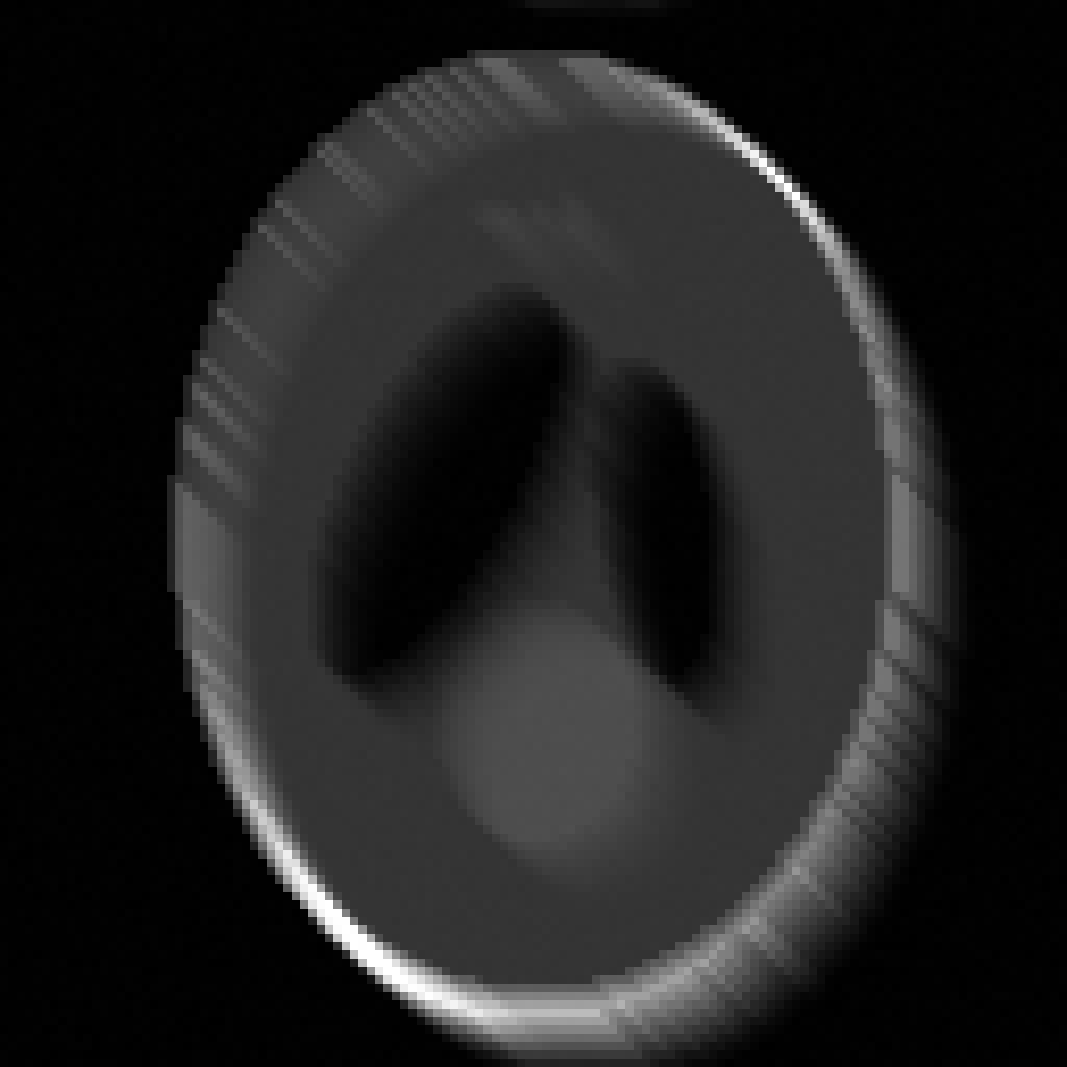}\\%(l)
	\end{minipage}
		\begin{minipage}{0.18\textwidth}
		\includegraphics[width=\textwidth]{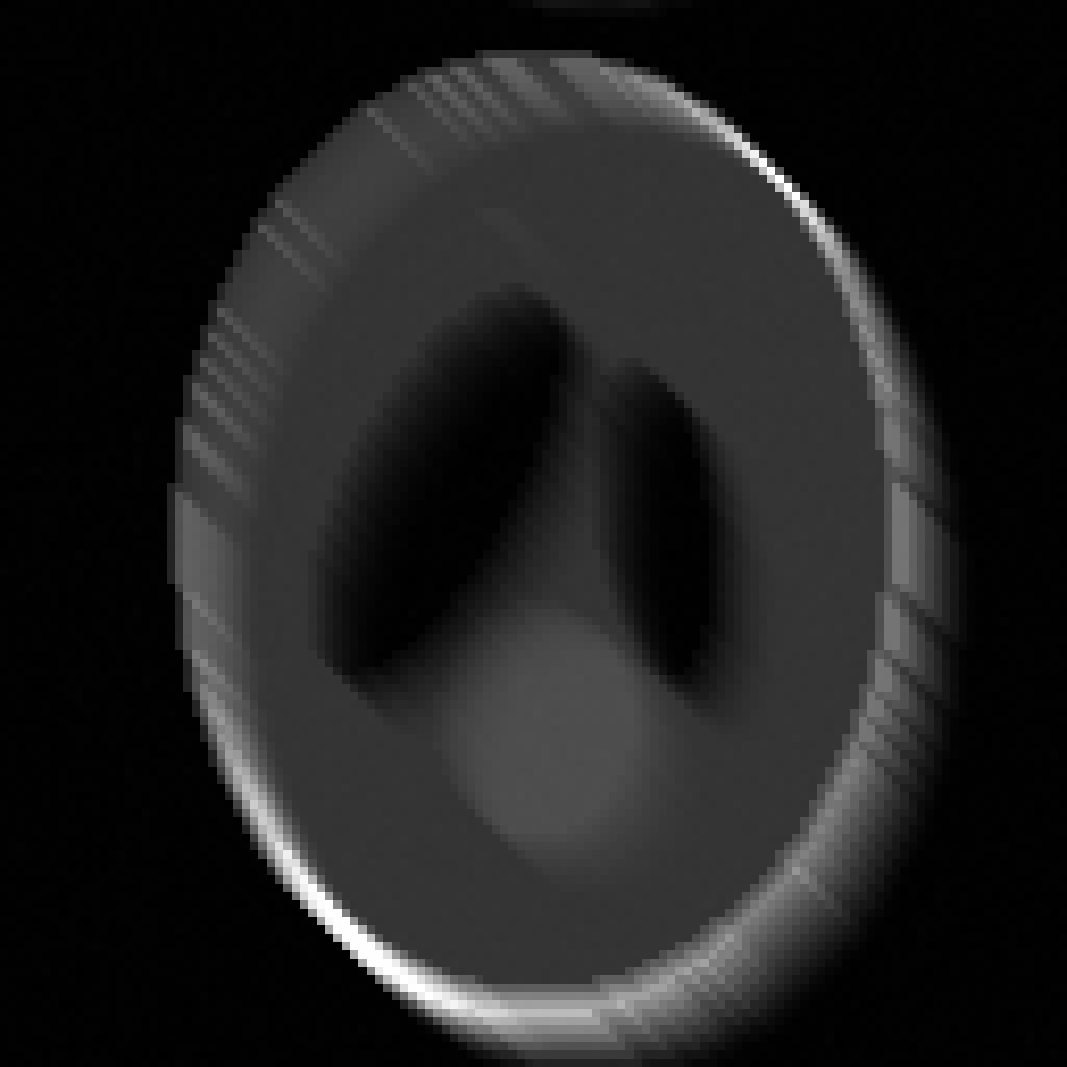}\\%(m)
	\end{minipage}

		\begin{minipage}{0.18\textwidth}
		\includegraphics[width=\textwidth]{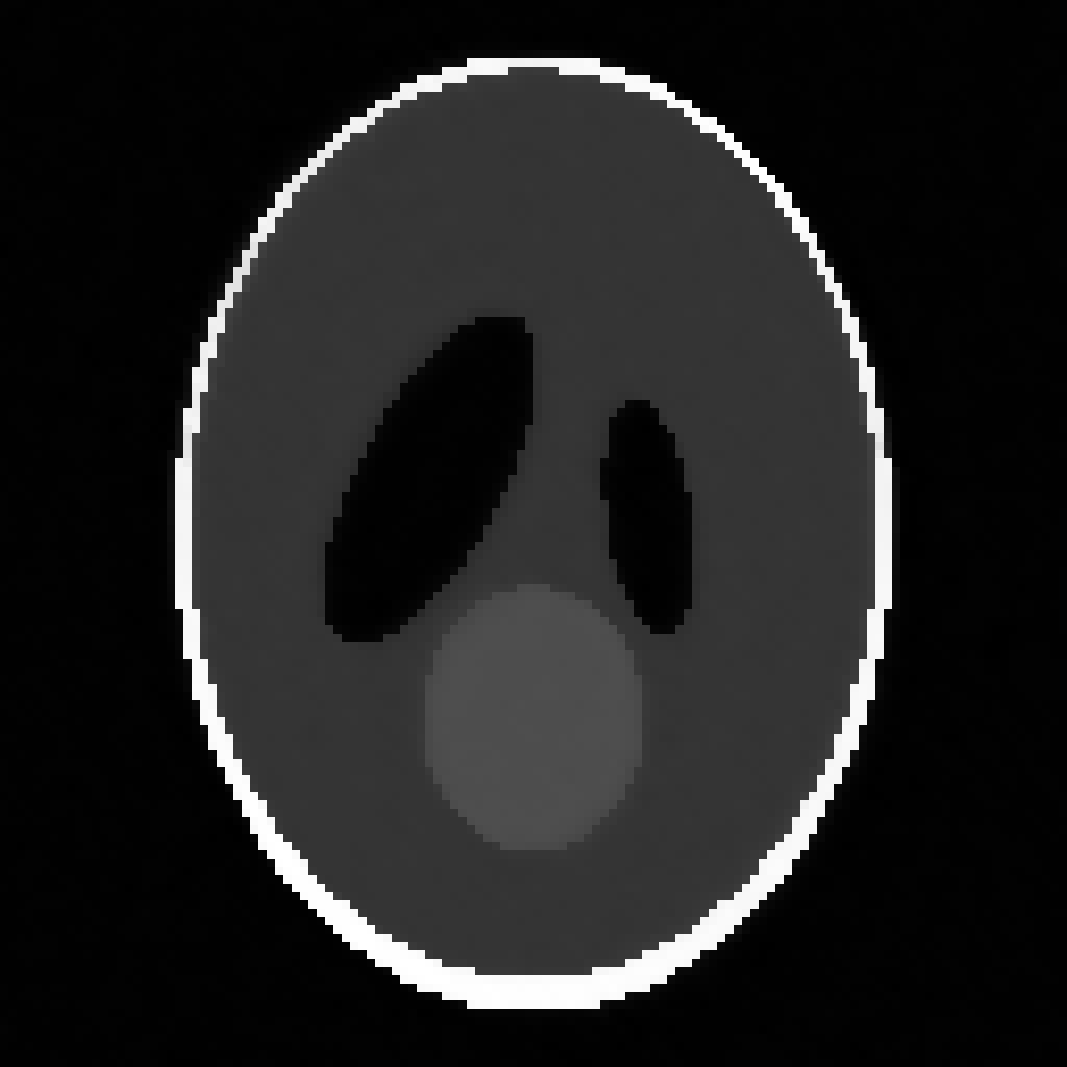}\\%(i)
	\end{minipage}
	\begin{minipage}{0.18\textwidth}
		\includegraphics[width=\textwidth]{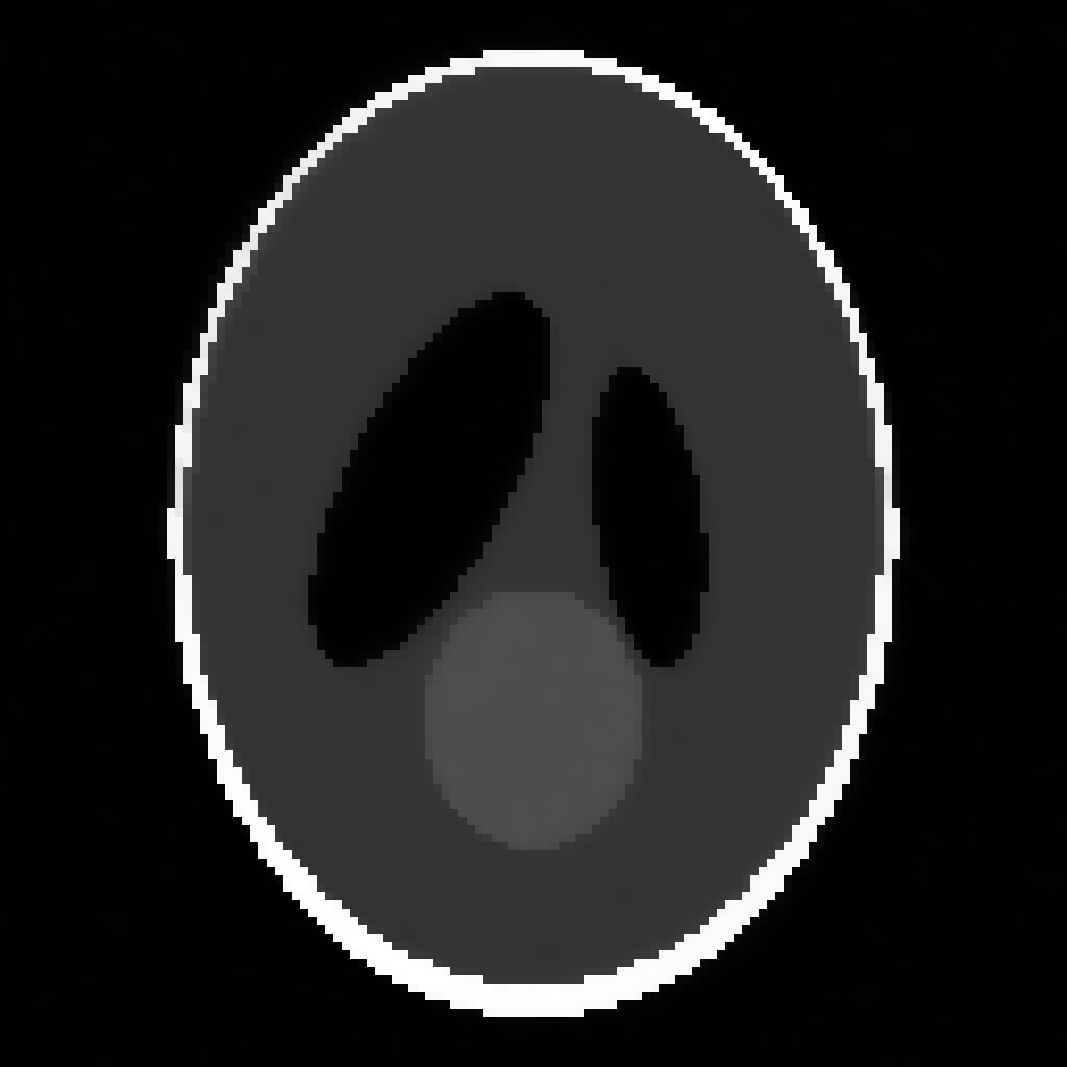}\\%(j)
	\end{minipage}
		\begin{minipage}{0.18\textwidth}
		\includegraphics[width=\textwidth]{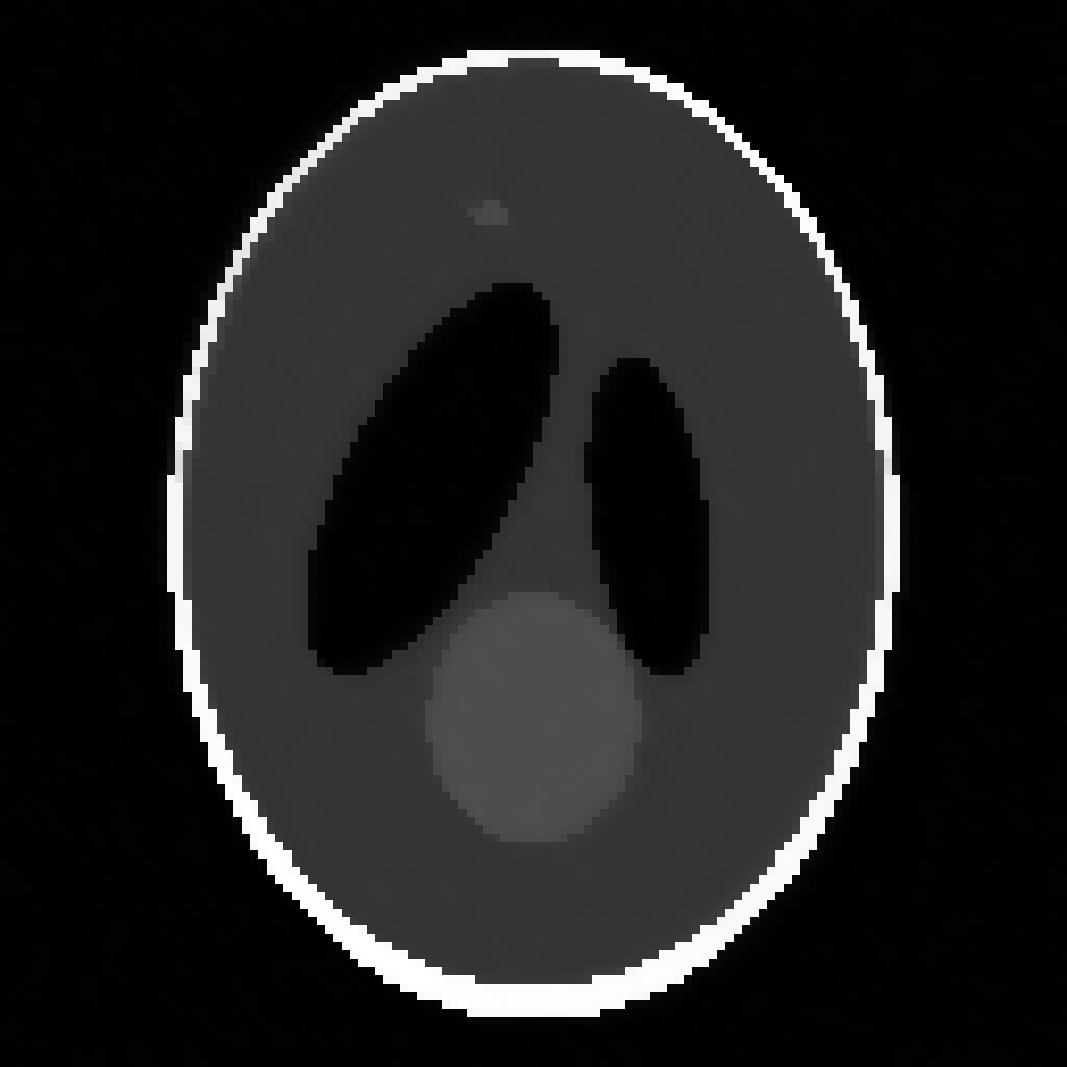}\\%(k)
	\end{minipage}
		\begin{minipage}{0.18\textwidth}
		\includegraphics[width=\textwidth]{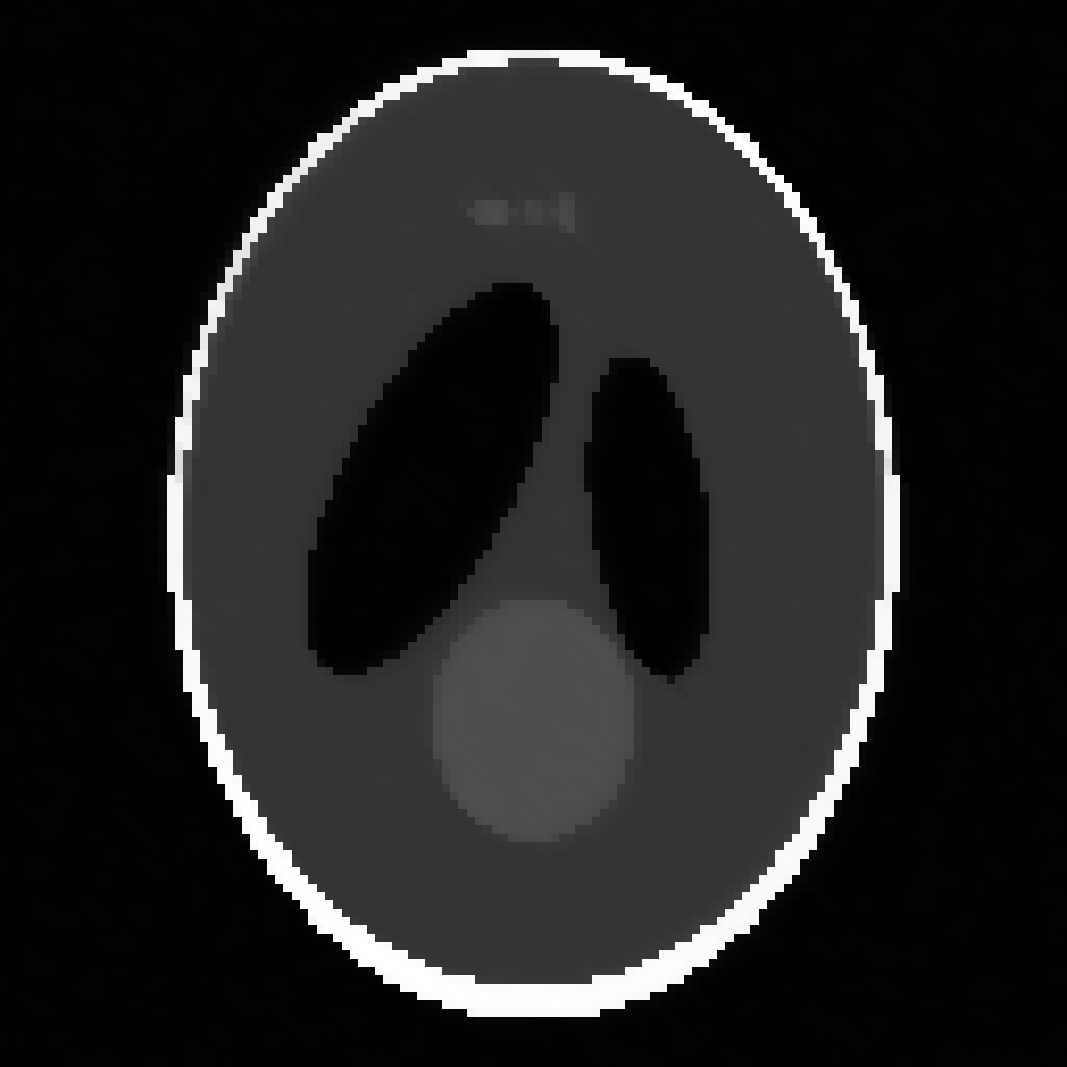}\\%(l)
	\end{minipage}
		\begin{minipage}{0.18\textwidth}
		\includegraphics[width=\textwidth]{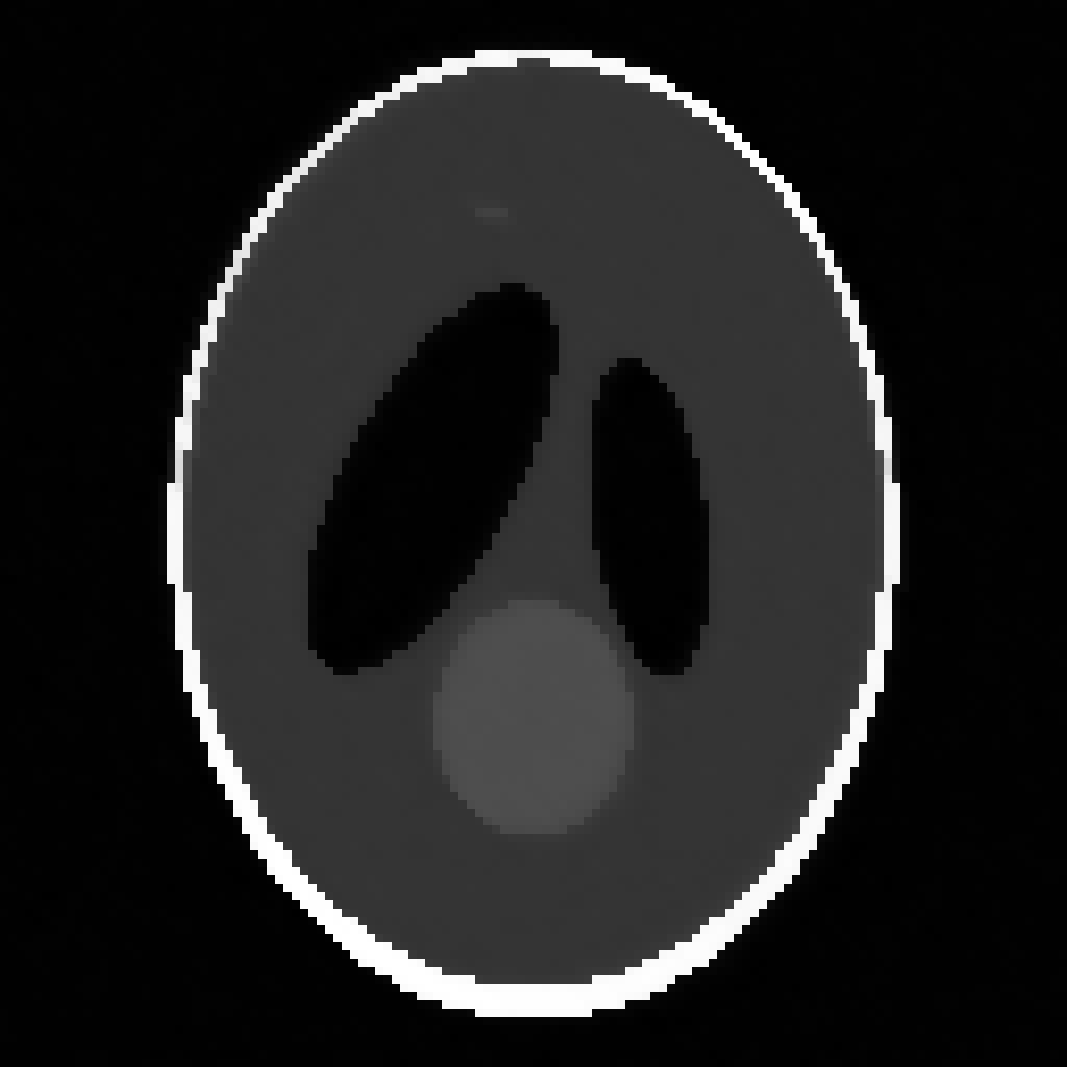}\\%(m)
	\end{minipage}

\end{center}
	\caption{Space-time image deblurring test problem: the first row represents a sample of true images at time steps $t = 1,3,5,6,7$. The second row shows the respective blurred any noisy images with medium blur and $1\%$ Gaussian noise. The third row shows the reconstructed images  $\bu^{(t)}$, $t = 1,3,5,6,7$ obtained by AnisoTV.}
	\label{Fig: ShepLogDeblur} 
\end{figure}

The goal here is to reconstruct a sequence of approximations of desired images from a sequence of blurry and noisy images. A sample of the true images is shown in the first row of Figure \ref{Fig: ShepLogDeblur}. The simulated available data are obtained by blurring 8 images of size $128\times 128$ with Gaussian point spread function with a medium blur using \cite{gazzola2019ir}. We consider all the operators $\bA = \bA^{(t)}\in \R^{16,384 \times 16,384}$, $t =1,2,\dots,8$ to be the same, resulting in the matrix $\bF = \bI_8 \otimes \bA \in \R^{131,072 \times 131,072}$. Therefore, the dynamic nature of the problem is characterized by images that change at different time instances and not from changes in the measurement process. The blurred images are perturbed with $1\%$ Gaussian noise and are shown in Figure \ref{Fig: ShepLogDeblur}. We solve~\eqref{eqn:dynamic} where $i = 1,2,\dots, 6$ corresponds to methods AnisoTV, TVplusTikhonov, IsoTV, Aniso3DTV, Iso3DTV, and GS, respectively. 
Some quantitative results are displayed in Figure~\ref{Fig: ErrorSL}. 

Figure~\ref{Fig: ErrorSL} (a) shows the RRE at each time point for each method. The RRE is computed at the iteration when the discrepancy principle \eqref{eq: DP} is first satisfied; the number of iterations  and the regularization parameter $\lambda$ that was chosen are displayed in Table~\ref{tab:ex1}; note that we estimate the corresponding regularization parameter at each MM-Krylov iteration.  
In Figure \ref{Fig: ErrorSL} (b), we show the convergence history for all the methods when each method is allowed to run for 150 iterations without considering any other stopping criteria. Solid diamond markers over the lines in Figure \ref{Fig: ErrorSL} (b) show the iteration and the value of the RRE when the discrepancy principle is satisfied. Notice that each line in Figure \ref{Fig: ErrorSL} (b) shows the convergence for each method for all images together, that is, the convergence for $\bu_{(k)}$.  
We observe that AnisoTV and GS outperform the other methods for this example. Moreover, as illustrated in Figure \ref{Fig: ErrorSL}. For methods IsoTV, Iso3DTV, and TVplusTikhonov we observe an increase of the RRE in the early iterations, but if the method is let to run enough iterations, then the convergence behaviour starts to stabilize. 
Reconstructions with AnisoTV at time steps $t = 1,3,5,6,7$ are shown in the third row of Figure \ref{Fig: ShepLogDeblur}.

\begin{table}[!ht]
    \centering
    \caption{Space-time image deblurring example: The number of iterations when the discrepancy principle is satisfied for the first time and the corresponding regularization parameters for AnisoTV, TVplusTikhonov, IsoTV, Aniso3DTV, Iso3DTV, and GS.}
    \label{tab:ex1}
    \begin{tabular}{l|c|c|c|c|c|c}
       -  &  AnisoTV&  TVplusTikhonov &  IsoTV & Aniso3DTV&  Iso3DTV & GS\\ \hline
      Iters   & 56 &  58 & 49 & 60 & 50 &  46 \\ 
      $\lambda$ & 0.24 &  0.4 & 0.45 & 0.36 &  0.27 & 0.3
    \end{tabular}
\end{table}
\subsection{Example 2: Dynamic photoacoustic tomography (PAT)} 
As a second example, we consider PAT that is a hybrid imaging modality that combines the rich contrast of optical imaging with the high resolution of ultrasound imaging. 

A continuous model for PAT reconstruction for a single reconstruction under motion is studied in \cite{chung2017motion} while the PAT model for a sequence of images in the Bayesian setting involving Mat\'ern spatial and temporal type of priors.

\begin{figure}[!ht]
    \centering

{ \includegraphics[width = 0.23\textwidth]{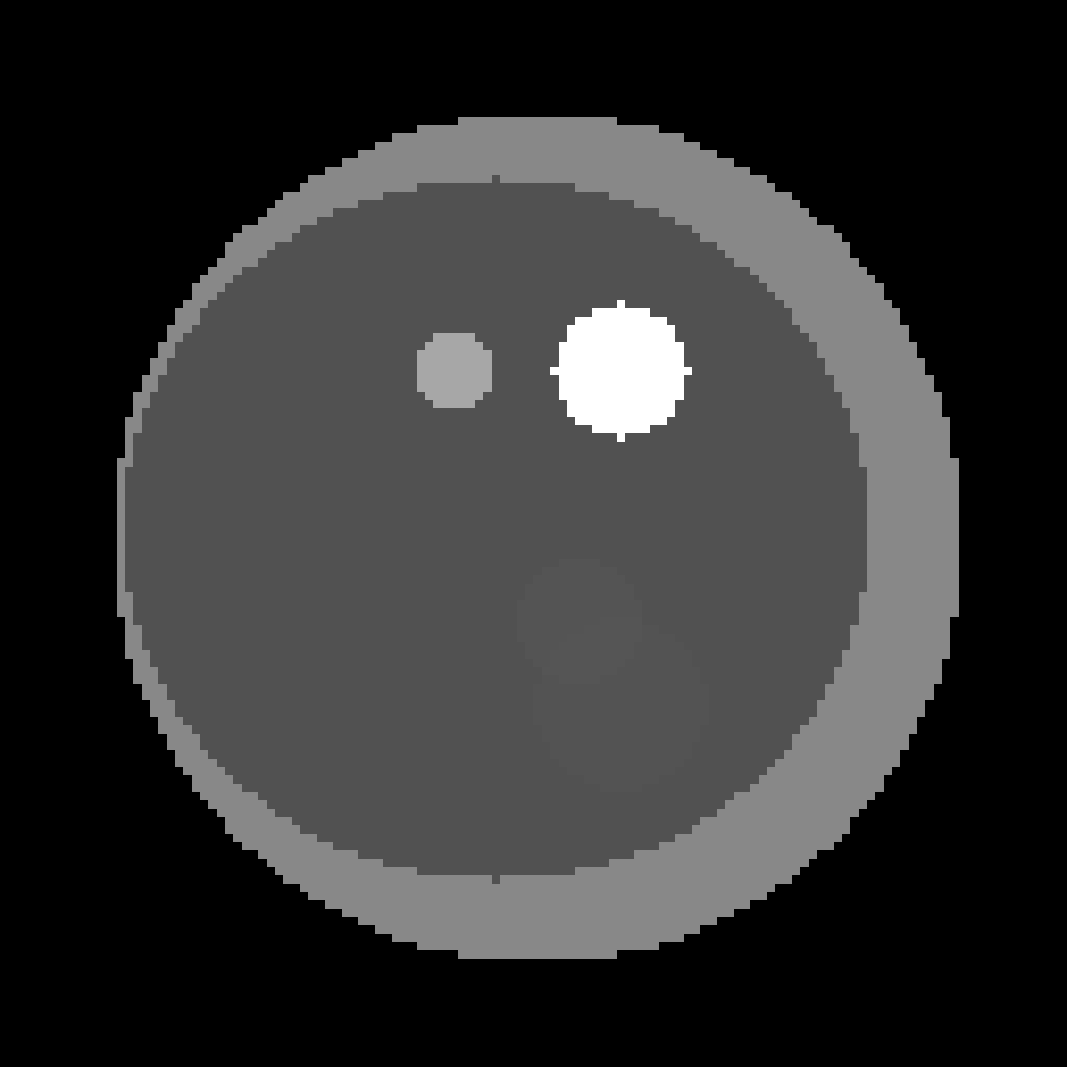}
    \includegraphics[width = 0.23\textwidth]{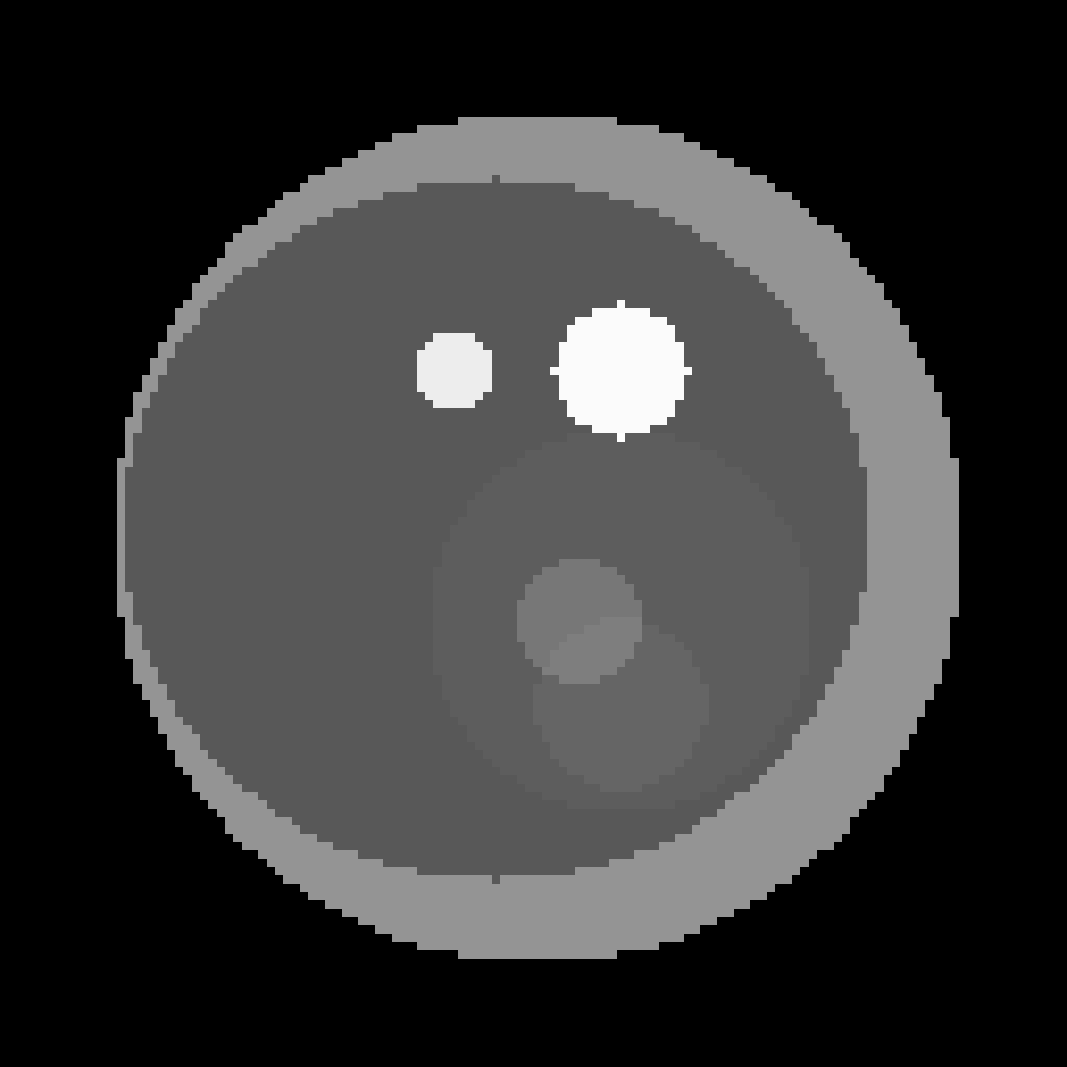} 
    \includegraphics[width = 0.23\textwidth]{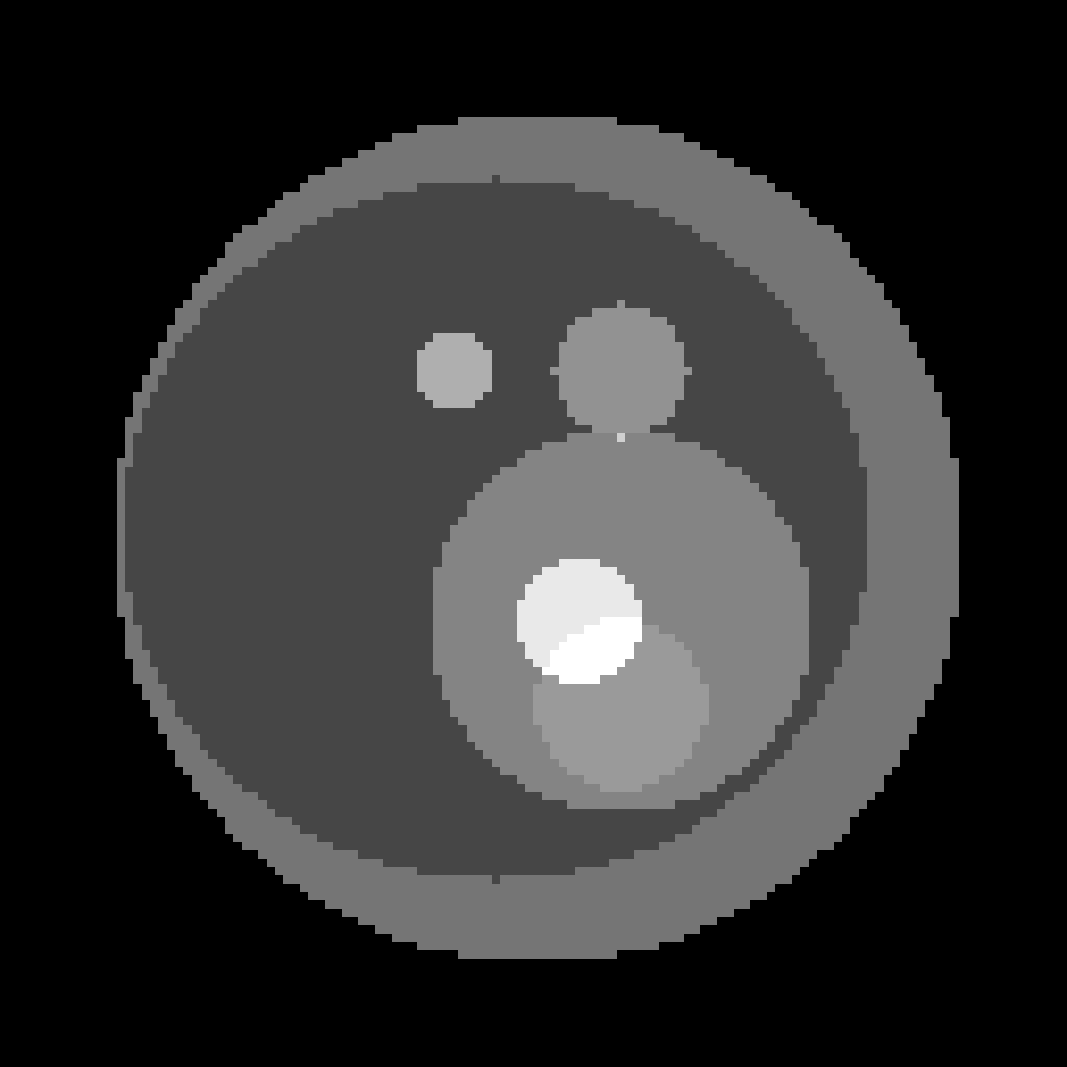}
    \includegraphics[width = 0.23\textwidth]{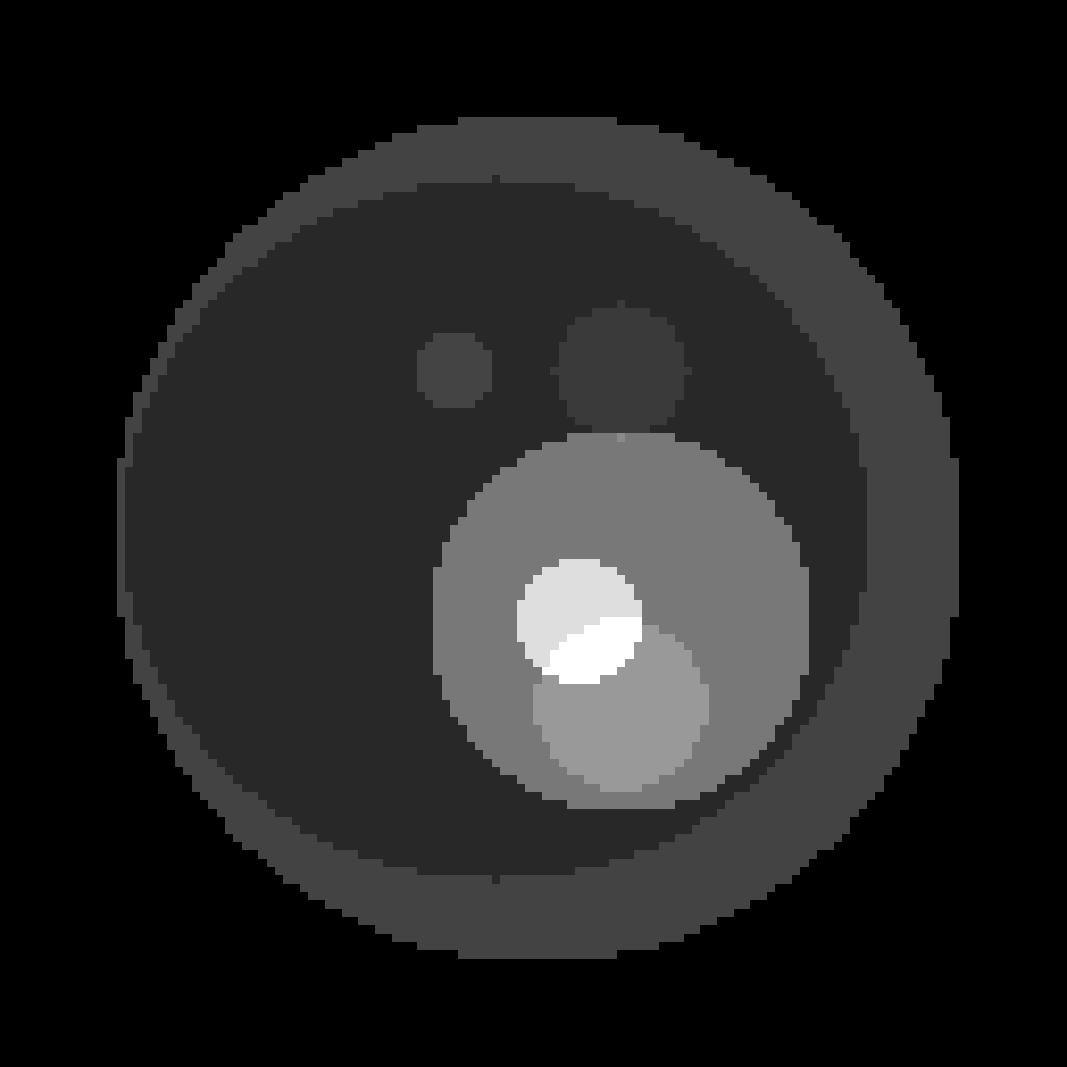}
    \label{fig:lorenz_extrinsic_inv_single}
    } %\\\\

{ \includegraphics[height=0.2\textwidth,width = 0.17\textwidth]{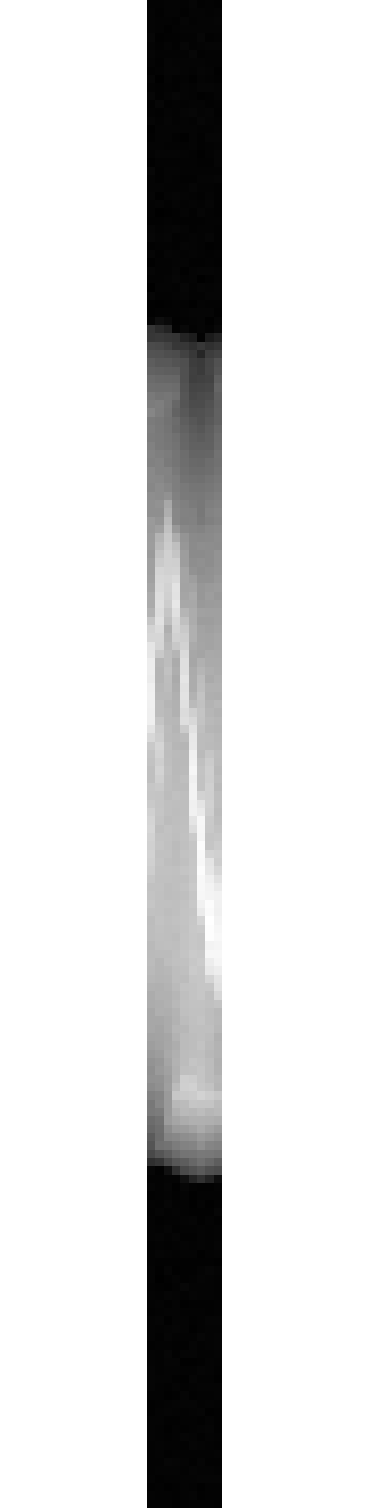}
\includegraphics[height=0.2\textwidth,width = 0.17\textwidth]{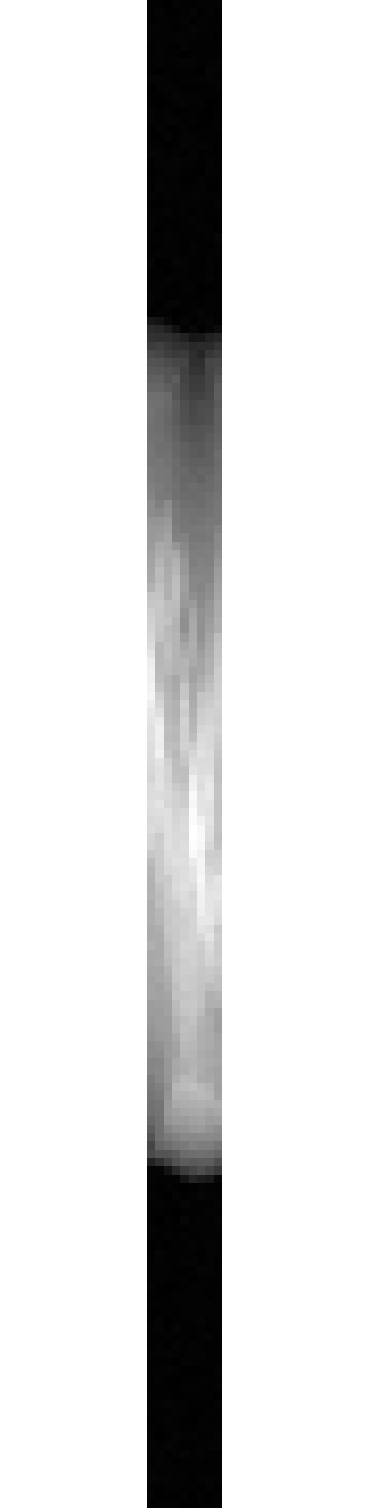}
\includegraphics[height=0.2\textwidth,width = 0.17\textwidth]{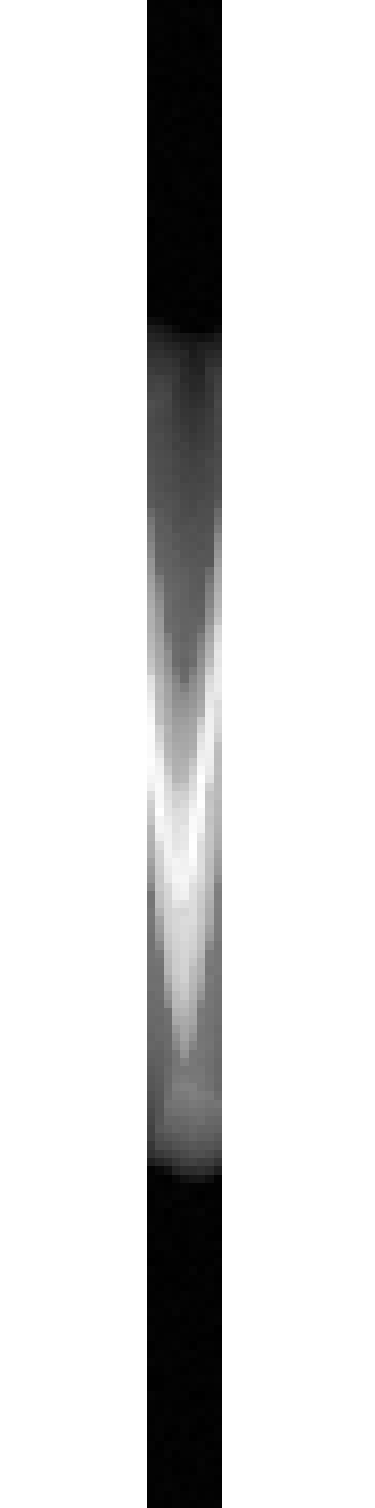}
\includegraphics[height=0.2\textwidth,width = 0.17\textwidth]{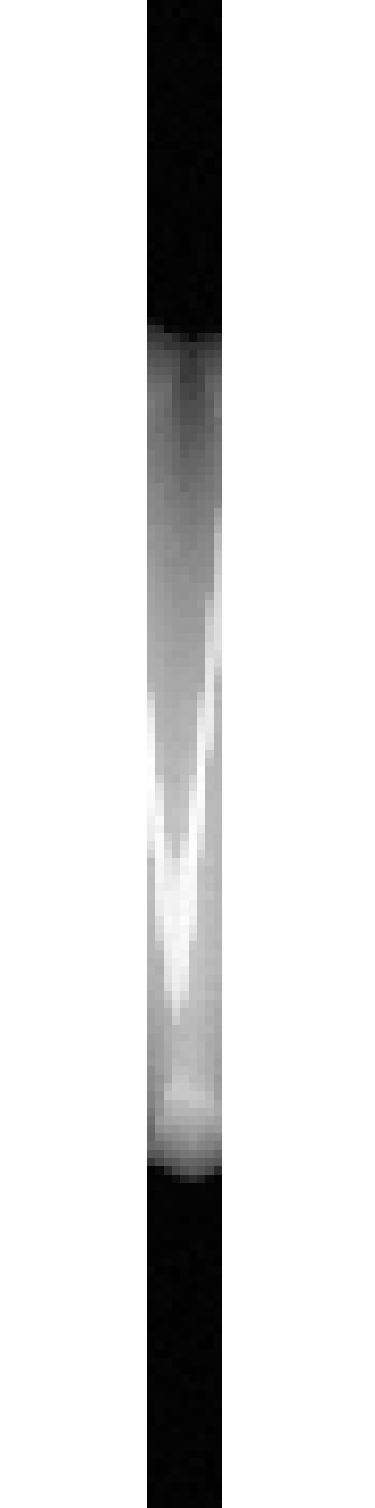}
\includegraphics[height=0.2\textwidth,width = 0.23\textwidth]{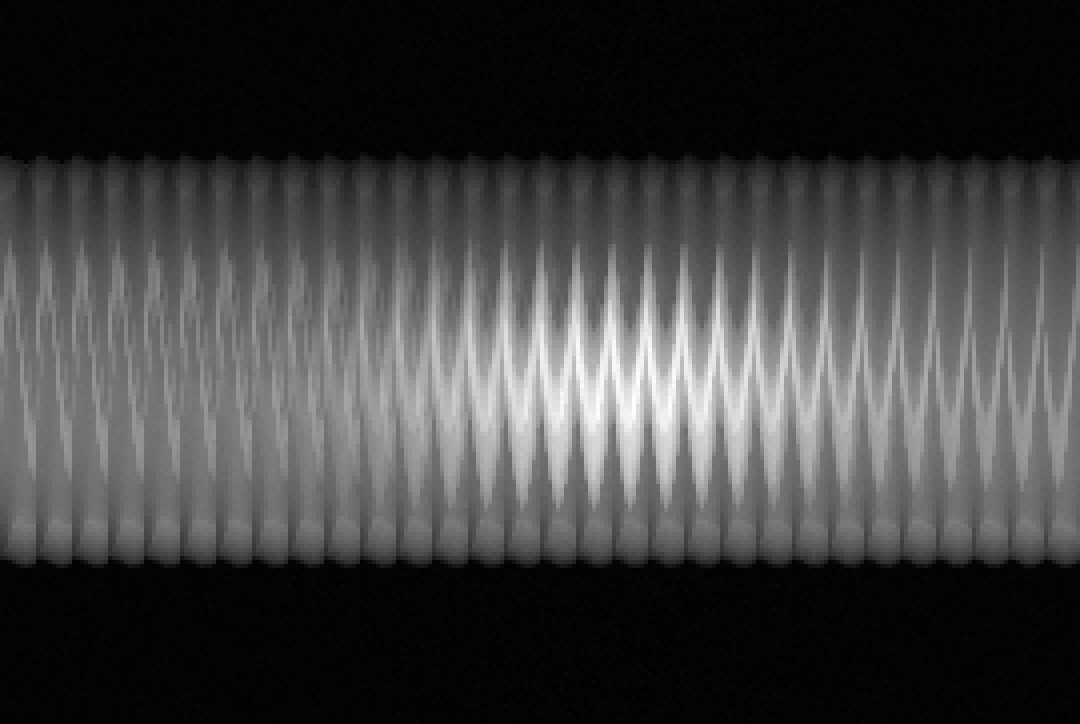}
\label{fig:lorenz_extrinsic_inv_mul}}
	\caption{PAT test problem: First row, from left to right: True images at time steps $t = 1,10,20,30$. Second row, from left to right: sample of sinograms at time steps $t=1,10,20,30$ and the full sinogram.}
\label{Fig: PatTrue}
\end{figure}

In this example, we consider a dynamic inverse problem in which the forward operator is time-dependent (see Section~\ref{sec:intro}), so that the operator $\bF$ has the blockdiagonal structure \eqref{eq: blockF} and the number of time points are $n_t=30$. The operator $\bA^{(t)}$ corresponds to the projection angles $t,t+30,\dots,t+269$, and for each angle there are $362$ measurements. Each image 
$\bU^{(t)}$
is of size $256\times 256$ and represents a superposition of six circular objects that are in motion. This implies that the total number of unknowns is $ n = 256\times 256\times 30 = 1,966,080$. A sample of true images at time instances $t = 1, 10, 20, 30$ is shown at the first row of Figure \ref{Fig: PatTrue}. We add $1$\% white Gaussian noise to the available measurements and the resulting noisy sinograms $\bd^{(t)} \in \R^{3258}$ at time steps $t = 1,10,20,30$ along with the total sinogram (obtained by concatenating all 30 available sinograms together) % of size $3,258 \times 30$
with a total of 
$m=97,740$ observations are shown in the second row of Figure \ref{Fig: PatTrue}.
Note that this inverse problem is severely underdetermined.

\begin{figure}[h!]
\centering
    \begin{minipage}{0.45\textwidth}
		\includegraphics[width=\textwidth]{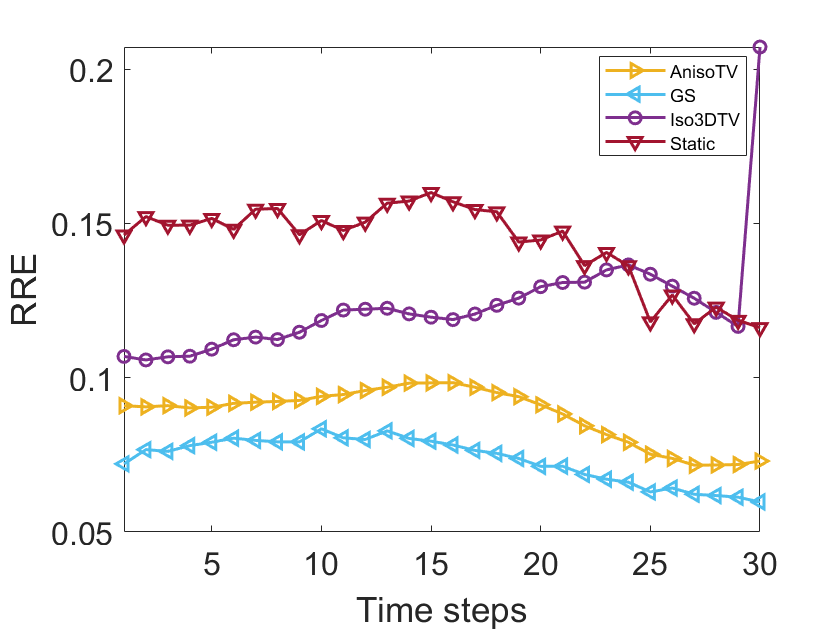}\\(a)
	\end{minipage}
\centering
	 \begin{minipage}{0.45\textwidth}
		\includegraphics[width=\textwidth]{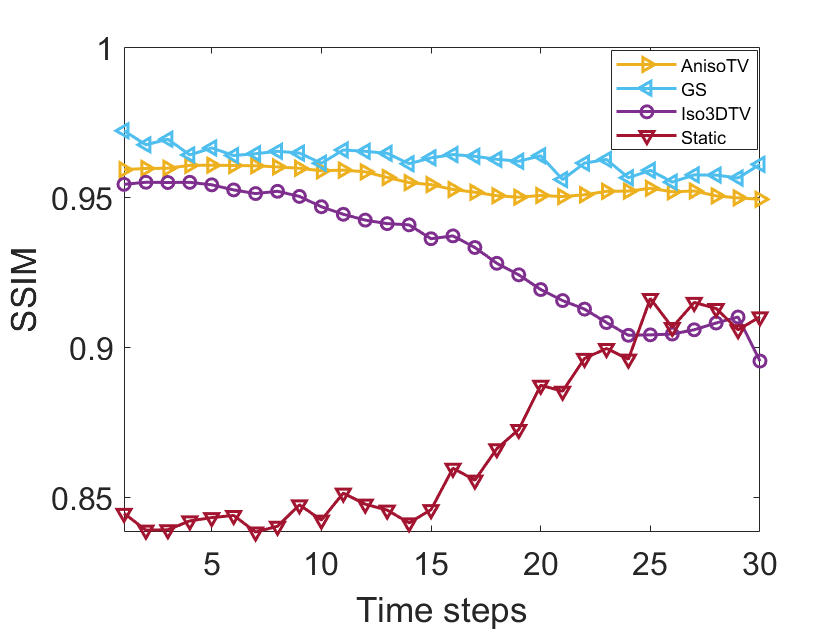}\\(b)
	\end{minipage}
	\caption{PAT test problem: a) RRE computed at the iteration when the DP \eqref{eq: DP} is first satisfied,  for each time step. b) SSIM computed at the iteration when the DP \eqref{eq: DP} is first satisfied,  for each time step. The methods considered here are AnisoTV (right-pointing triangle line), Iso3DTV (dotted line), and GS (left-pointing triangle line). The asterisk line represents the figures of merit for the static problem.}
   \label{Fig: RREPAT} 
\end{figure}

We carry out the following numerical experiments:
\begin{enumerate}[(i)]
    \item Solve the large-scale dynamic inverse problem~\eqref{eqn:dynamic} with $j = 1, 4, 6$. More specifically, we choose AnisoTV from anisotropic-type methods, Iso3DTV from isotropic-type methods, and GS.  
    \item  Solve the static inverse problem~\eqref{eq: static}  with the regularization term $\mathcal{R}(\bu) = \|\bL_s\bu^{(t)}\|_1$, at $t =1,2,\dots,30$.
\end{enumerate}

\begin{figure}[h!]
	\centering
	\begin{minipage}{0.2\textwidth}
		\centering
		\includegraphics[width=\textwidth]{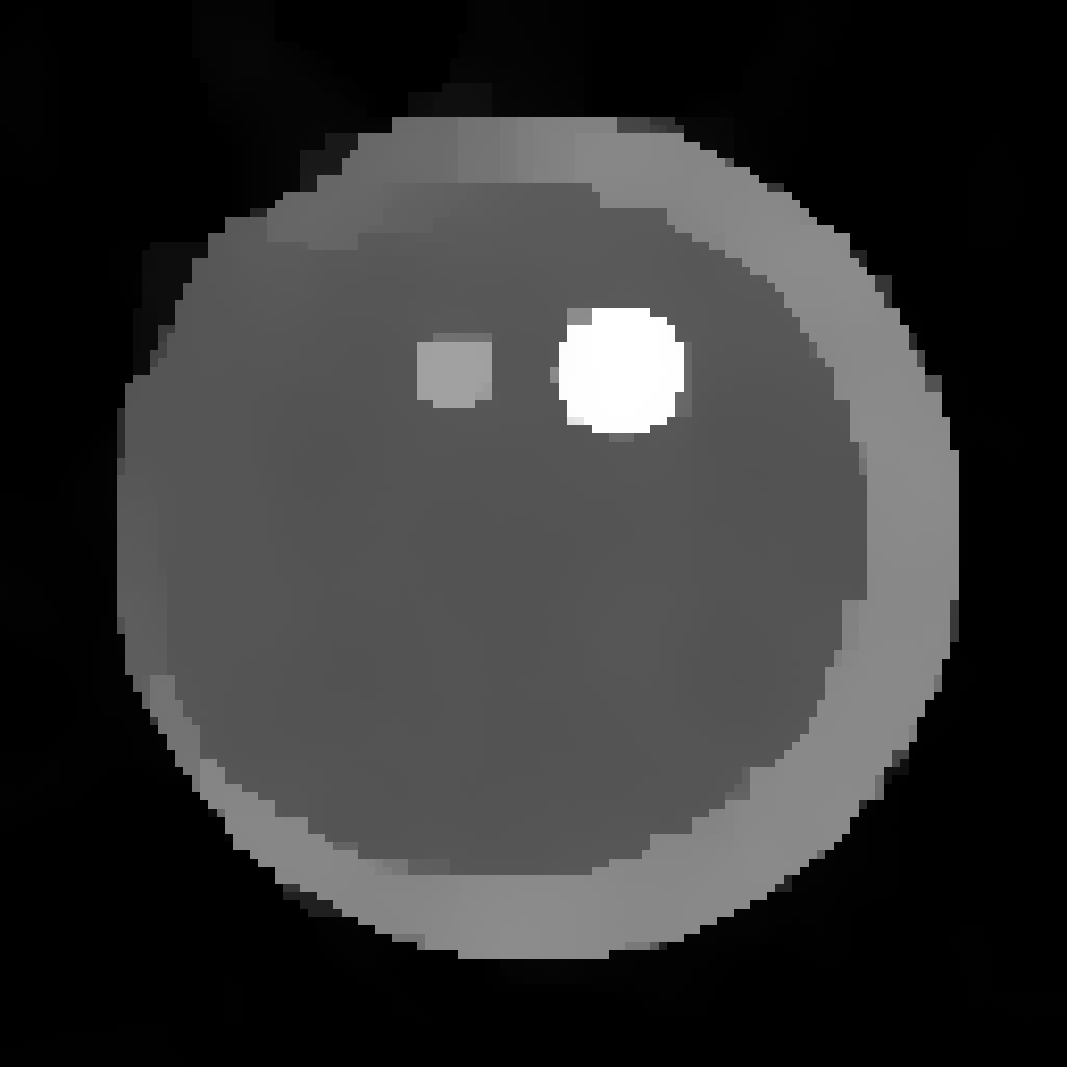}\\%(a)
	\end{minipage}
	\begin{minipage}{0.2\textwidth}
		\centering
		\includegraphics[width=\textwidth]{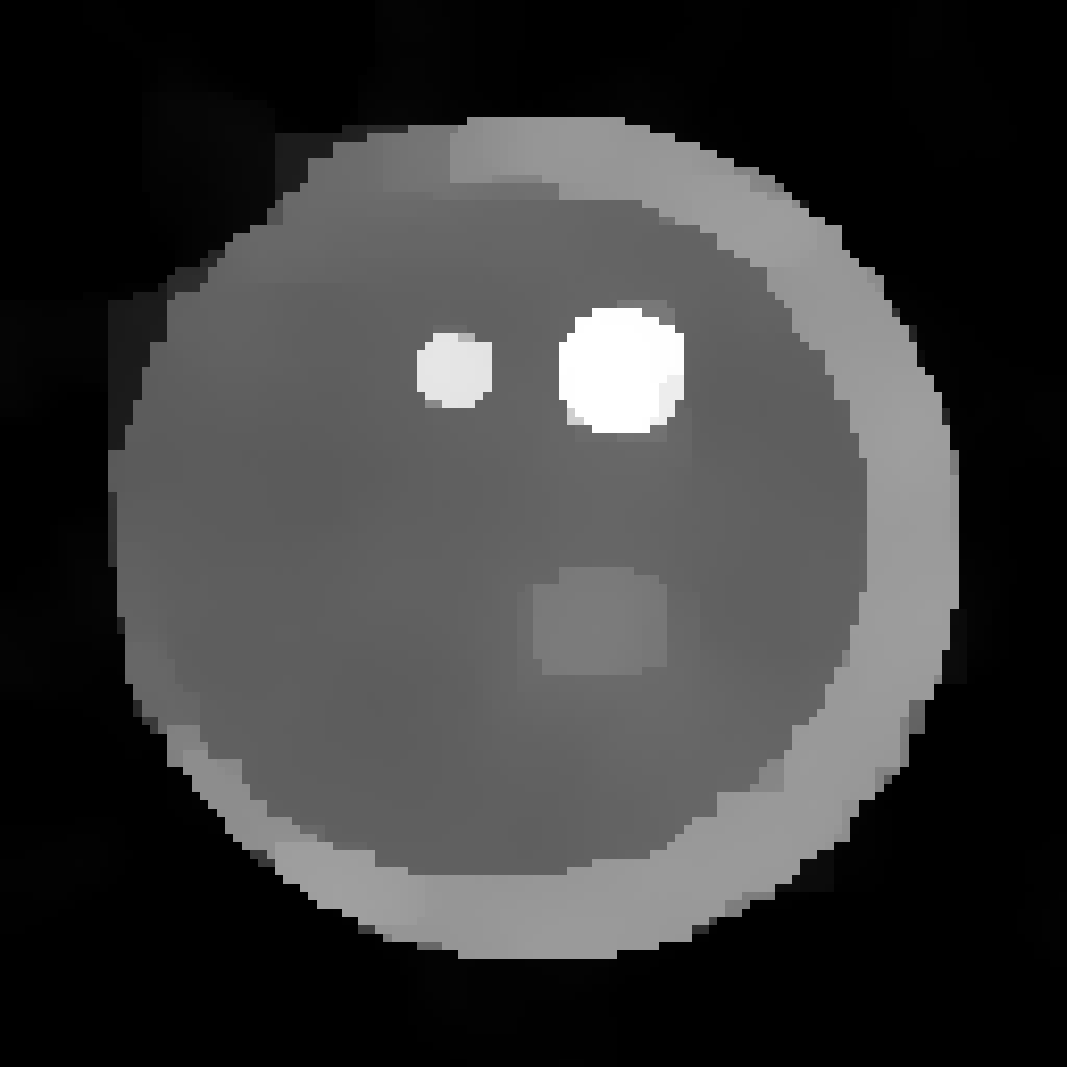}\\%(b)
	\end{minipage}
	\begin{minipage}{0.2\textwidth}
		\centering
		\includegraphics[width=\textwidth]{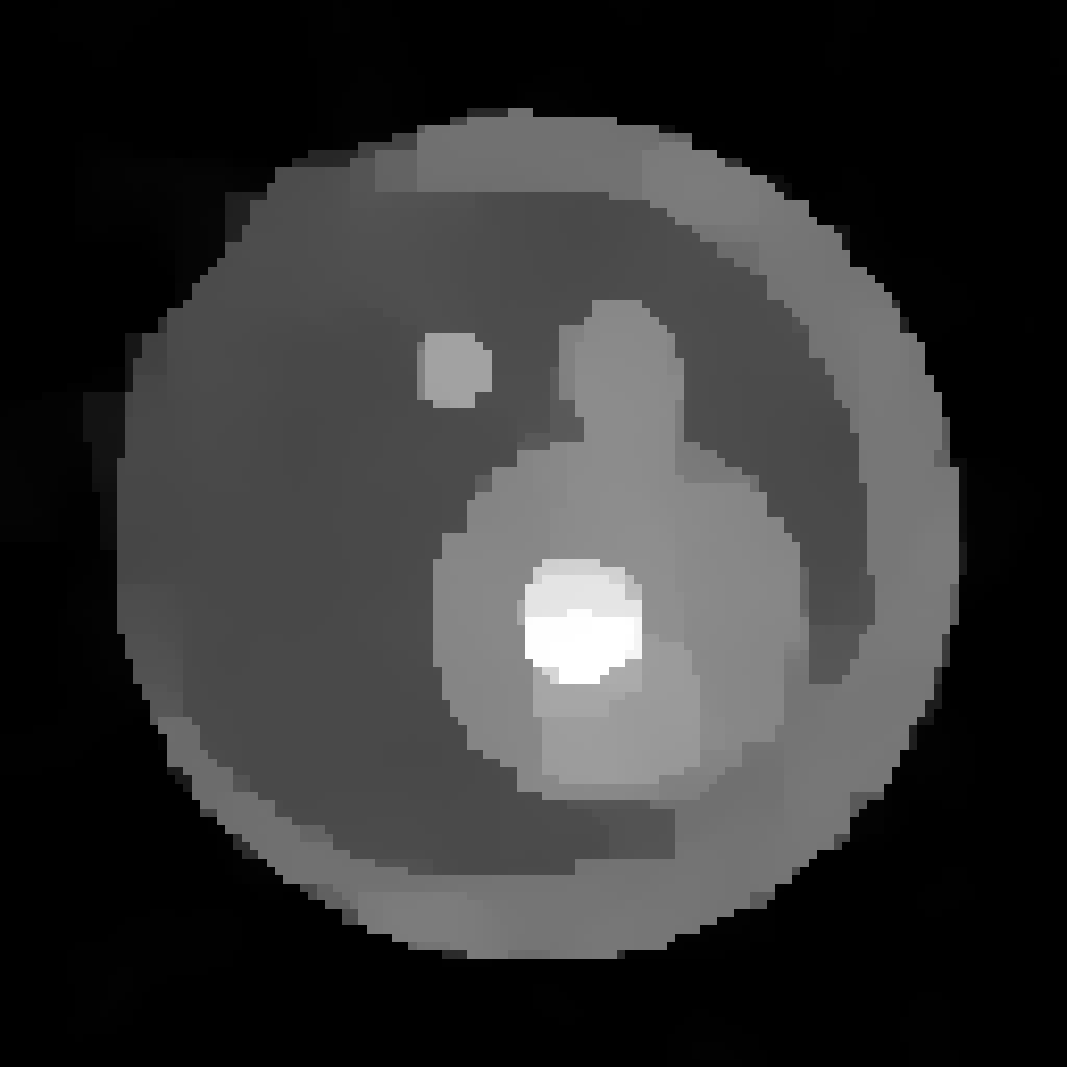}\\%(c)
	\end{minipage}
		\begin{minipage}{0.2\textwidth}
		\centering
		\includegraphics[width=\textwidth]{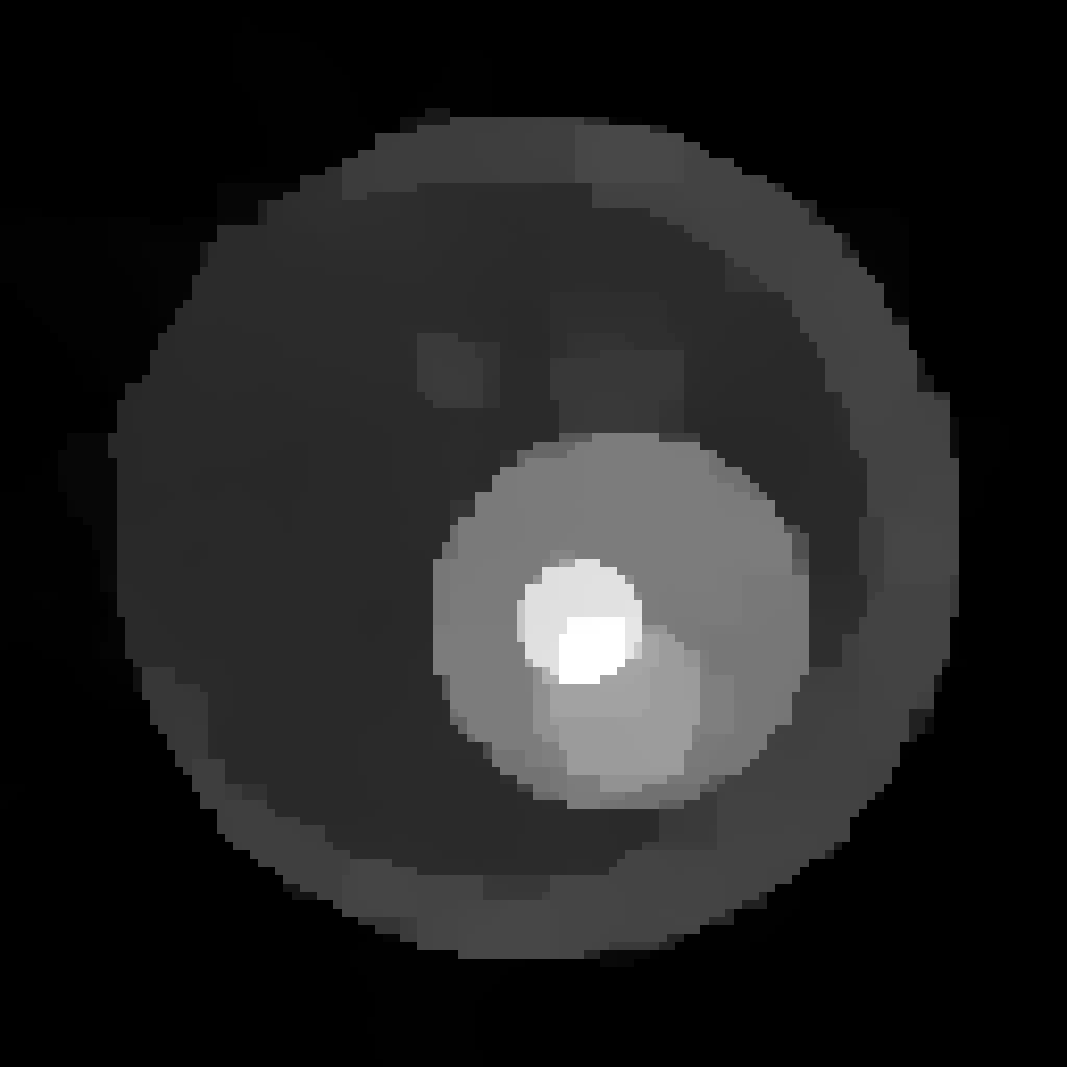}\\%(d)
	\end{minipage}

		\begin{minipage}{0.2\textwidth}
		\centering
		\includegraphics[width=\textwidth]{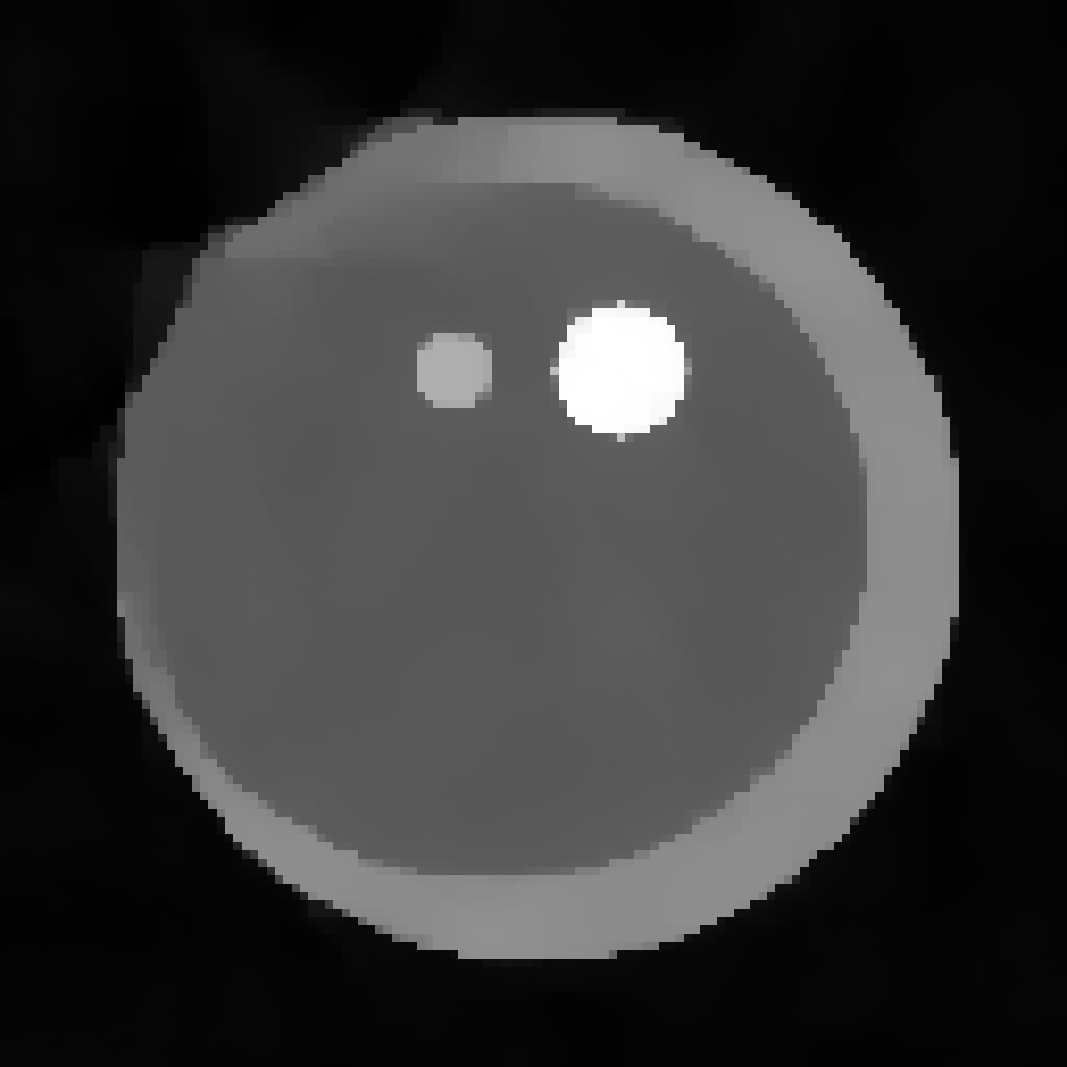}\\%(i)
	\end{minipage}
	\begin{minipage}{0.2\textwidth}
		\centering
		\includegraphics[width=\textwidth]{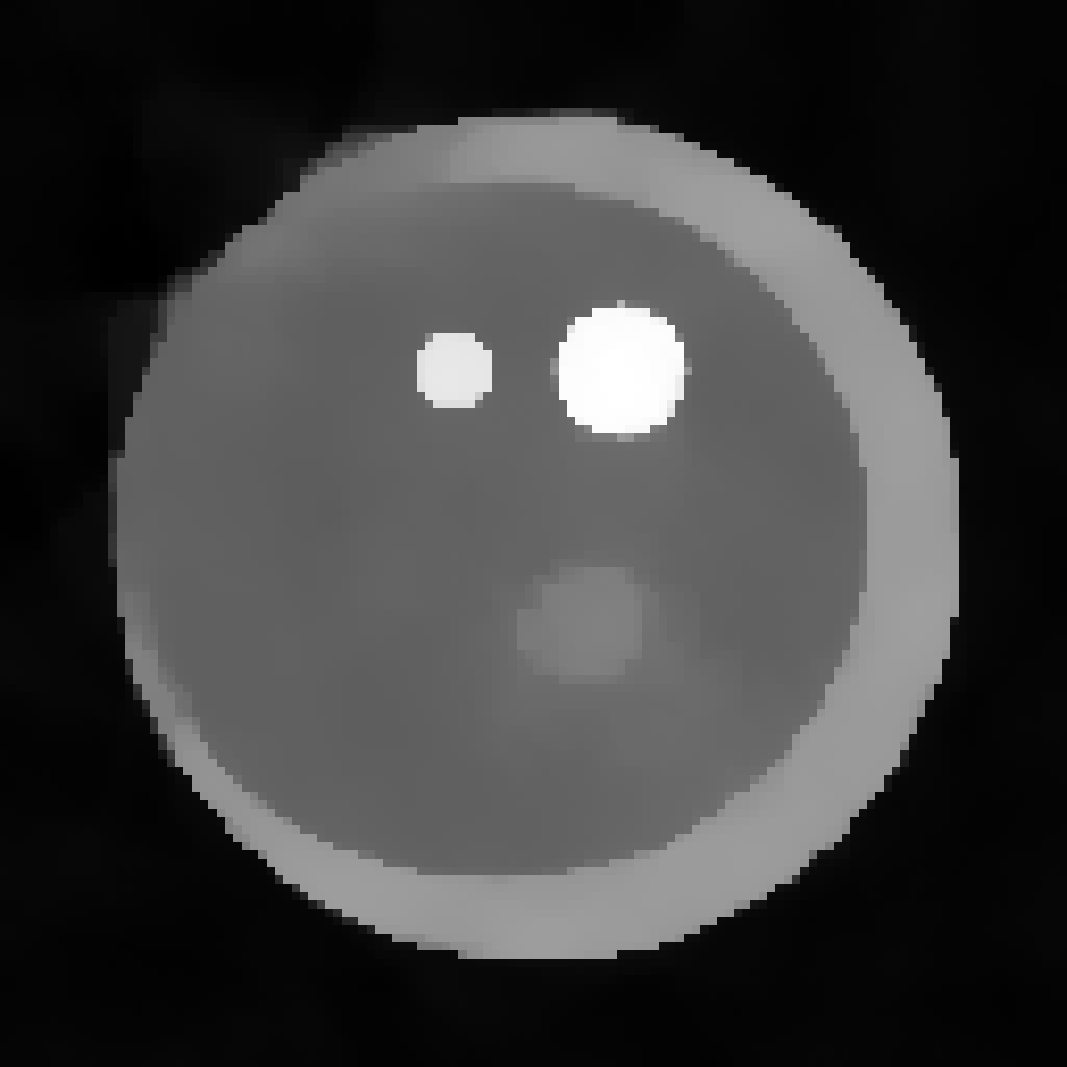}\\%(j)
	\end{minipage}
	\begin{minipage}{0.2\textwidth}
		\centering
		\includegraphics[width=\textwidth]{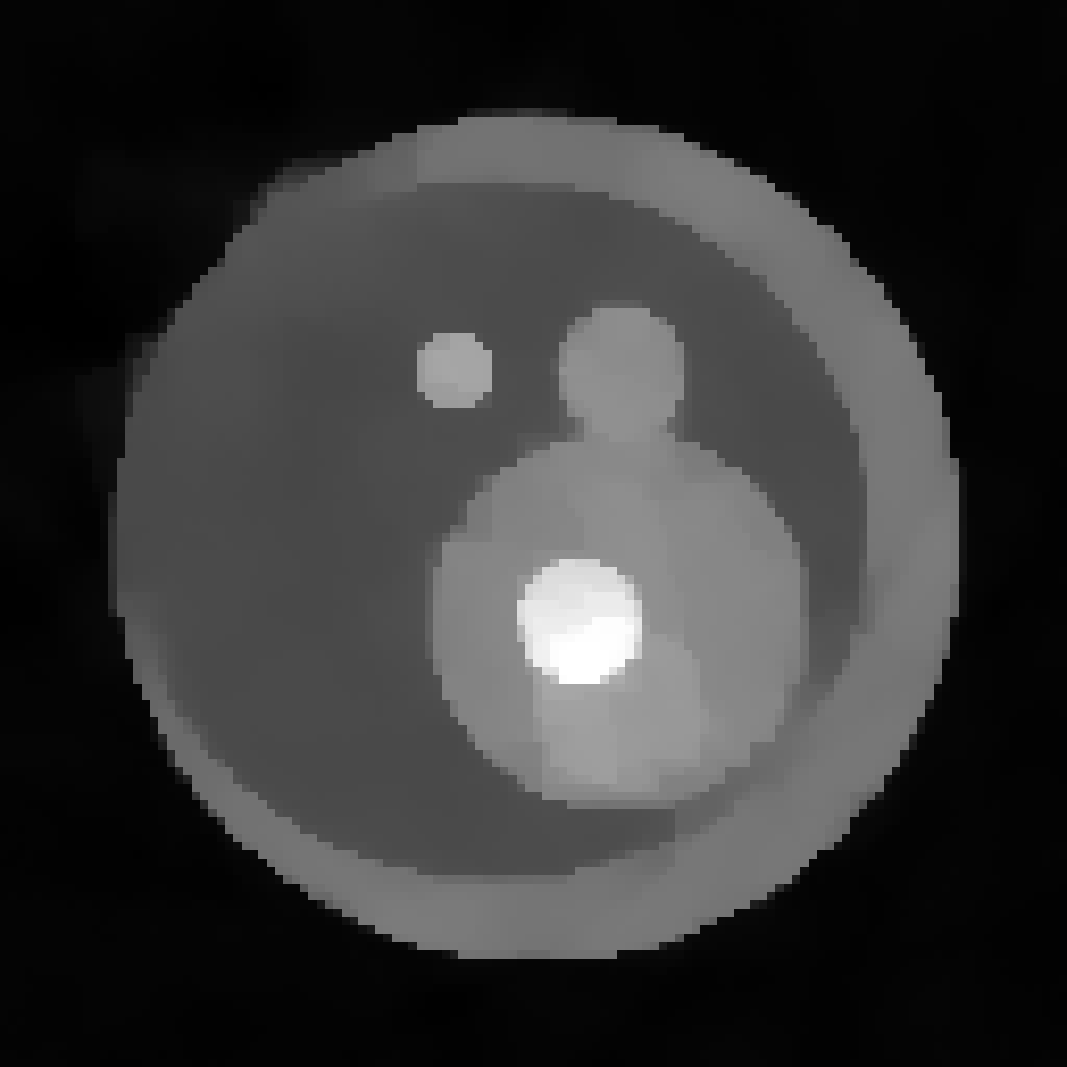}\\%(k)
	\end{minipage}
	\begin{minipage}{0.2\textwidth}
		\centering
		\includegraphics[width=\textwidth]{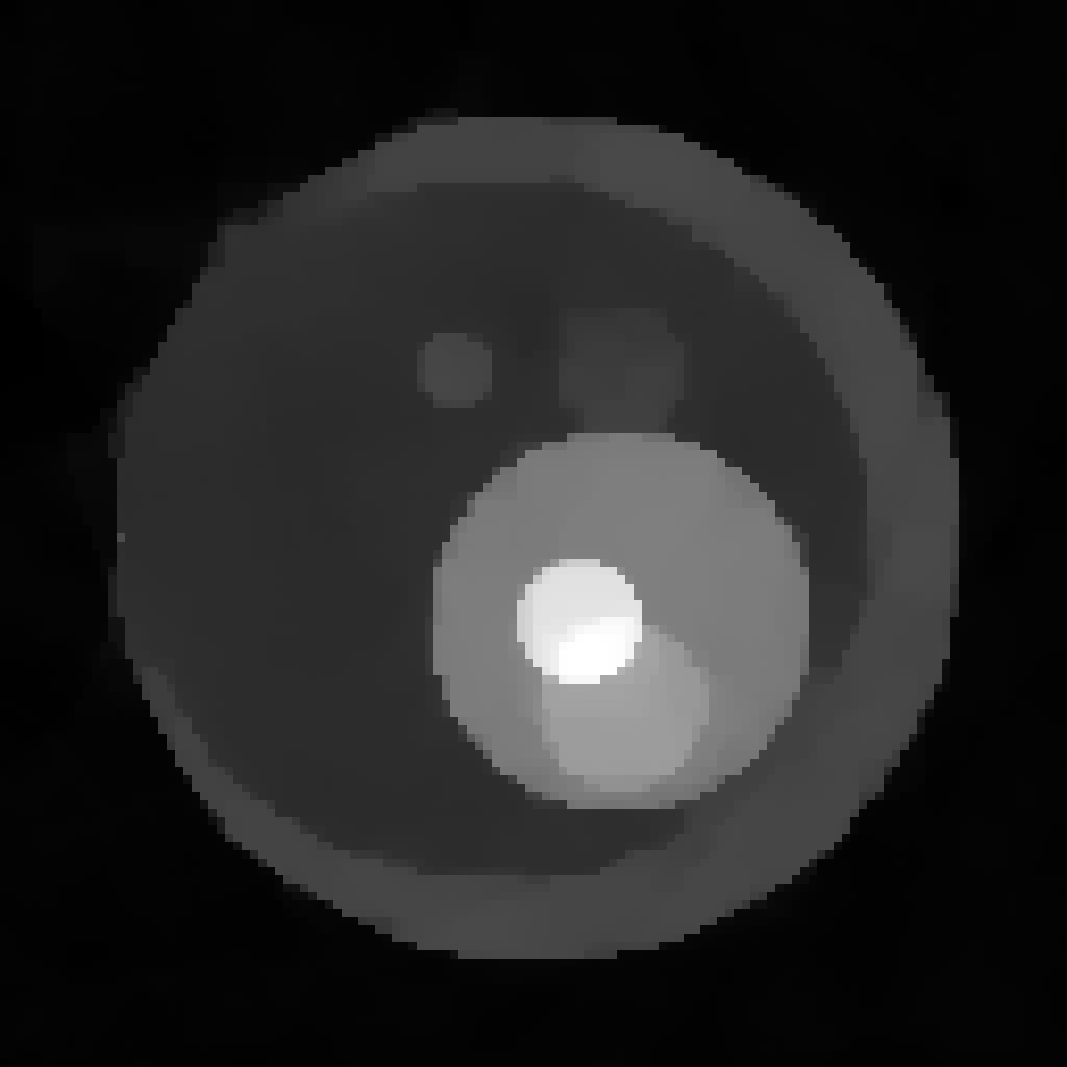}\\%(l)
	\end{minipage}
	
		\begin{minipage}{0.2\textwidth}
		\centering
		\includegraphics[width=\textwidth]{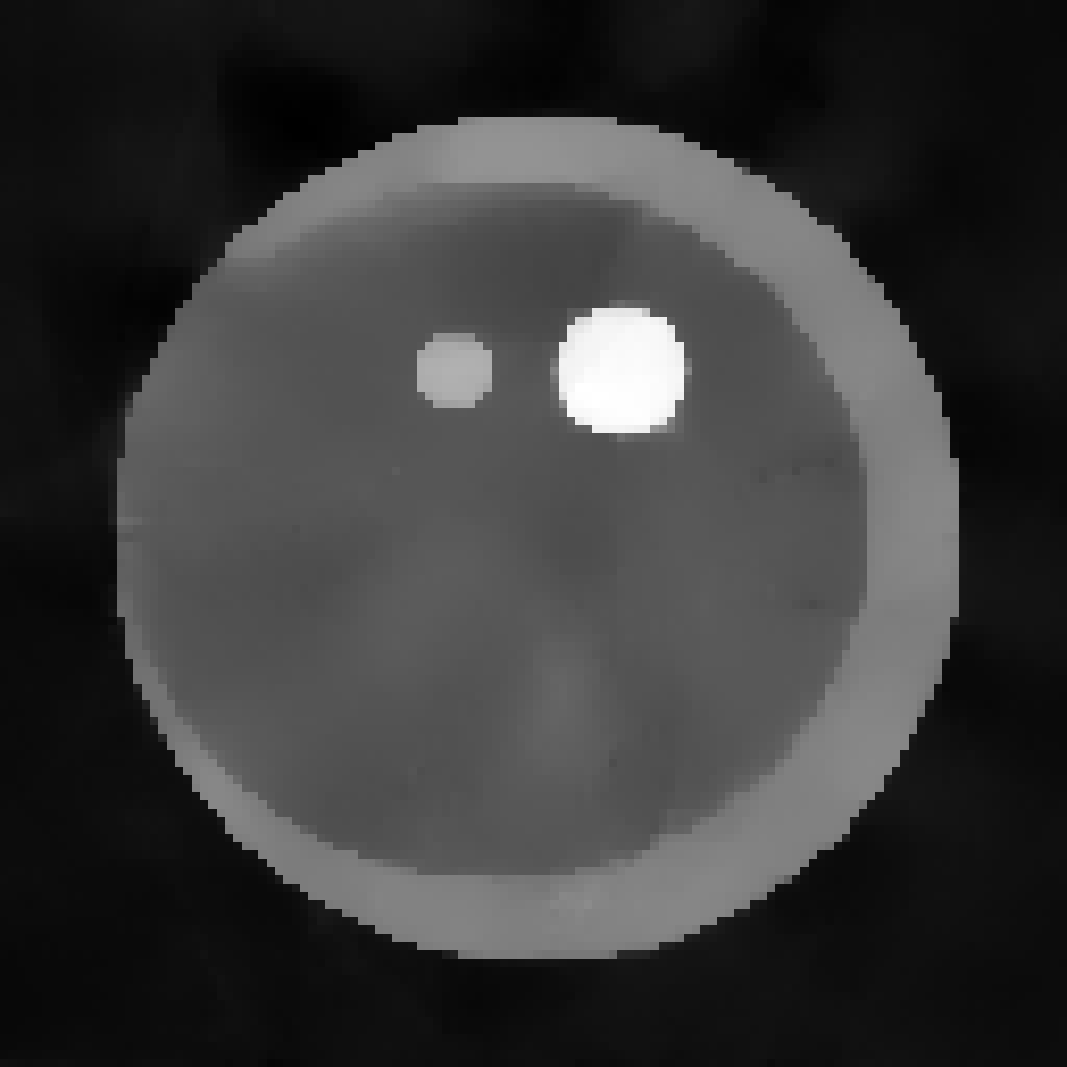}\\%(e)
	\end{minipage}
	\begin{minipage}{0.2\textwidth}
		\centering
		\includegraphics[width=\textwidth]{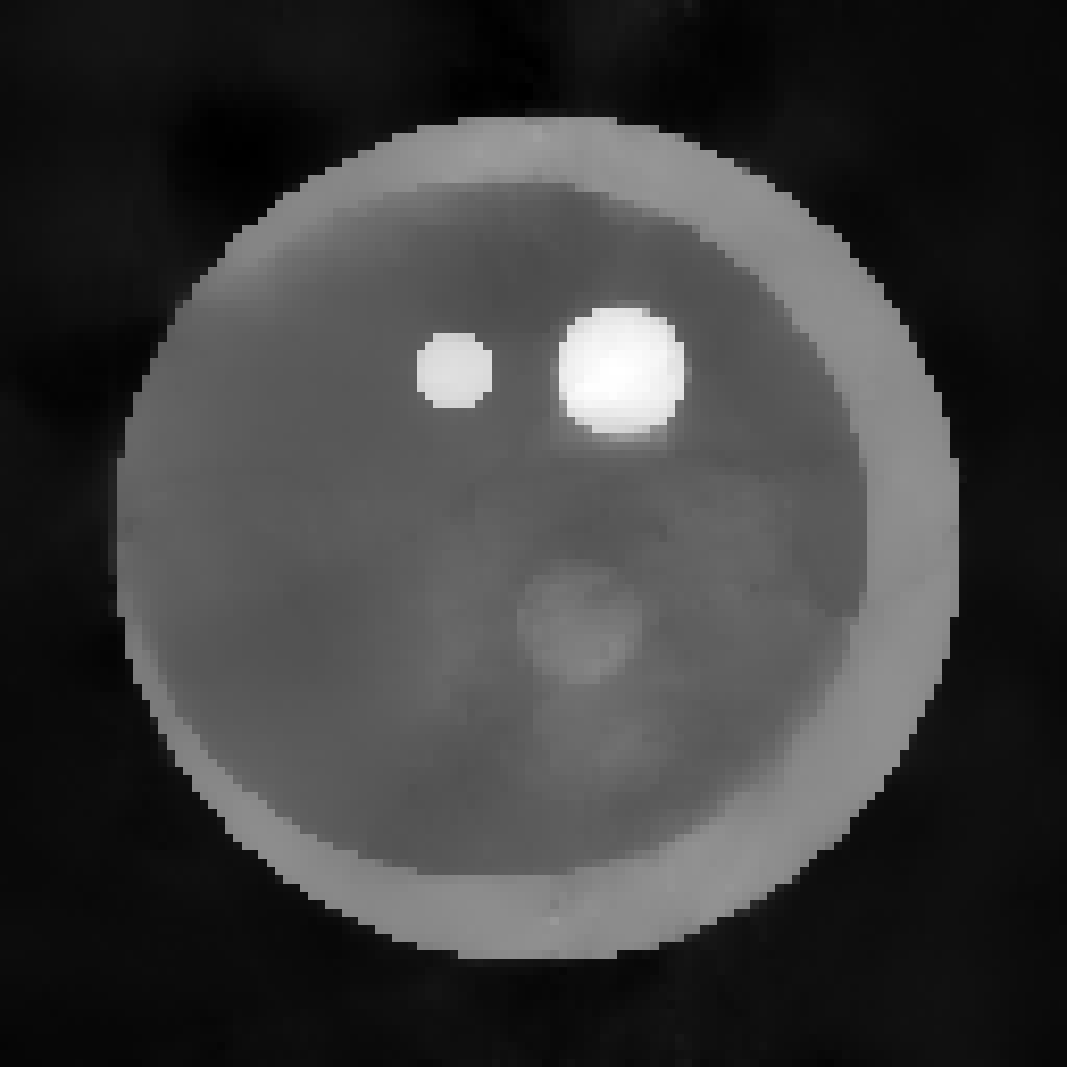}\\%(f)
	\end{minipage}
	\begin{minipage}{0.2\textwidth}
		\centering
		\includegraphics[width=\textwidth]{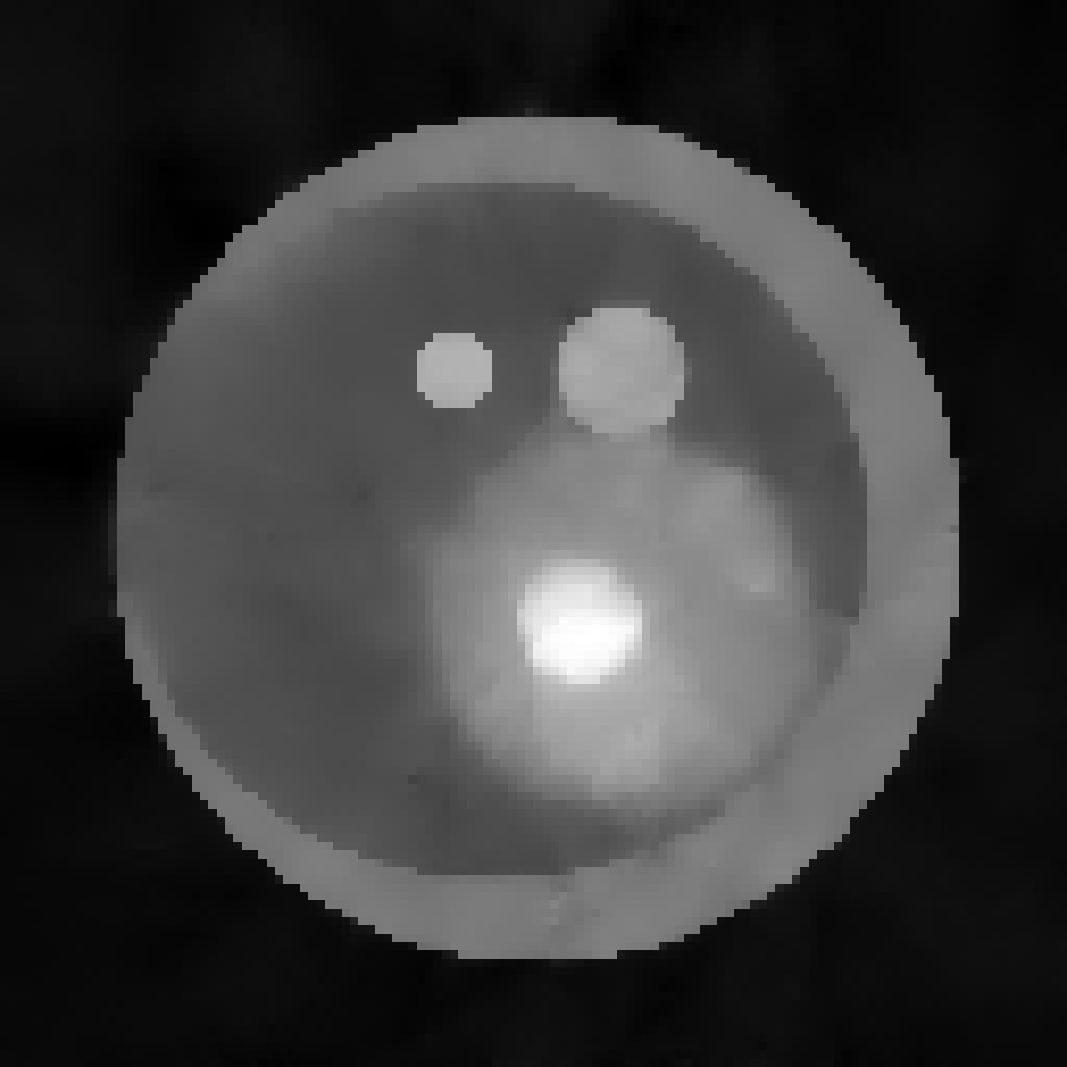}\\%(g)
	\end{minipage}
	\begin{minipage}{0.2\textwidth}
		\centering
		\includegraphics[width=\textwidth]{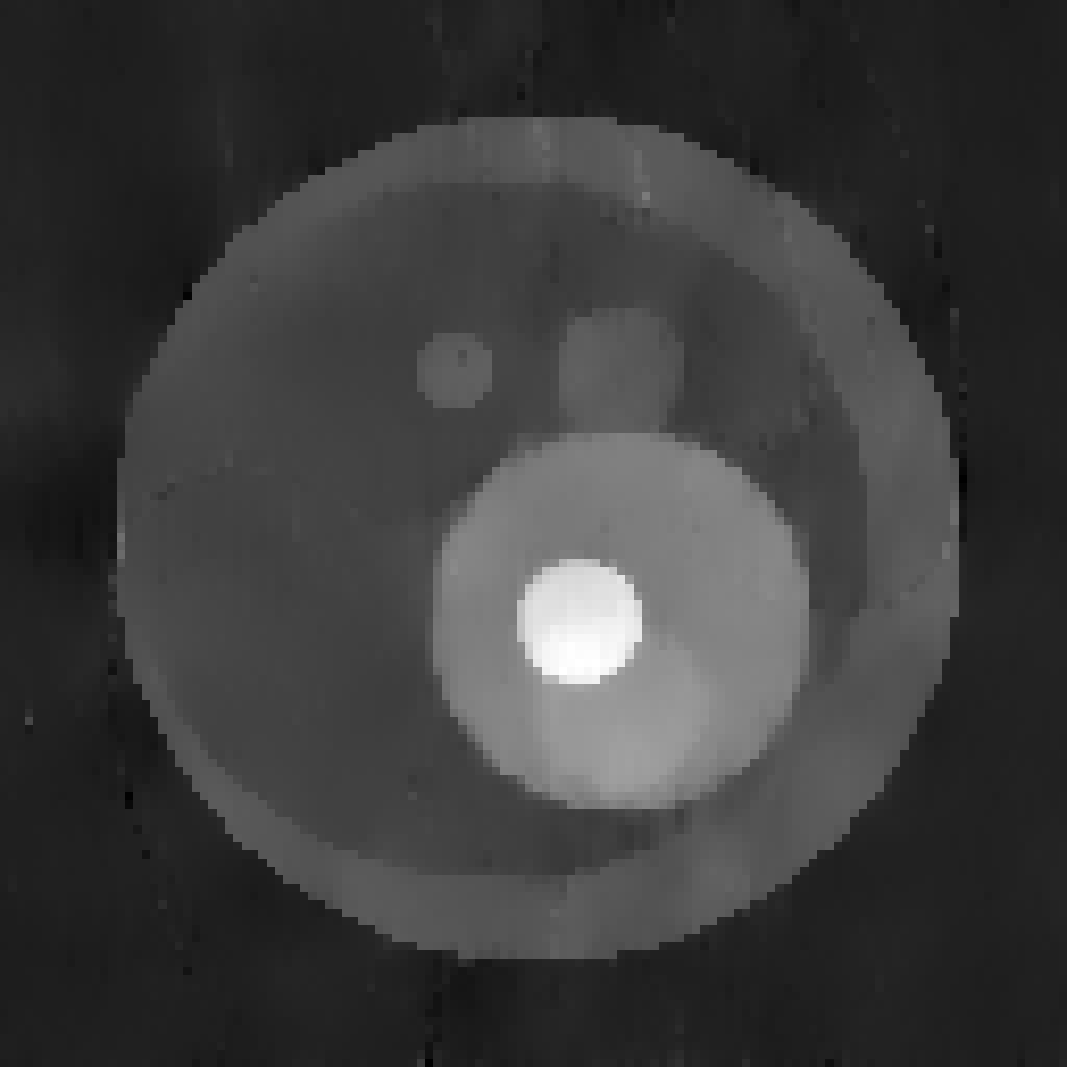}\\%(h)
	\end{minipage}
	
	\begin{minipage}{0.2\textwidth}
		\centering
		\includegraphics[width=\textwidth]{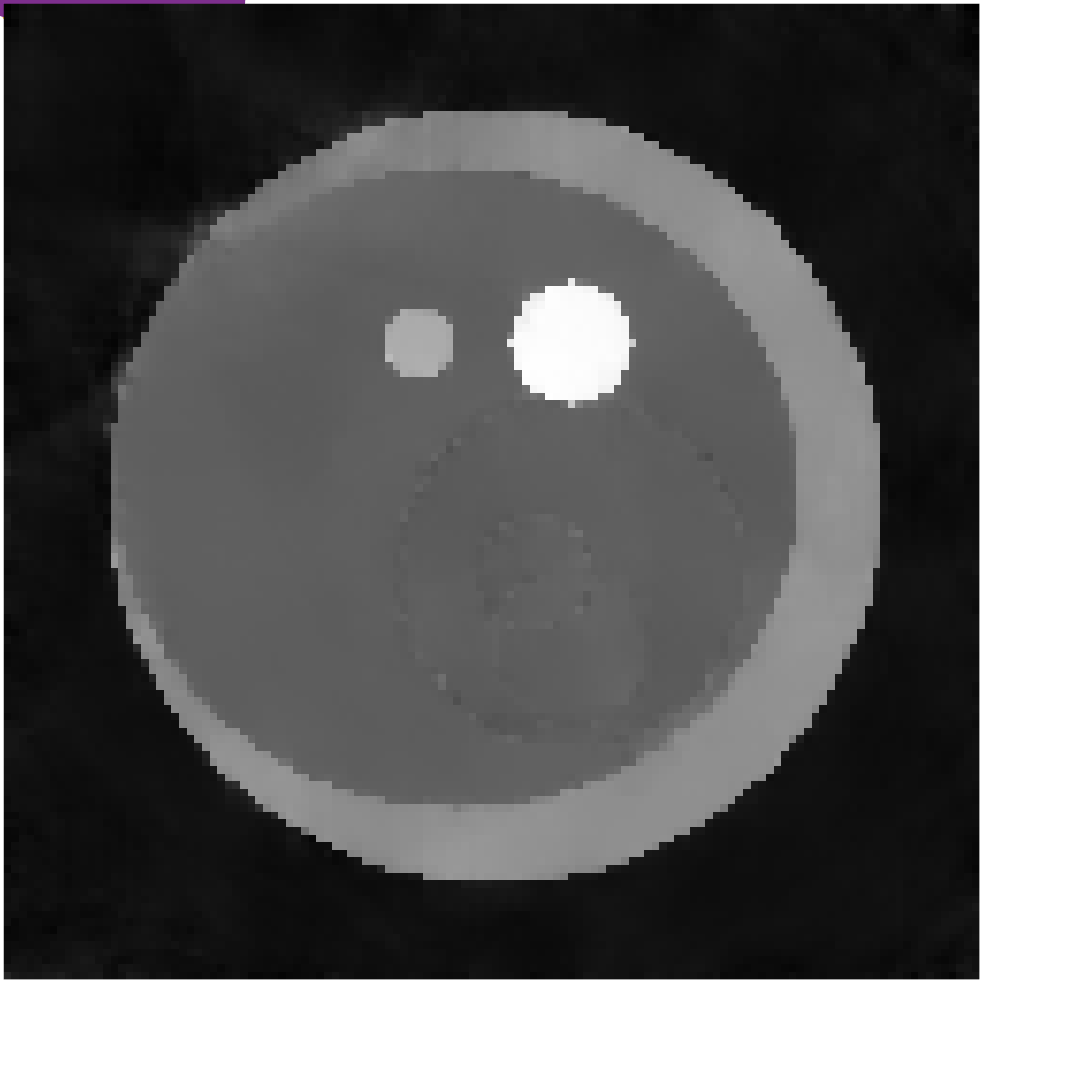}\\%(m)
	\end{minipage}
	\begin{minipage}{0.2\textwidth}
		\centering
		\includegraphics[width=\textwidth]{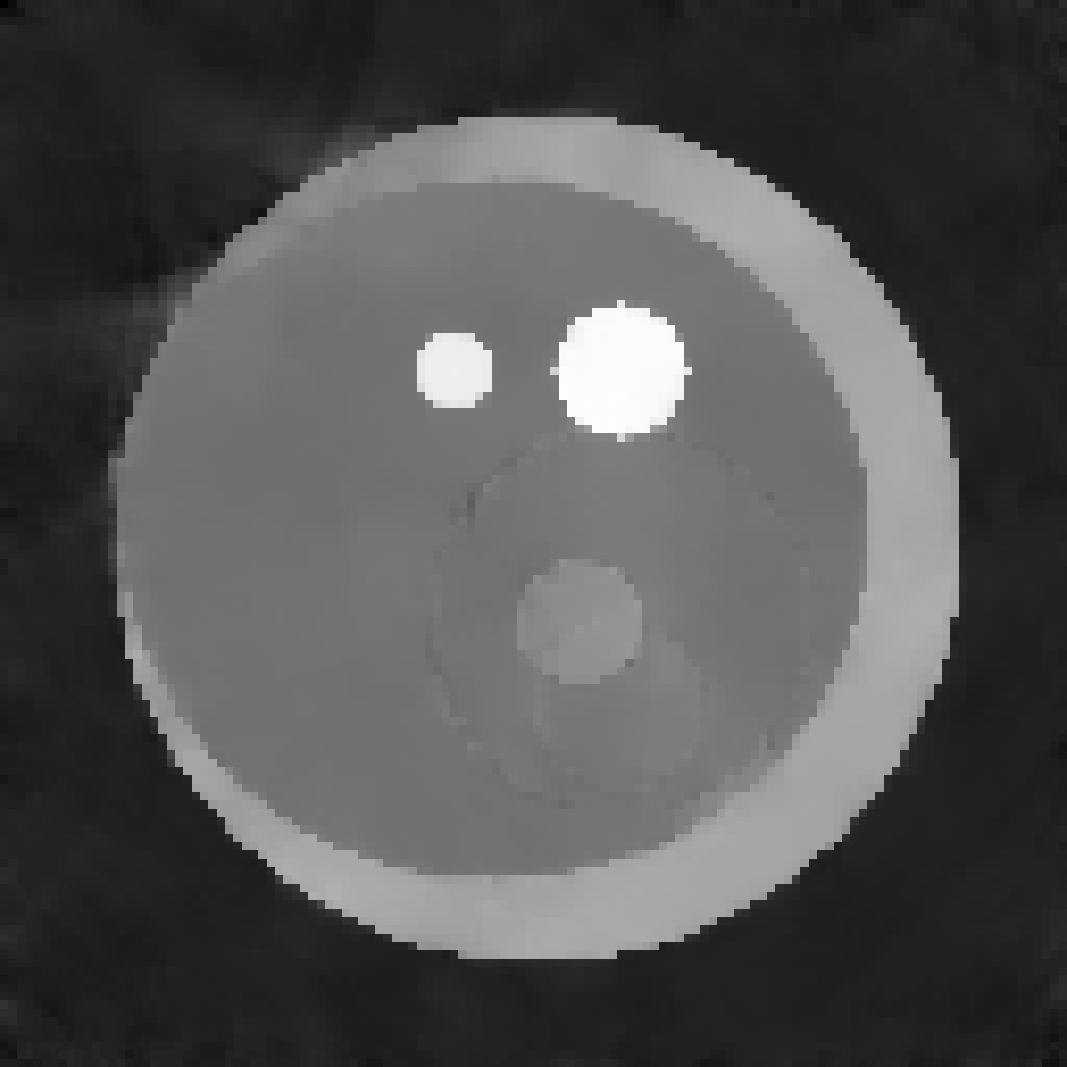}\\%(n)
	\end{minipage}
	\begin{minipage}{0.2\textwidth}
		\centering
		\includegraphics[width=\textwidth]{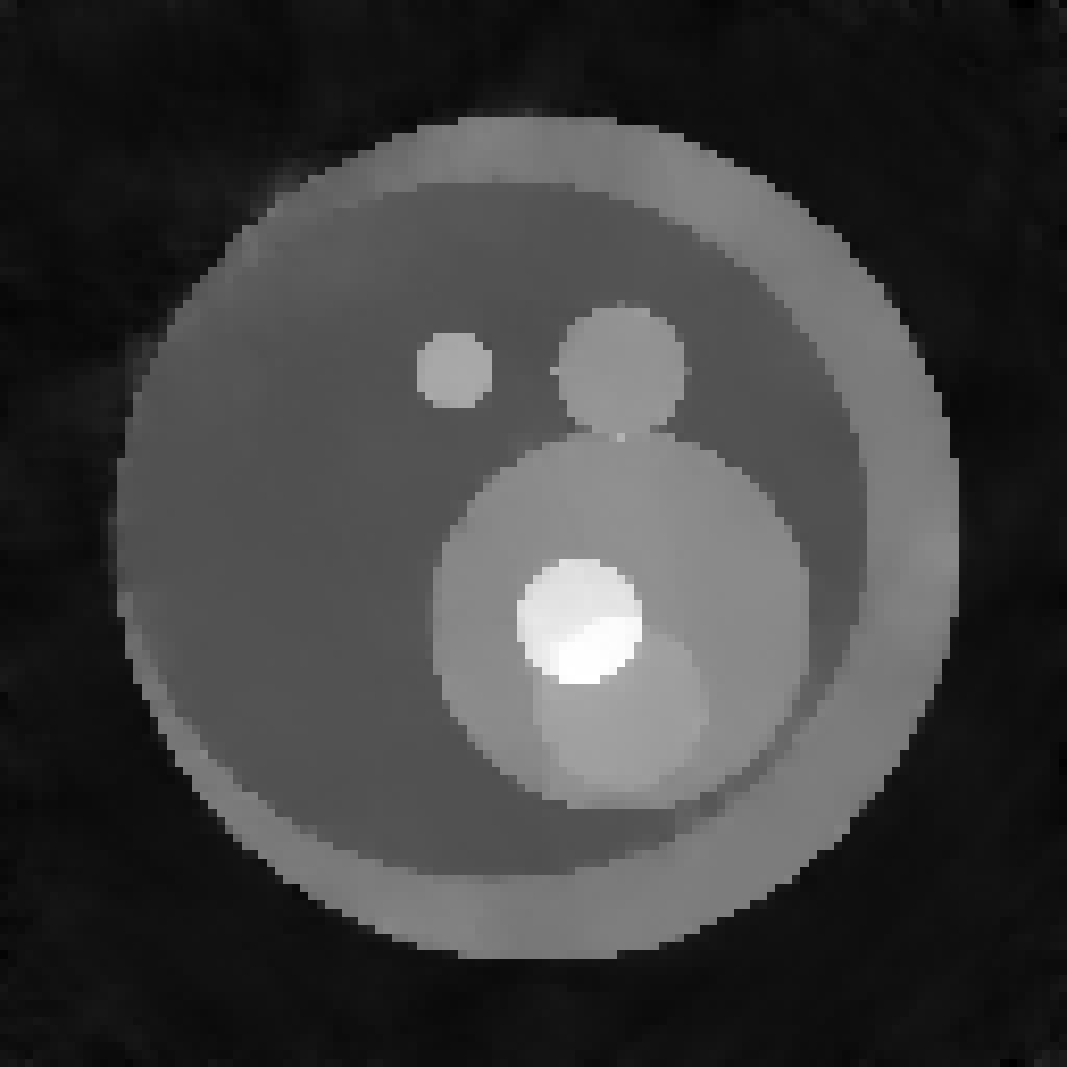}\\%(o)
	\end{minipage}
	\begin{minipage}{0.2\textwidth}
		\centering
		\includegraphics[width=\textwidth]{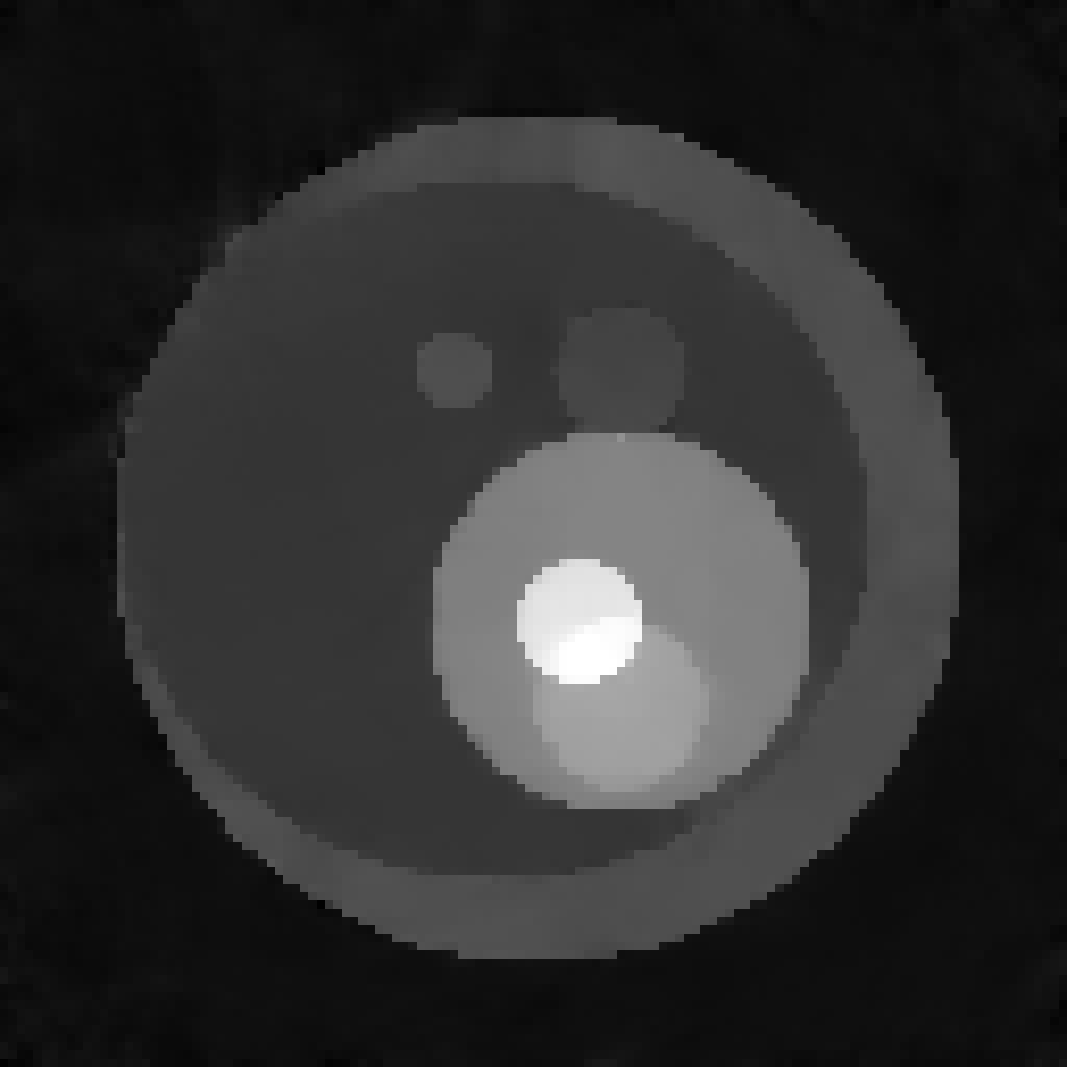}\\%(p)
	\end{minipage}
	\caption{PAT test problem: Panels in the first row show the reconstructions by solving the static inverse problem, panels in the second, third, and fourth rows show reconstructions with AnisoTV, Iso3DTV, and GS at time steps $t = 1, 10, 20, 30$ from left to right, respectively.}
	\label{fig: PATRec}
\end{figure}

We compute the average RRE$(\bu^{(k)},\bu_\mathrm{true})$ as well as the SSIM for both experimental setups as described in (i) and (ii) above and we report the results in Figure~\ref{Fig: RREPAT} when the discrepancy principle is satisfied for the first time. The number of the iterations and the corresponding regularization parameter $\lambda$ when the discrepancy principle is satisfied are reported in Table \ref{tab:ex2}. GS outperforms all the methods in this experimental setup, followed by AnisoTV, illustrated in both RRE and SSIM in Figure \ref{fig: PATRec}. Notice here that Iso3DTV is the least accurate method. We report the reconstructions and the quantitative figures of merit at the iteration 150 since the discrepancy principle is not satisfied. We remark that the similar reconstruction quality is manifested from the Aniso3DTV as well that we do not report here. 

In Figure \ref{fig: PATRec} we report the reconstructions at times steps $t = 1,10, 20, 30$ from left to right respectively. Different rows correspond to reconstructions with different methods.
The first row shows the reconstructions obtained by solving the static inverse problem~\eqref{eq: static} where we observe that even though the method is able to provide the locations of the inclusions, the detailed information of the inclusions are missing. The second row shows the reconstructions with Iso3DTV, where certainly the artifacts around the circular inclusions are present and the background is perturbed as well. Improved reconstructions are observed in the third and the fourth rows of Figure \ref{fig: PATRec}, obtained by AnisoTV and GS respectively.

\begin{table}[ht!]
    \centering
    \caption{Dynamic photoacoustic tomography (PAT) example: The number of iterations when the discrepancy principle is satisfied for the first time and the regularization parameters at those iterations for AnisoTV, Iso3DTV, and GS. Here superscript $(*)$ means that the maximum number of iterations were reached before convergence is satisfied using the discrepancy principle.}
   \label{tab:ex2}
    \begin{tabular}{l|c|c|c}
       -  &  AnisoTV&  Iso3DTV &  GS \\ \hline
      Iters   & 60 & 150$^{(*)}$ & 75\\ 
      $\lambda$ & 0.013 &  0.016 & 0.02 
    \end{tabular}
\end{table}

\subsection{Example 3: Dynamic X-Ray Tomography- 3D Emoji Data}
In this example we test our methods on real data of an emoji phantom measured at the
University of Helsinki \cite{meaney2018tomographic}.

\begin{figure}[h!]
\begin{centering}
    	\begin{minipage}{0.2\textwidth}
		\includegraphics[width=\textwidth]{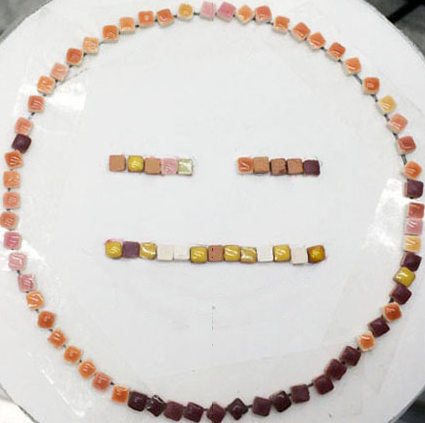}
	\end{minipage}
	\begin{minipage}{0.2\textwidth}
		\includegraphics[width=\textwidth]{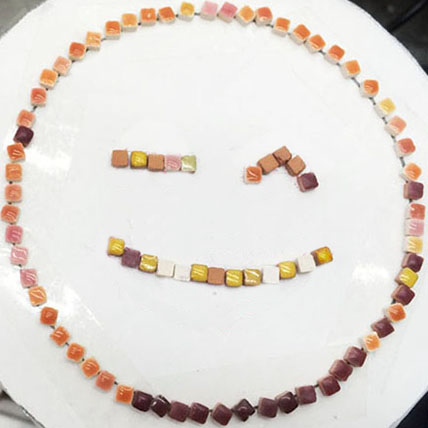}
	\end{minipage}
	\begin{minipage}{0.2\textwidth}
		\includegraphics[width=\textwidth]{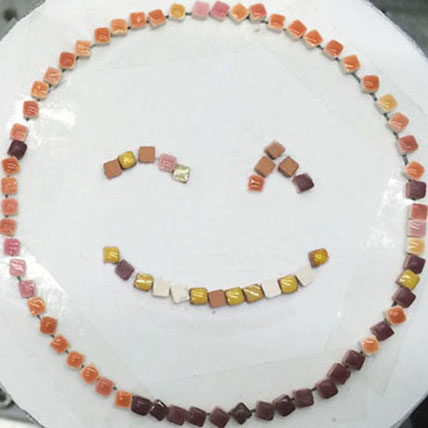}
	\end{minipage}
	\begin{minipage}{0.2\textwidth}
		\includegraphics[width=\textwidth]{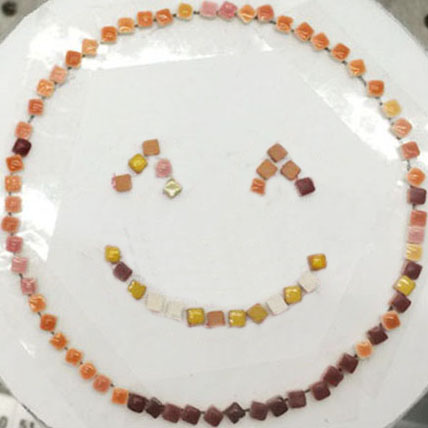}
	\end{minipage}\\
	\begin{minipage}{0.2\textwidth}
		\includegraphics[width=\textwidth,angle =-90]{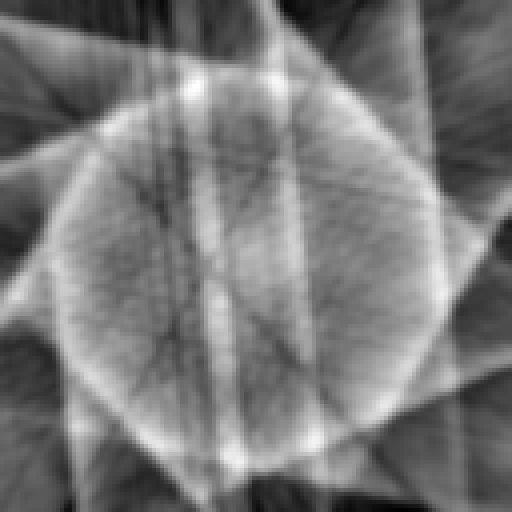}
	\end{minipage}
	\begin{minipage}{0.2\textwidth}
		\includegraphics[width=\textwidth,angle =-90]{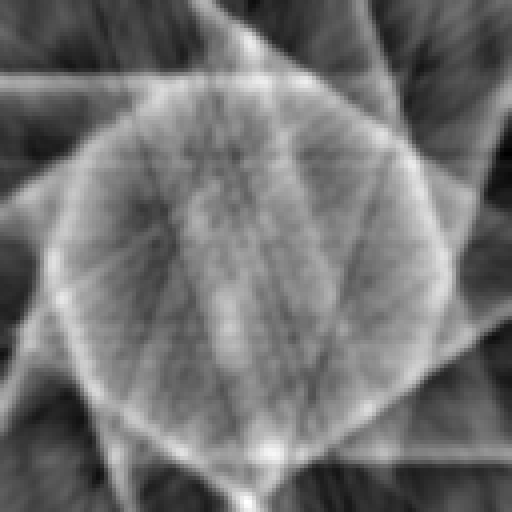}
	\end{minipage}
	\begin{minipage}{0.2\textwidth}
		\includegraphics[width=\textwidth,angle =-90]{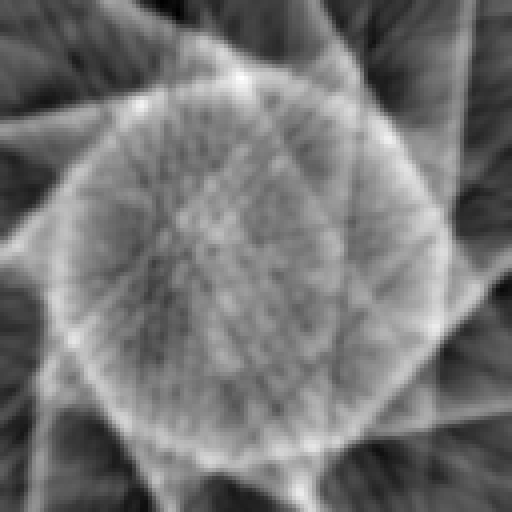}
	\end{minipage}
	\begin{minipage}{0.2\textwidth}
		\includegraphics[width=\textwidth,angle =-90]{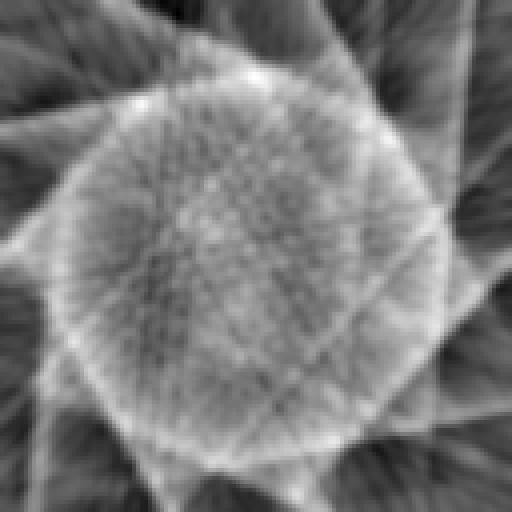}
	\end{minipage}\\
	\begin{minipage}{0.2\textwidth}
		\includegraphics[width=\textwidth]{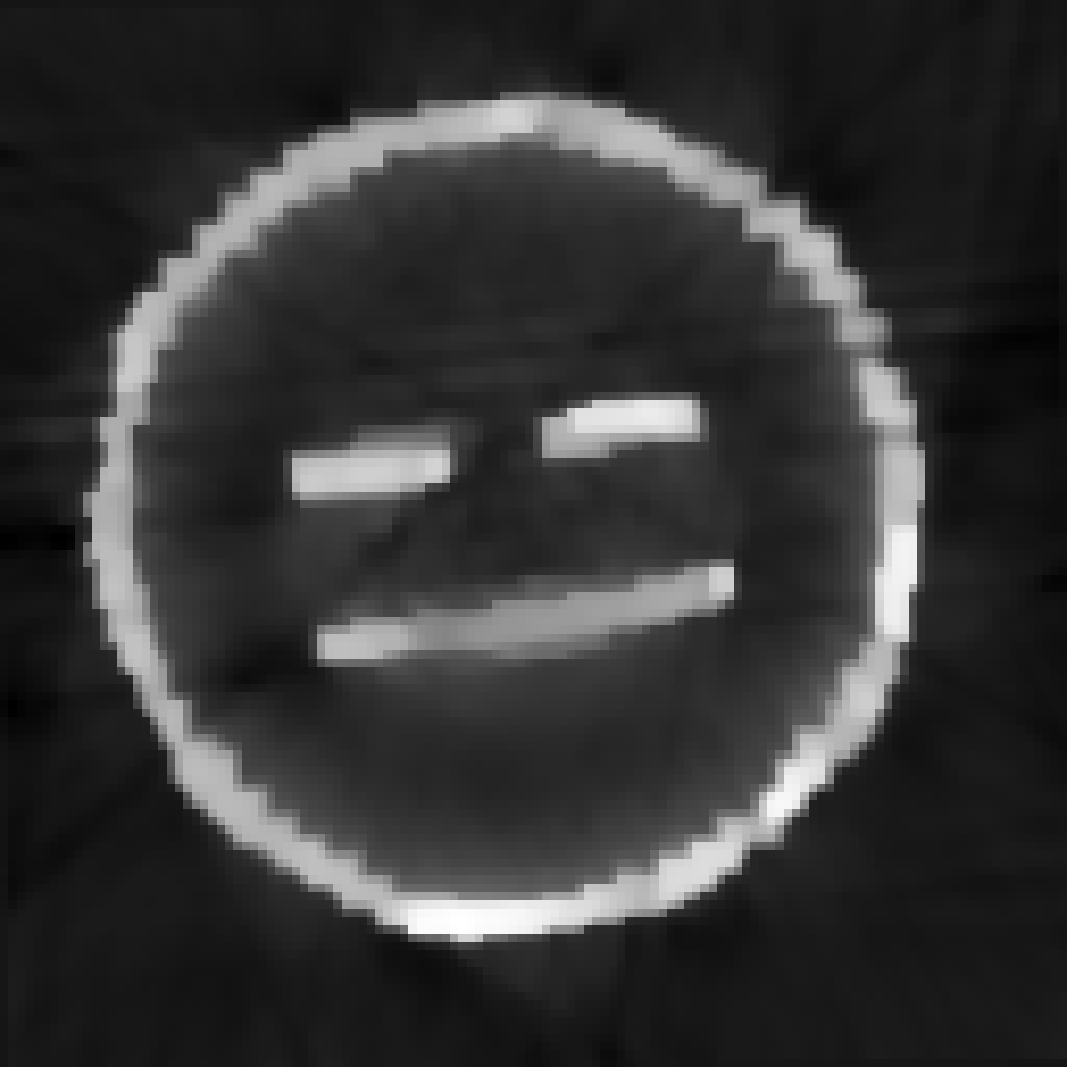}
	\end{minipage}
	\begin{minipage}{0.2\textwidth}
		\includegraphics[width=\textwidth]{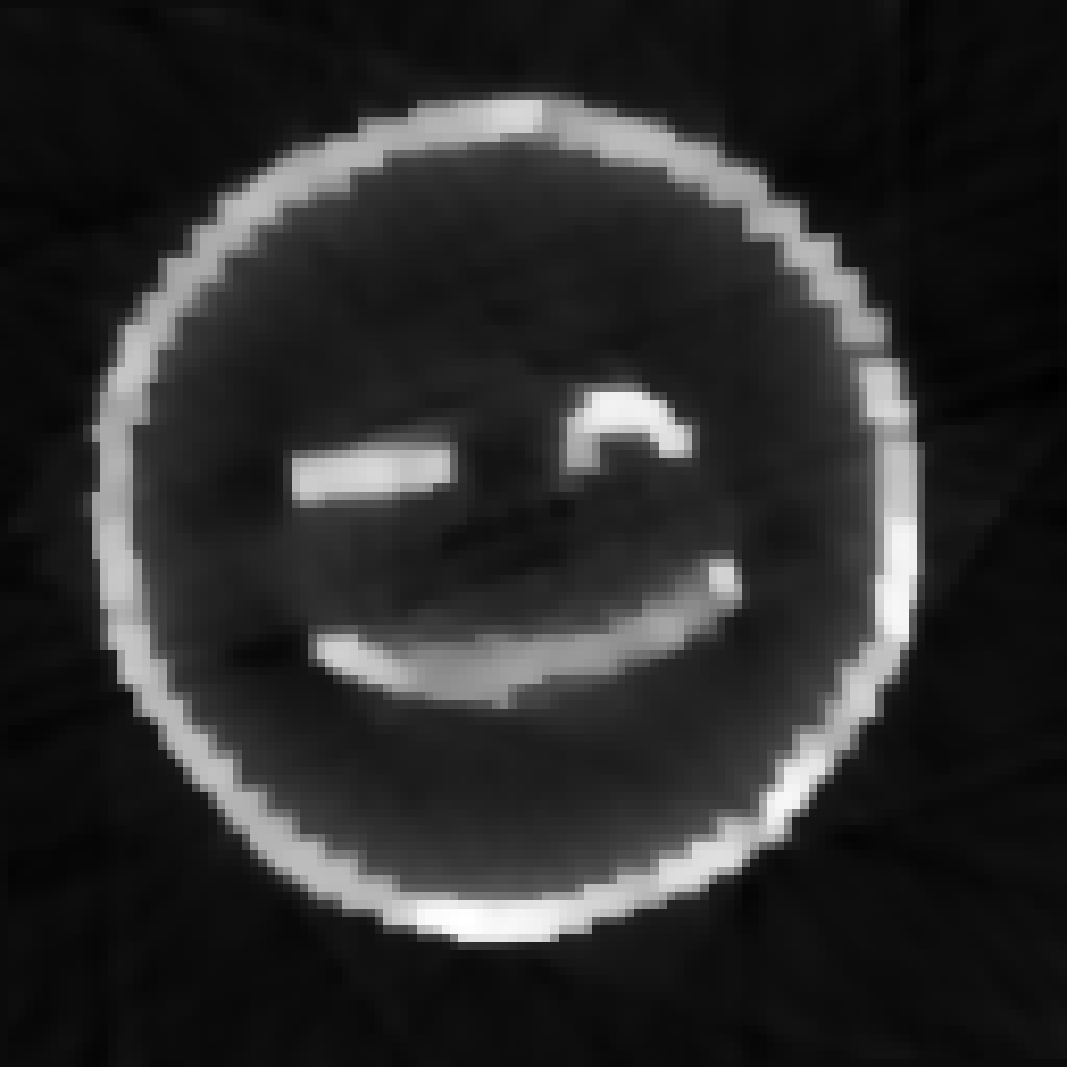}
	\end{minipage}
	\begin{minipage}{0.2\textwidth}
		\includegraphics[width=\textwidth]{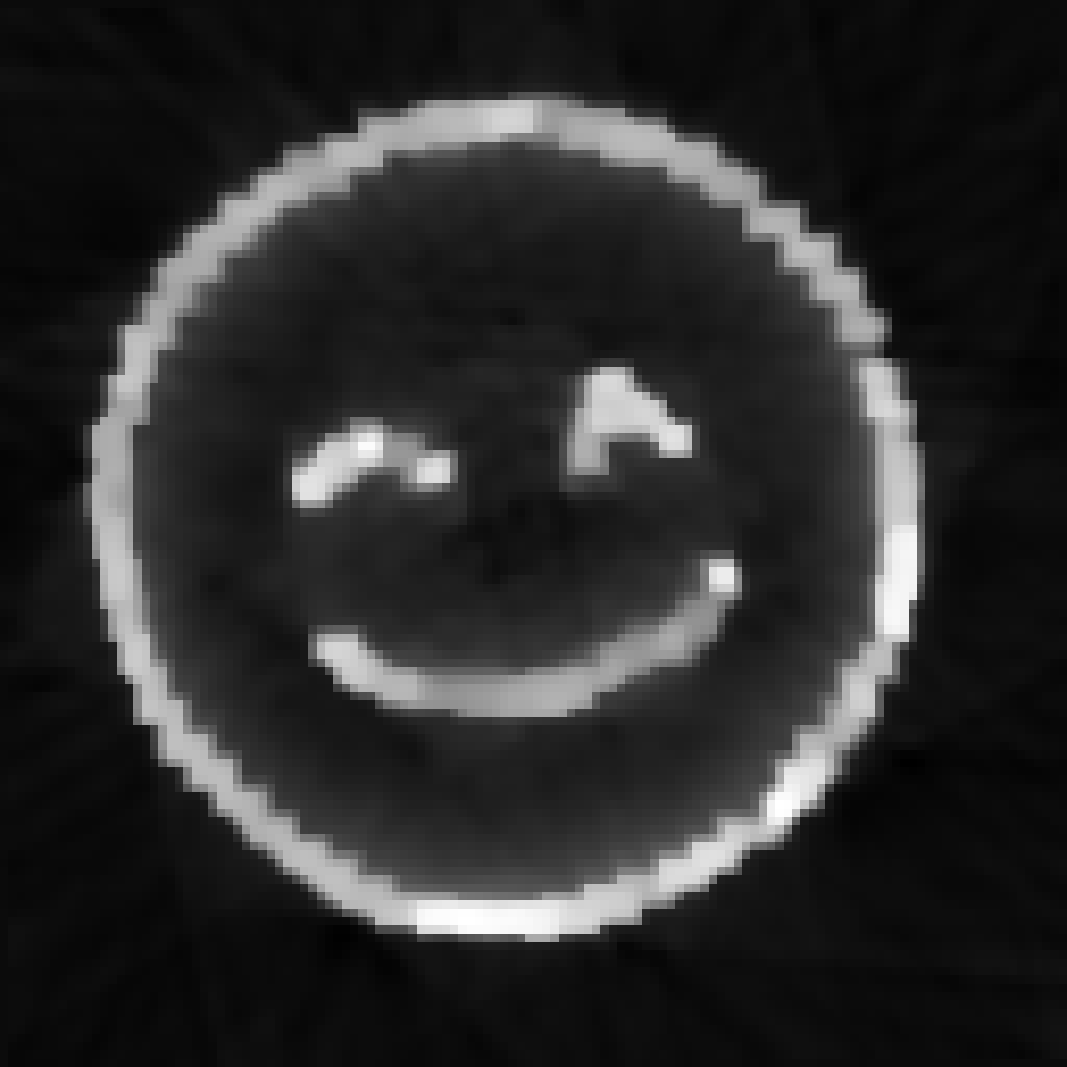}
	\end{minipage}
	\begin{minipage}{0.2\textwidth}
		\includegraphics[width=\textwidth]{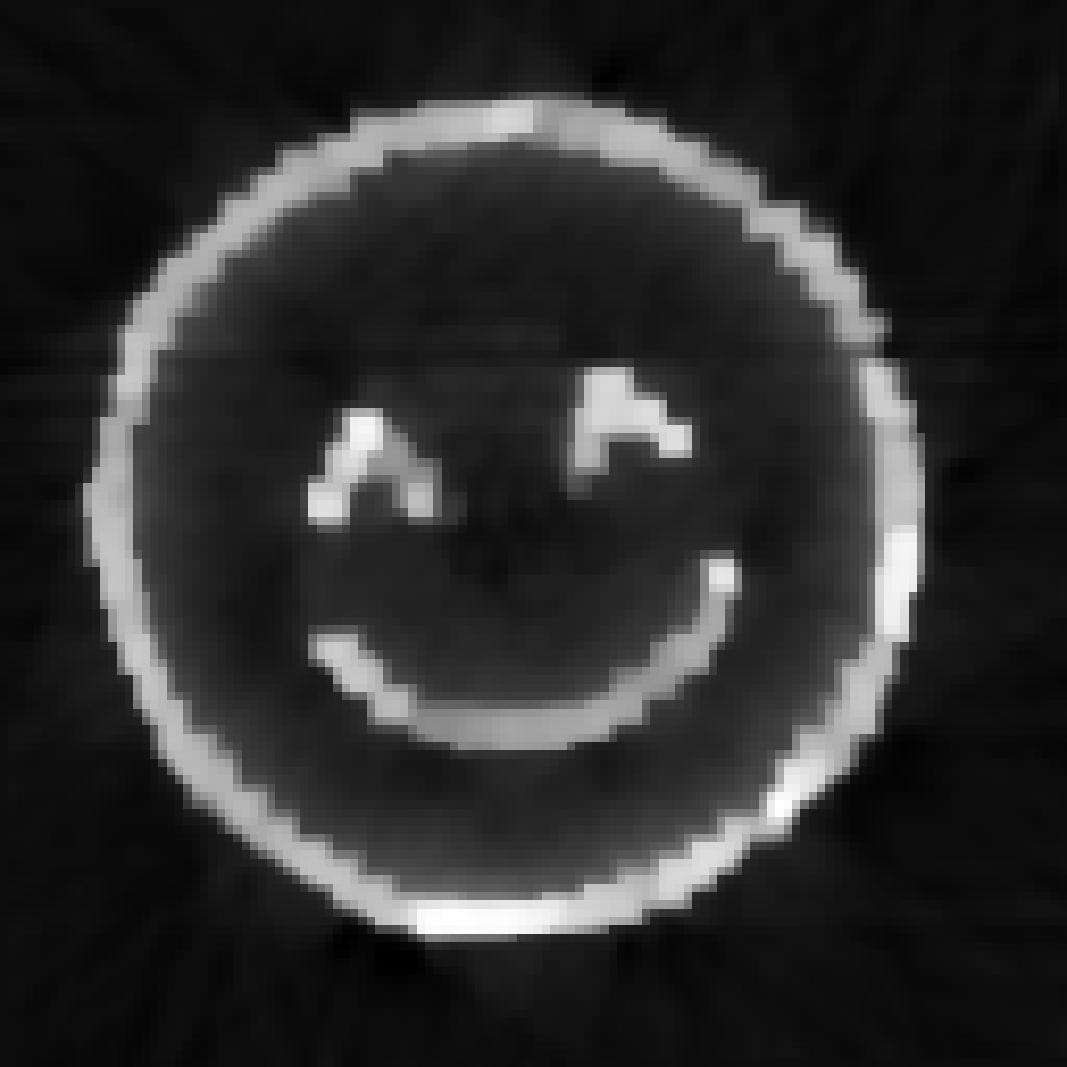}
	\end{minipage}\\
	\begin{minipage}{0.2\textwidth}
		\includegraphics[width=\textwidth]{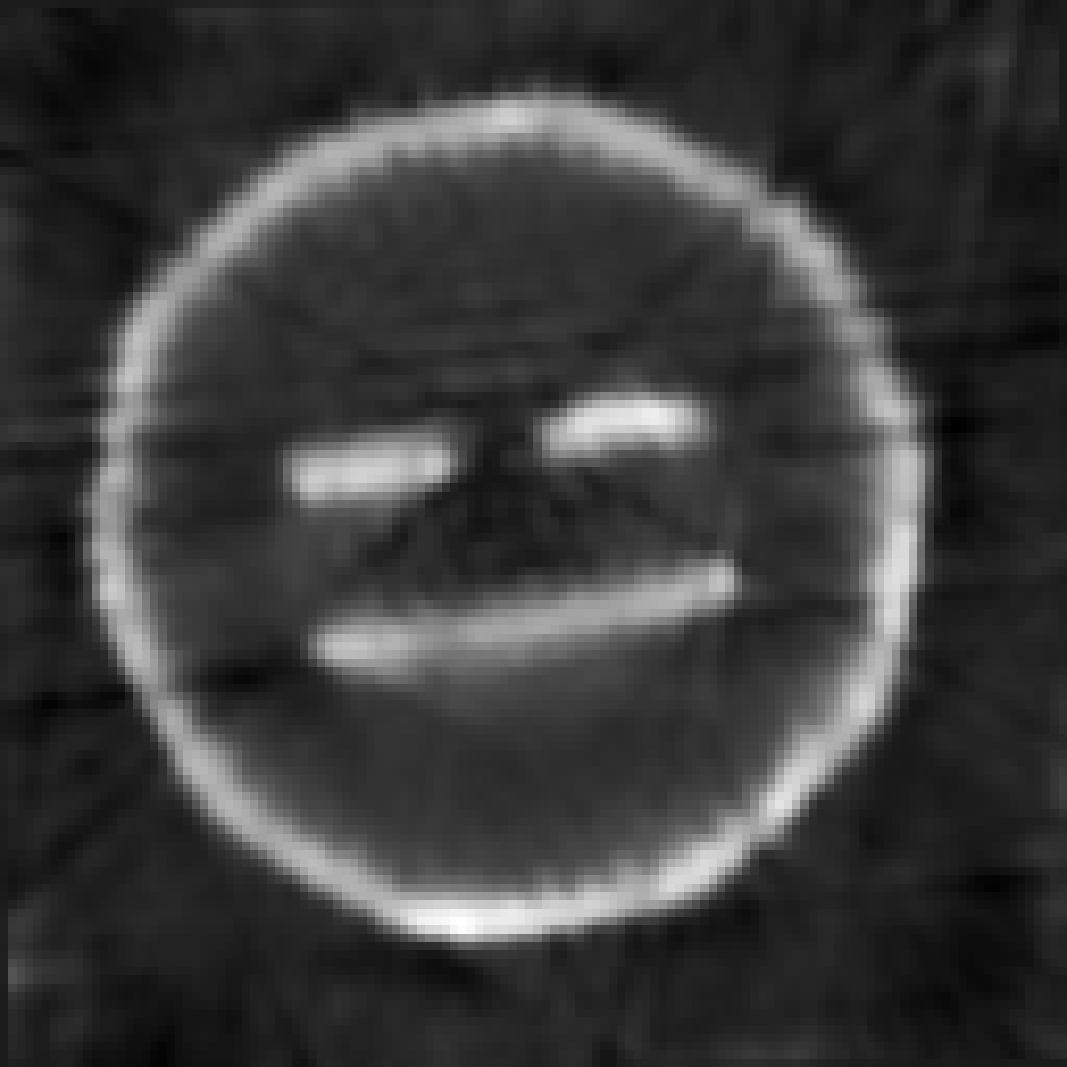}
	\end{minipage}
	\begin{minipage}{0.2\textwidth}
		\includegraphics[width=\textwidth]{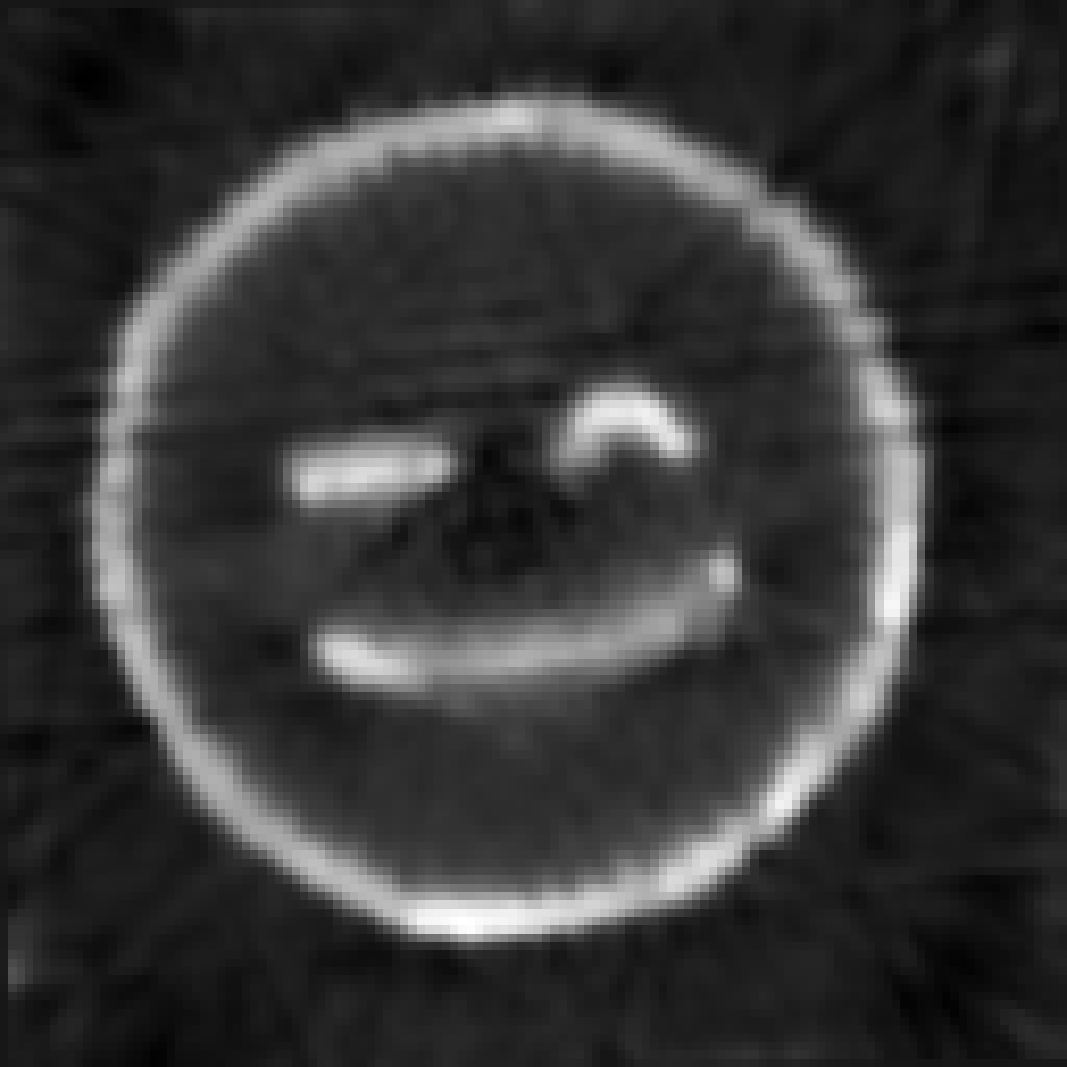}
	\end{minipage}
	\begin{minipage}{0.2\textwidth}
		\includegraphics[width=\textwidth]{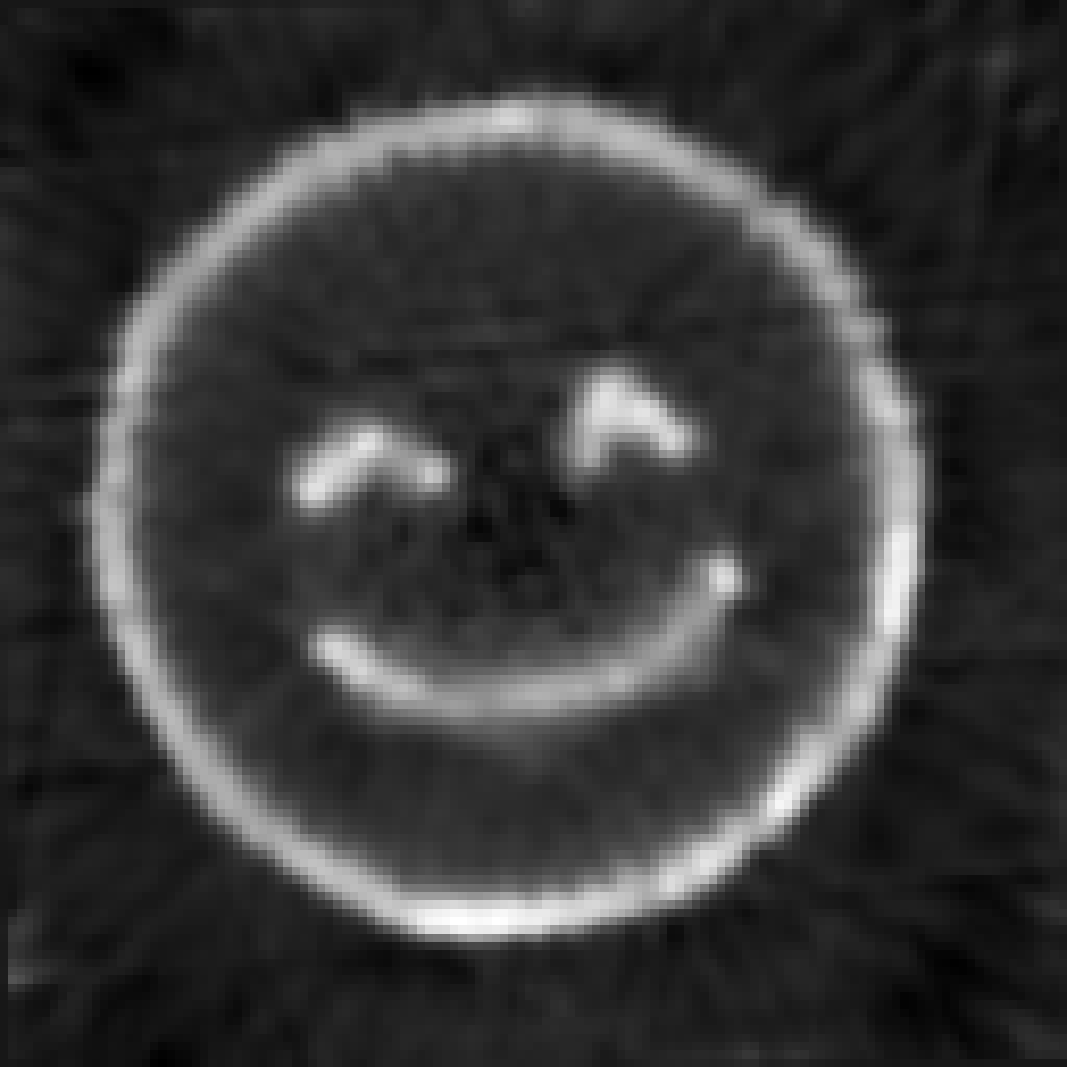}
	\end{minipage}
	\begin{minipage}{0.2\textwidth}
		\includegraphics[width=\textwidth]{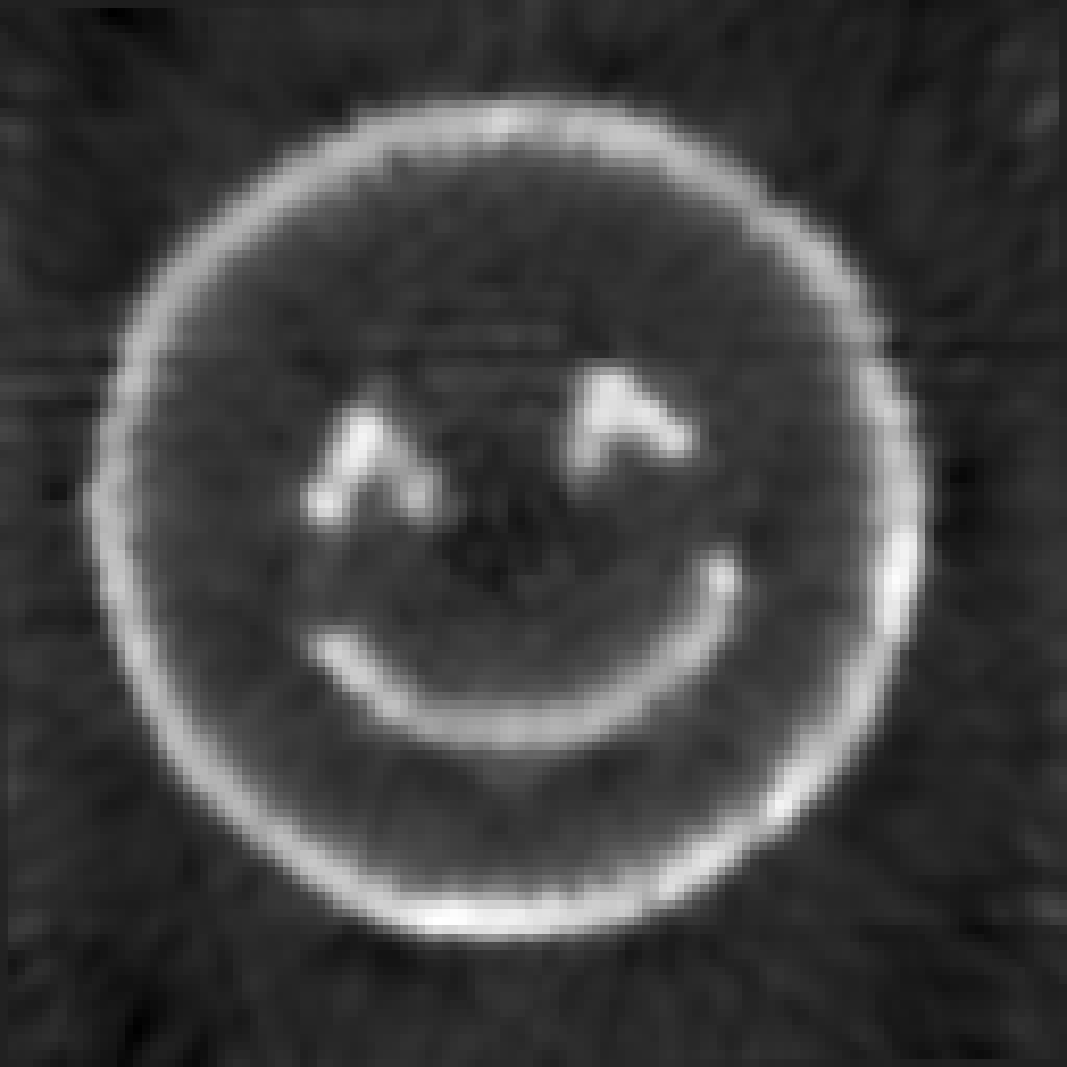}
	\end{minipage}\\
\end{centering}
	\caption{Reconstruction results for the emoji test problem with $n_a =10$. The rows represent (from top to bottom): the original images, the reconstructions when images are considered independently, the reconstructions by AnisoTV,  the reconstructions by Iso3DTV, at time steps $t = 2, 10, 18, 31$ (from left to right).}
%	$n_t = 1,5, 9, 15$
	\label{Fig: smile_10angles}
\end{figure}

In particular, we consider the file \verb|DataDynamic_128x30.mat|.
The available data represents $n_t=33$ time steps of a series of the X-ray sinogram of  ``emojis'' made of small ceramic stones obtained by shining 217 projections from $n_a =30$ angles. 
We are interested in reconstructing a sequence of images $\bU^{(t)}$, $t=1,2,\dots, 33$ is $128 \times 128$, of size $n_h\times n_v$, where $n_h = n_v = 128$, from low-dose observations measured from a limited number of angles $n_a$. Hence $\bu \in \R^{540,672}$, with $\bu=\left[(\bu^{(1)})^T,\\(\bu^{(2)})^T,\\\dots\\,(\bu^{(33)})^T\right]^T$
representing the dynamic sequence of the emoji changing from an expressionless face with closed eyes and a straight mouth to a face with smiling eyes and mouth where the circular shape does not change. See the first row of Figure \ref{Fig: smile_10angles} for a sample of 4 images at time steps $t = 6,14,20, 26$.
%$n_t = 1,5,9,15$. 
The low-dose available observations can be modeled
by the measurement matrix $\bF$ which describes the forward model of the Radon transform that represents line integrals. In this case, we have a block matrix $\bF$ as in \eqref{eq: blockF} with 33 blocks. Although the ground truth is not available to compare the qualitative
results, we can observe the visual results from different numbers of projections as illustrated in Case 1 and Case 2 below. 
We do not know which kind of noise contaminates the data if any, but we artificially add $1\%$ of white Gaussian noise in the available sinogram.
The experiment that we perform here focuses on the visual inspection of the reconstructions from different number of angels $n_a = 10$ and $n_a = 30$, highlighting the effect of the number of the projection angles and also the  visual differences in the reconstruction when static sub-problems \eqref{eq: static}
are solved independently and when the dynamic inverse problem %$\bF\bu + \be = \bd$ 
\ref{eqn:dynamic}
is solved. 
\paragraph{Case 1: Consider $n_a =10$ projection angles}
First, we limit the number of angles $n_a$ to 10 from the dataset \verb|DataDynamic_128x30.mat|.
In this way we generate underdetermined problems $\bA^{(t)}\bu^{(t)} + \be^{(t)} = \bd^{(t)}$, $t = 1,2,\dots, 33$ where $\bA^{(t)} \in \R^{2,170\times 16,384}$. Therefore $\bF \in \R^{71,610 \times 540,672}$ and the measurement vector $\bd \in \R^{71,610}$ contains the measured sinograms $\bd^{(t)} \in \R^{2,170}$ obtained from 217 projections around 10 equidistant angles. 

Figure~\ref{Fig: smile_10angles} displays some reconstructions (see also supplementary materials for an animation). 

\begin{figure}[h!]
\centering
	\begin{minipage}{0.2\textwidth}
		\includegraphics[width=\textwidth]{Figures/smile_true_1.png}
	\end{minipage}
	\begin{minipage}{0.2\textwidth}
		\includegraphics[width=\textwidth]{Figures/smile_true_5.png}
	\end{minipage}
	\begin{minipage}{0.2\textwidth}
		\includegraphics[width=\textwidth]{Figures/smile_true_9.png}
	\end{minipage}
	\begin{minipage}{0.2\textwidth}
		\includegraphics[width=\textwidth]{Figures/smile_true_15.png}
	\end{minipage}\\
	\begin{minipage}{0.2\textwidth}
		\includegraphics[width=\textwidth]{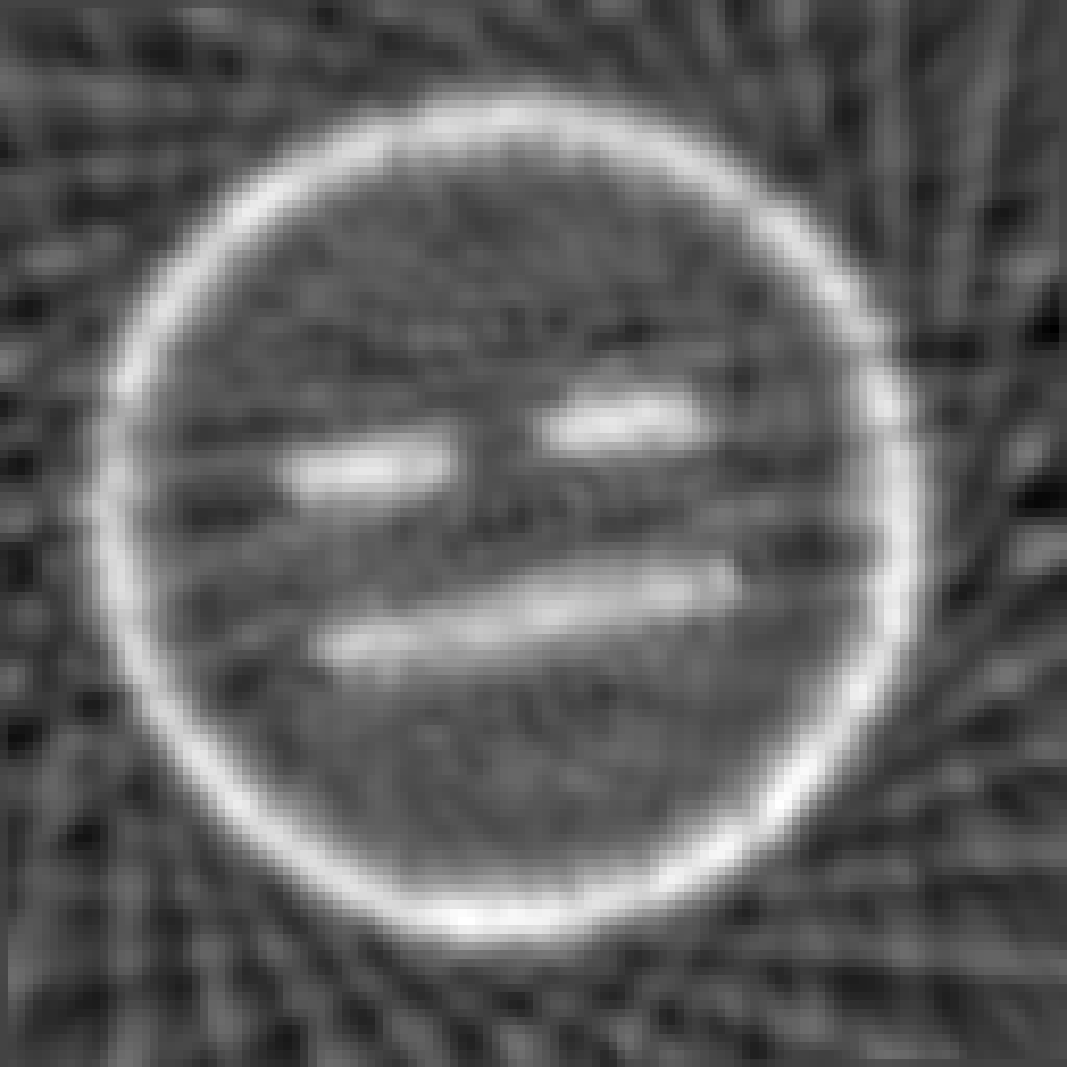}
	\end{minipage}
	\begin{minipage}{0.2\textwidth}
		\includegraphics[width=\textwidth]{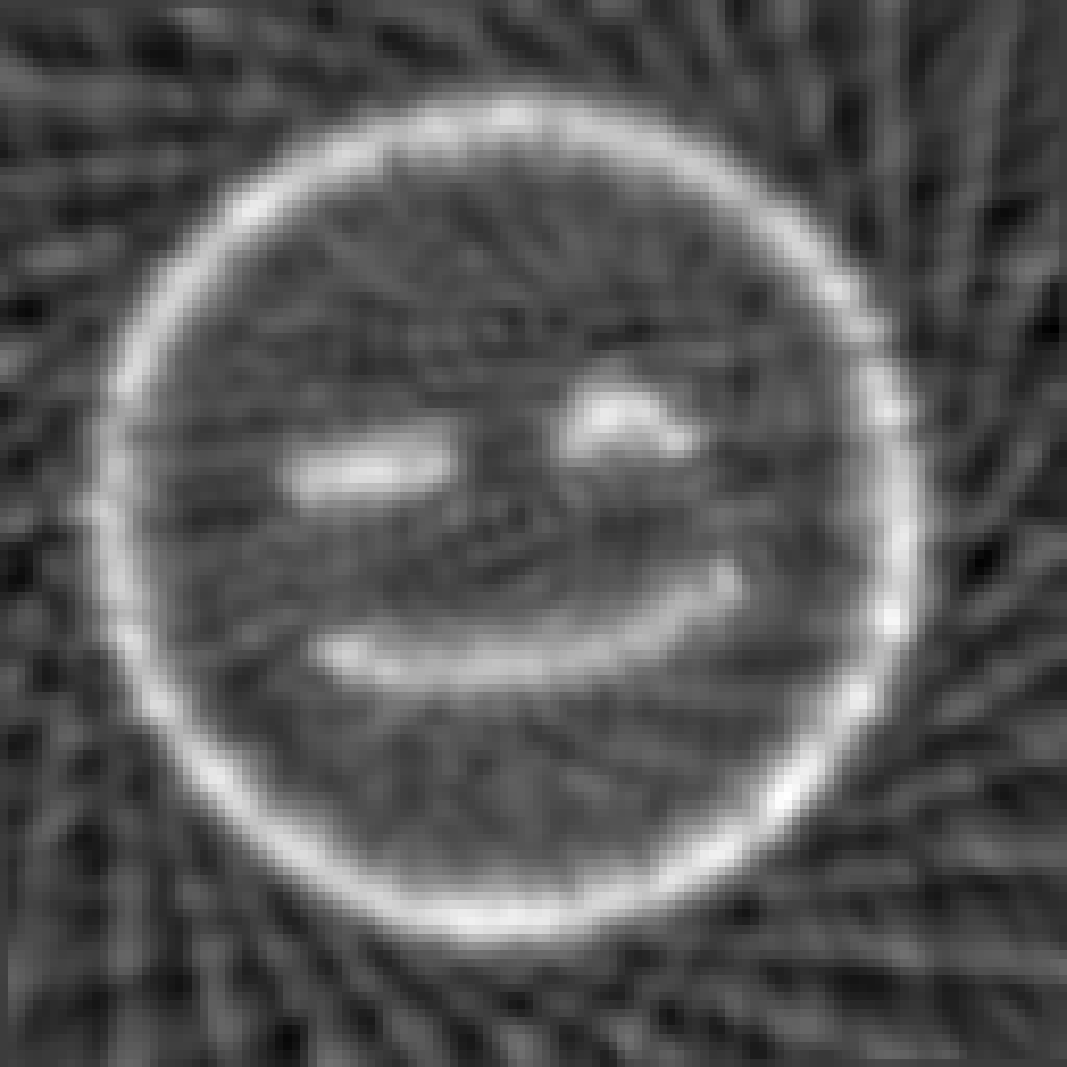}
	\end{minipage}
	\begin{minipage}{0.2\textwidth}
		\includegraphics[width=\textwidth]{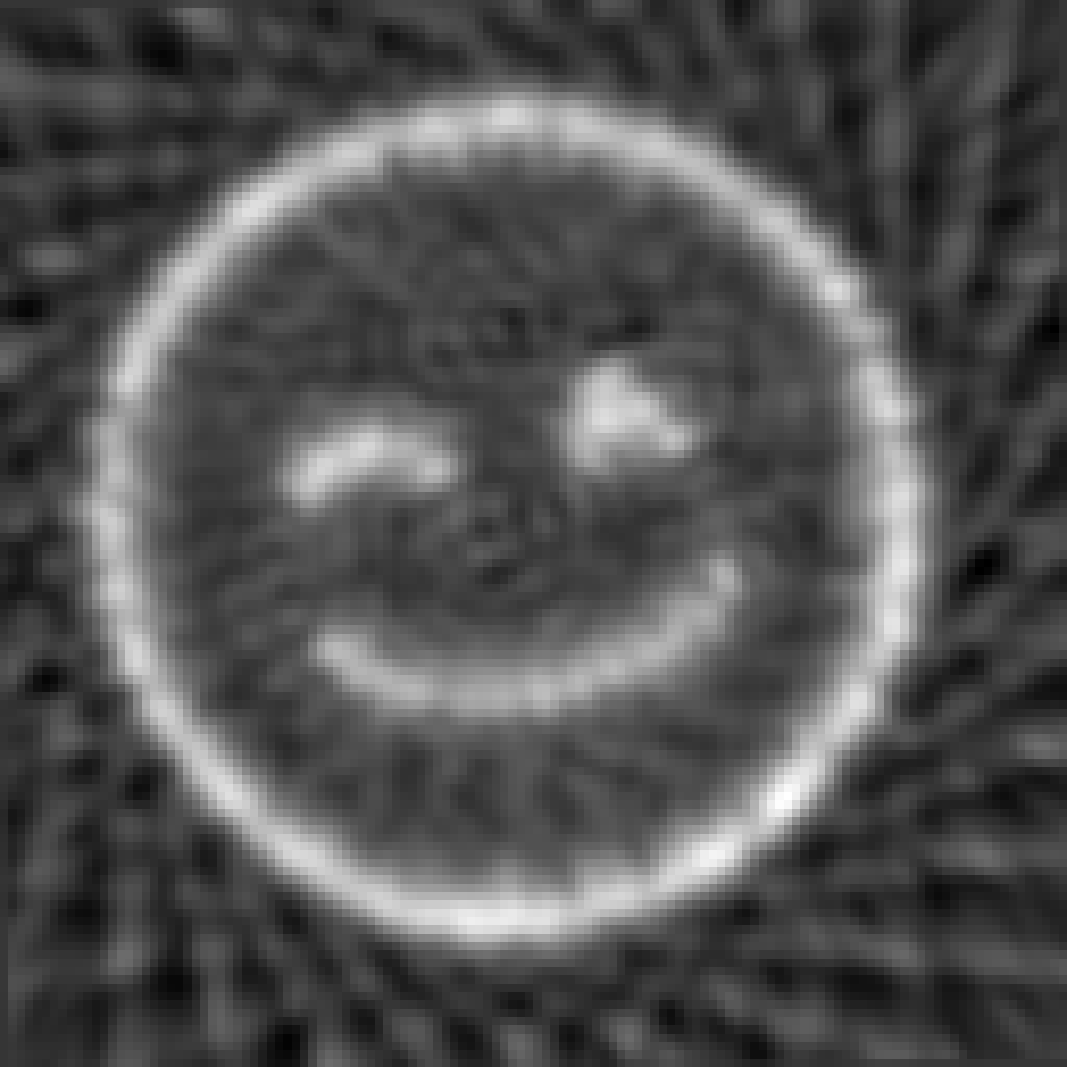}
	\end{minipage}
	\begin{minipage}{0.2\textwidth}
		\includegraphics[width=\textwidth]{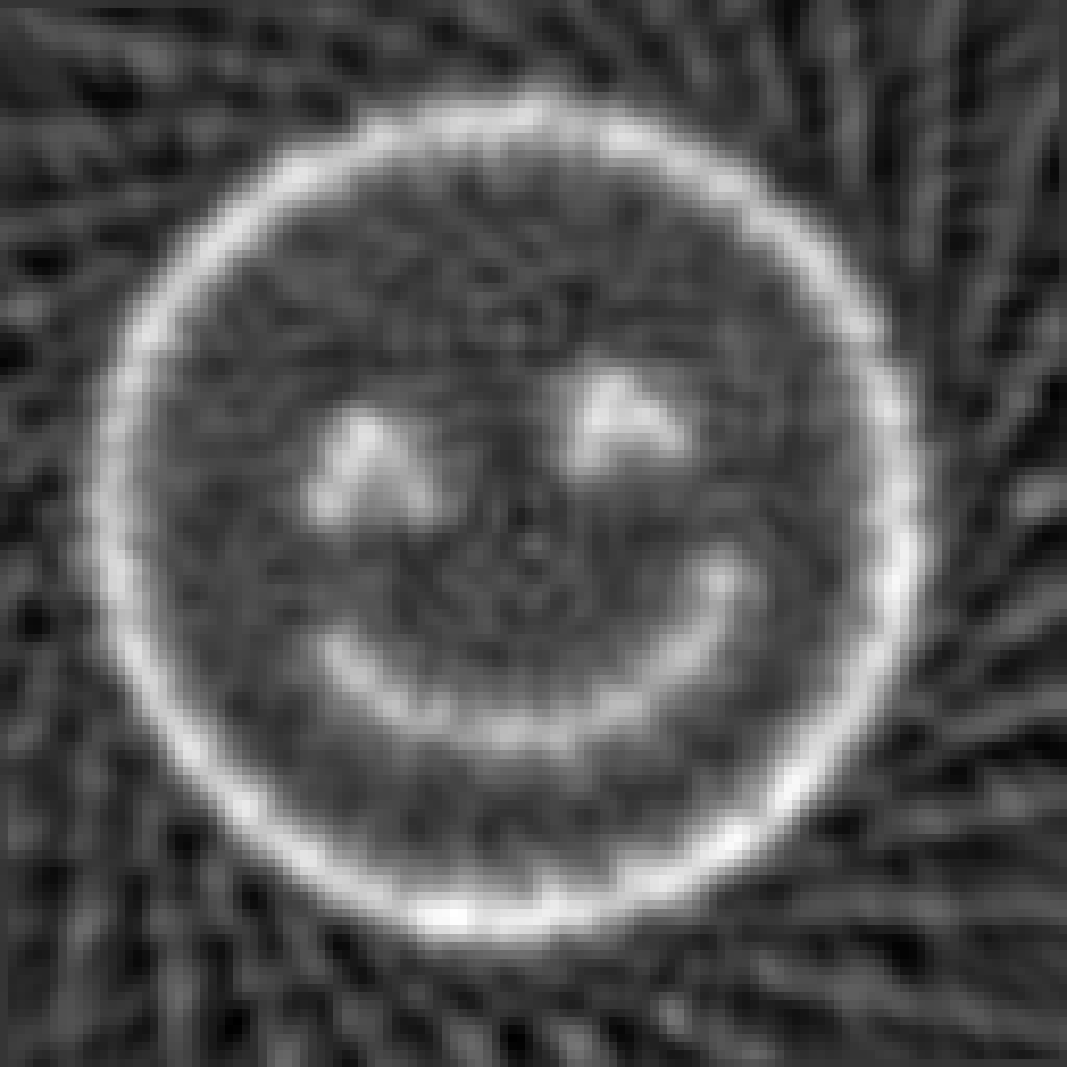}
	\end{minipage}\\
		\begin{minipage}{0.2\textwidth}
		\includegraphics[width=\textwidth]{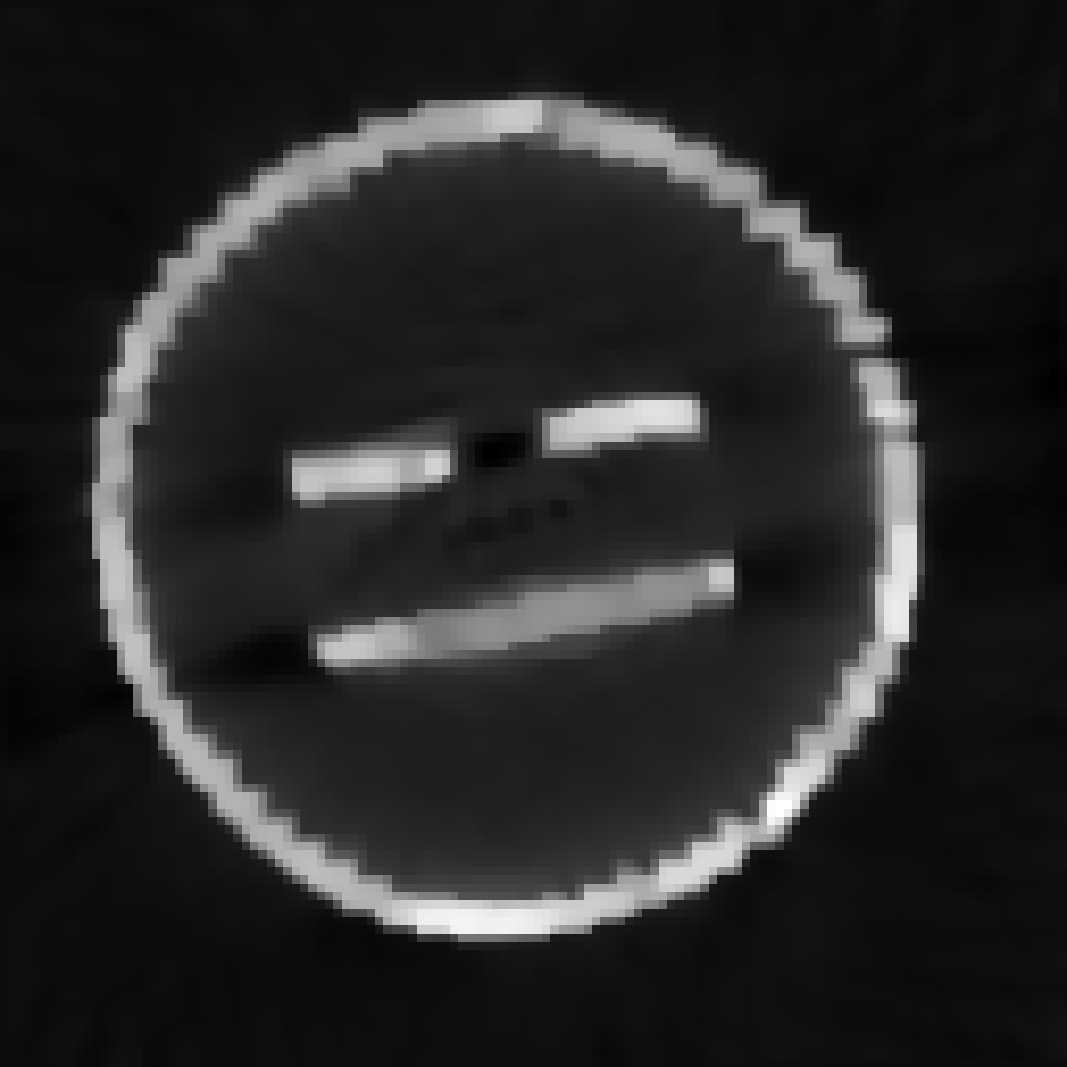}
	\end{minipage}
	\begin{minipage}{0.2\textwidth}
		\includegraphics[width=\textwidth]{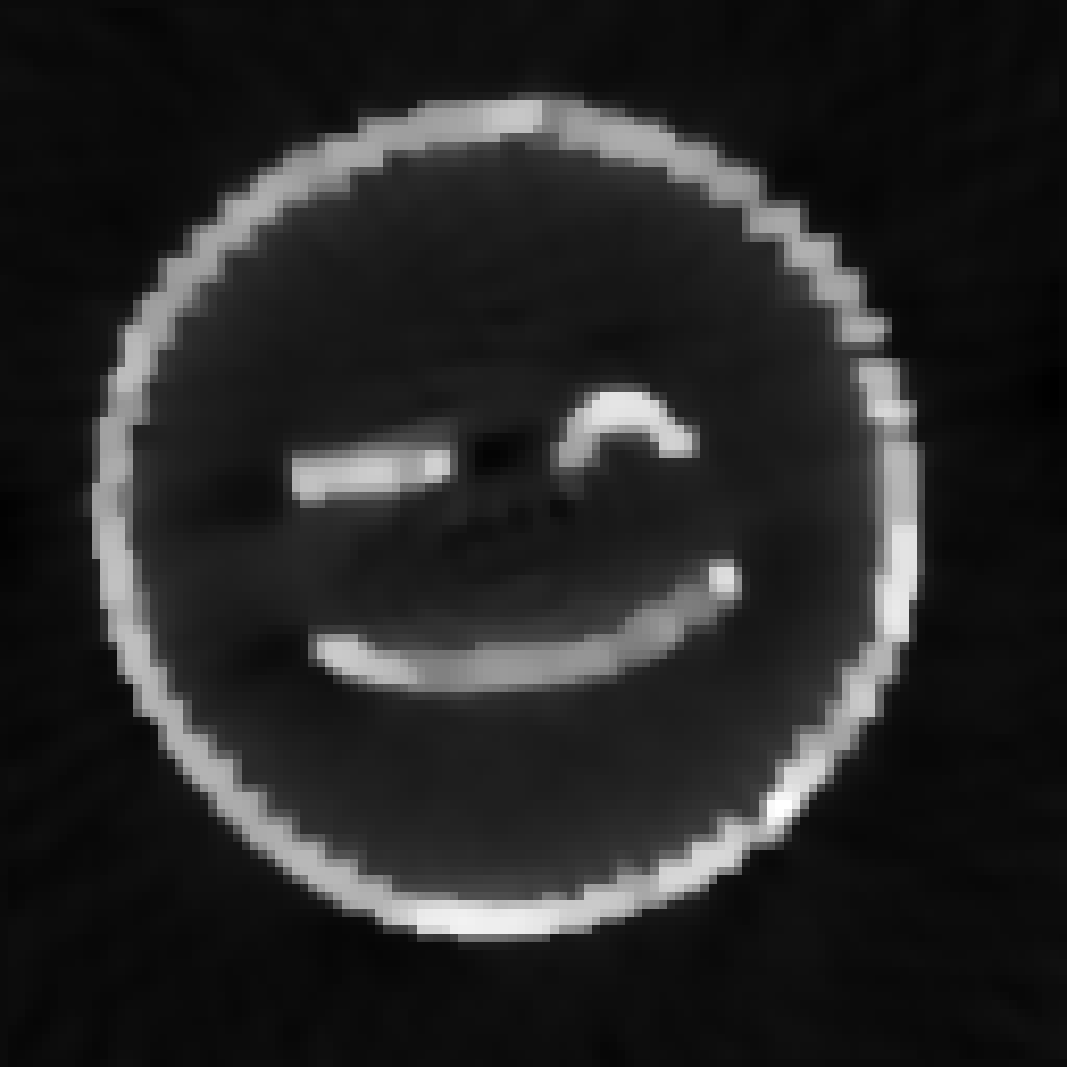}
	\end{minipage}
	\begin{minipage}{0.2\textwidth}
		\includegraphics[width=\textwidth]{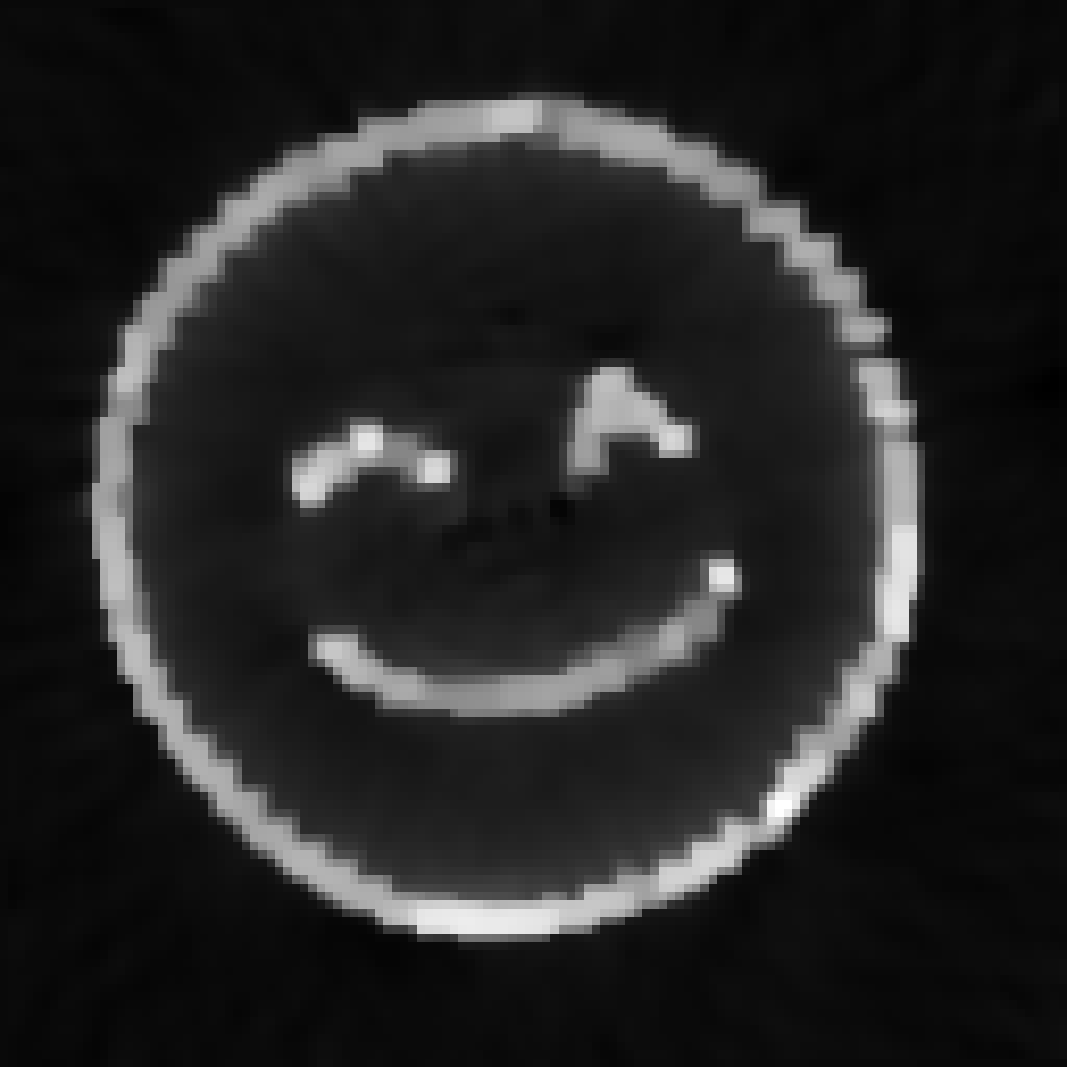}
	\end{minipage}
	\begin{minipage}{0.2\textwidth}
		\includegraphics[width=\textwidth]{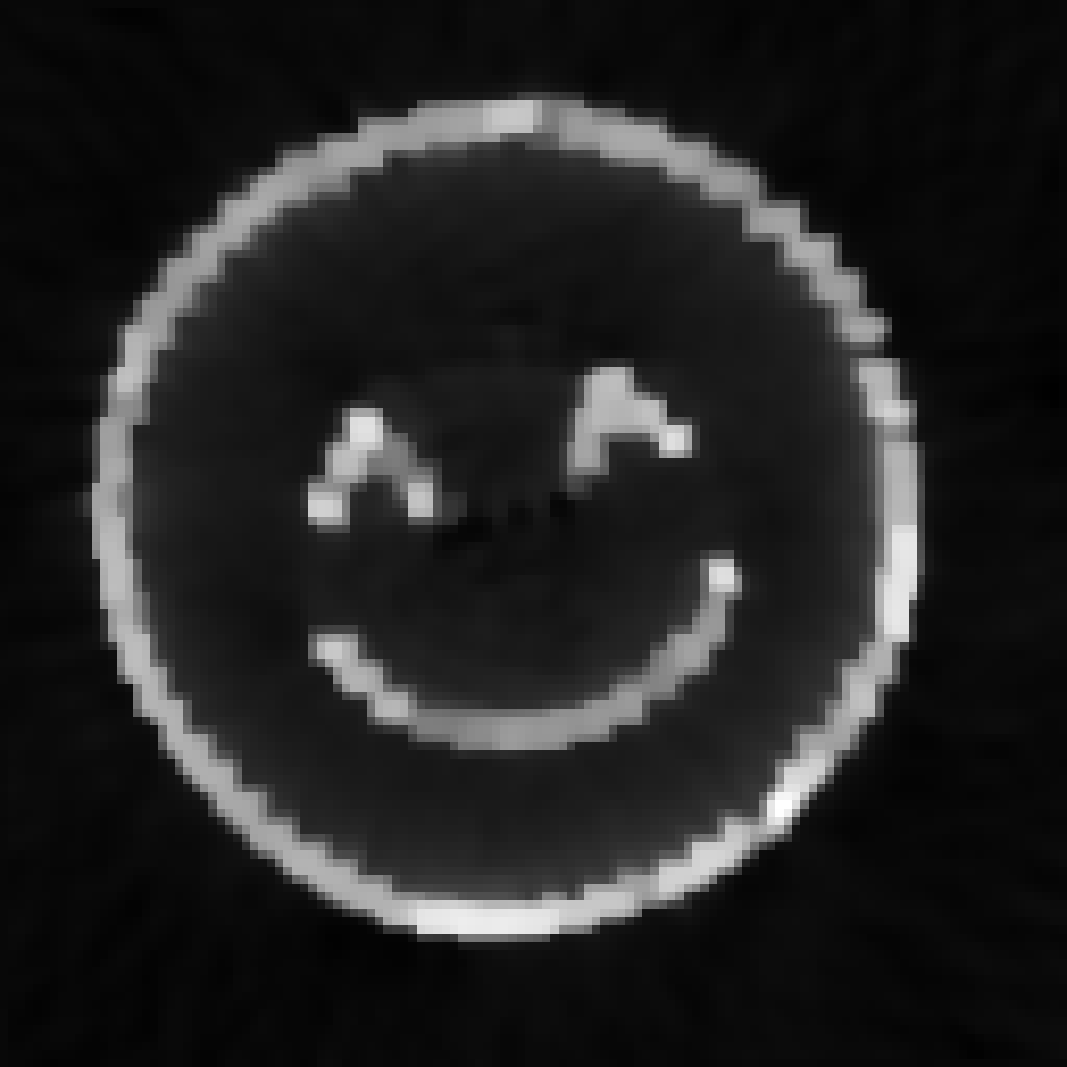}
	\end{minipage}\\
	\begin{minipage}{0.2\textwidth}
		\includegraphics[width=\textwidth]{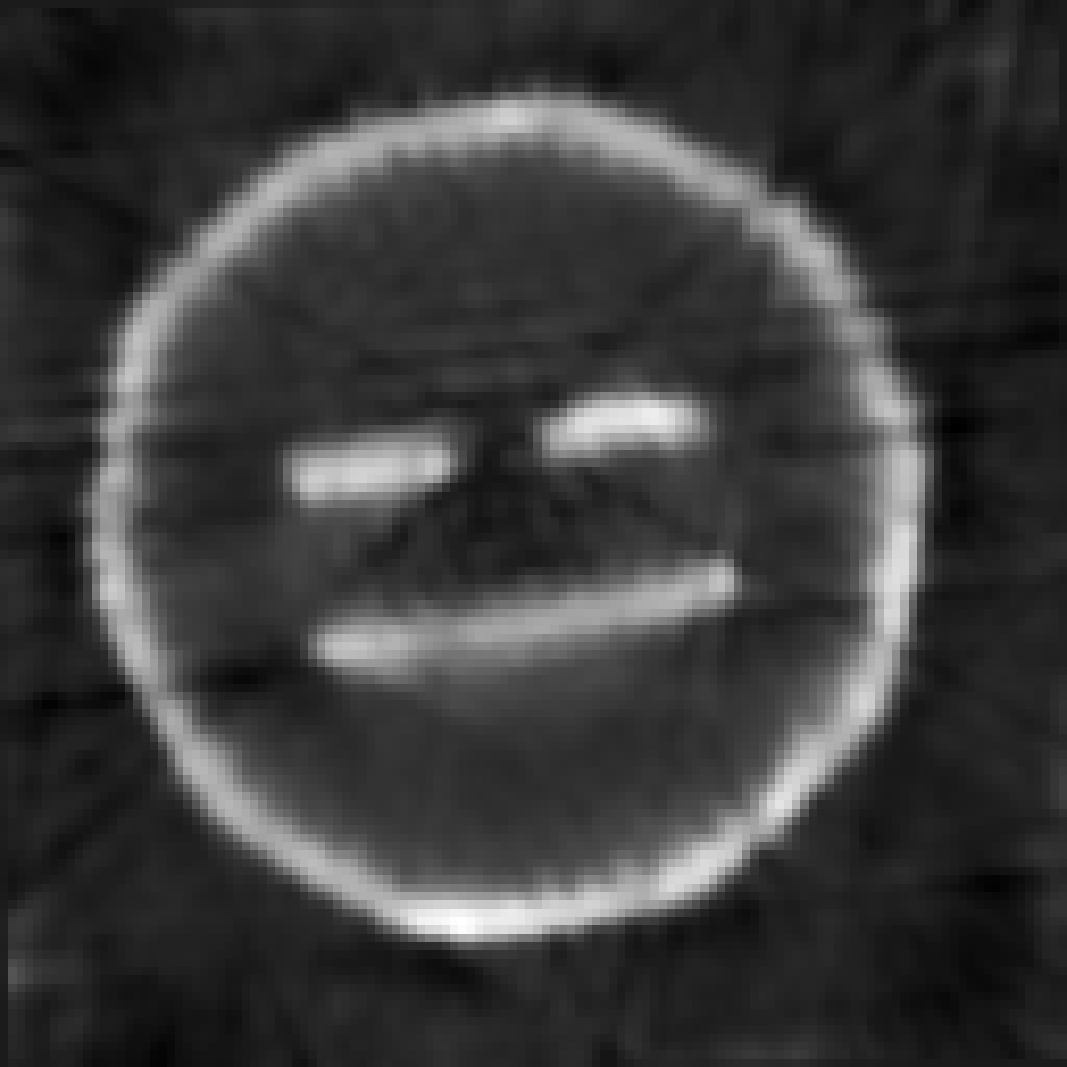}
	\end{minipage}
	\begin{minipage}{0.2\textwidth}
		\includegraphics[width=\textwidth]{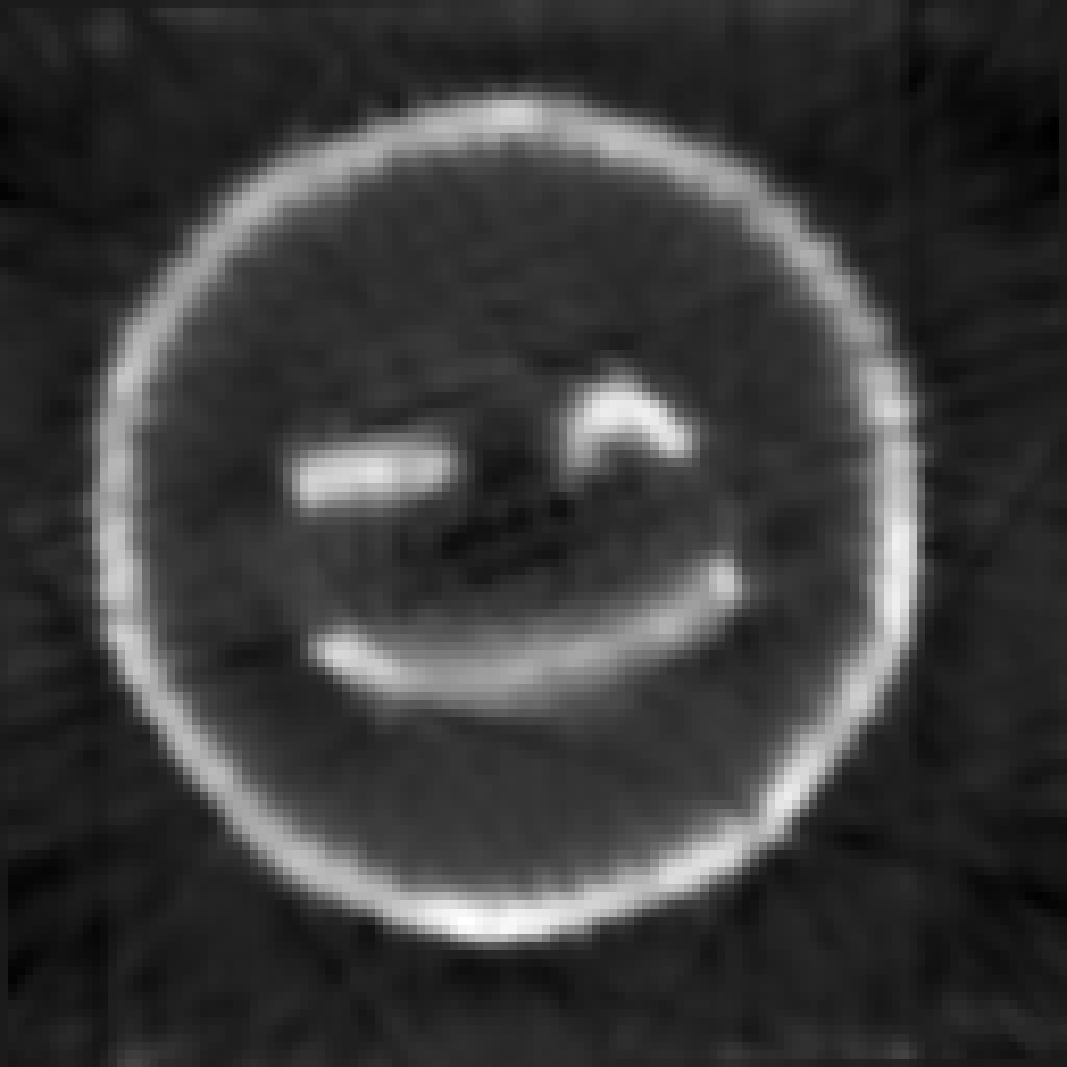}
	\end{minipage}
	\begin{minipage}{0.2\textwidth}
		\includegraphics[width=\textwidth]{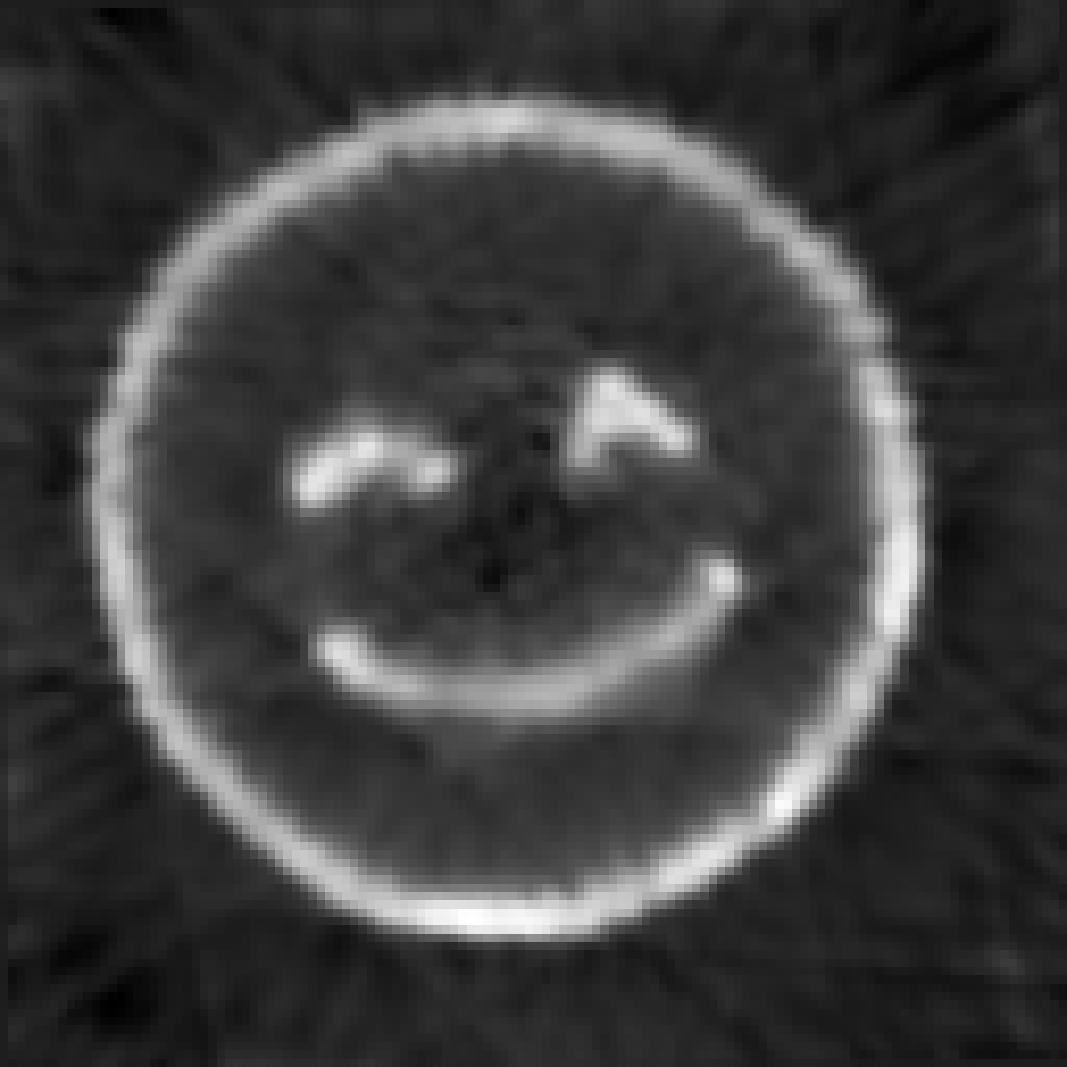}
	\end{minipage}
	\begin{minipage}{0.2\textwidth}
		\includegraphics[width=\textwidth]{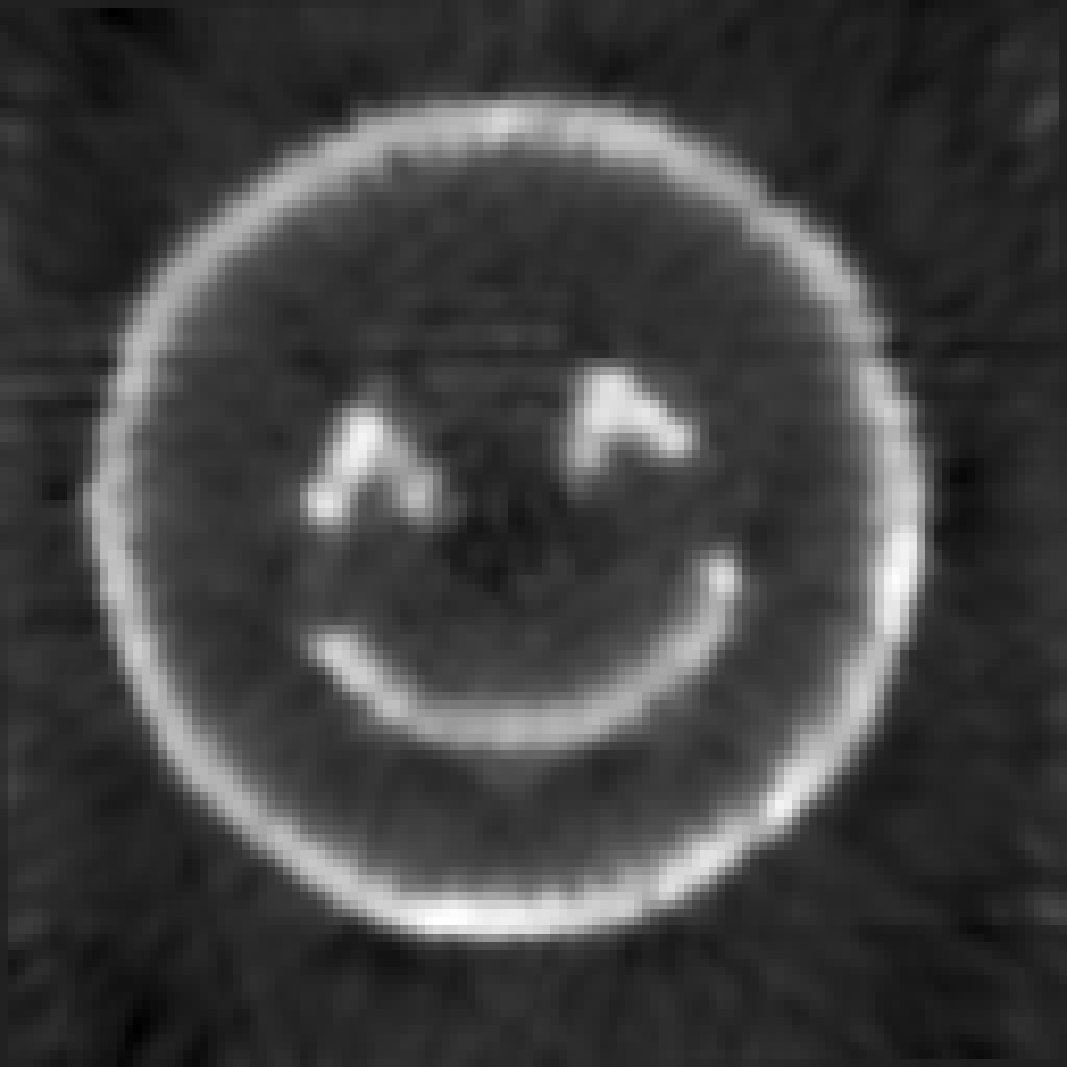}
	\end{minipage}\\
	\caption{Reconstruction results for the emoji test problem with $n_a =30$. The rows represent (from top to bottom): the original images, the reconstructions when images are considered independently, the reconstructions by AnisoTV,  the reconstructions by Iso3DTV, at time steps $t = 2, 10, 18, 31$ (from left to right).}
	\label{Fig: smile_30angles}
\end{figure}

It is evident (second row) that insufficiency of the information caused by the limited number of projection angles results in poor reconstructions, where the important details (features of the face) are missing as observed in the second row of Figure \ref{Fig: smile_10angles}. Solving the dynamic inverse problem (third and fourth rows) enhances the quality of the reconstruction. In particular, by considering the new regularization terms, we are able to reconstruct the edges clearly. This is observed in the third and the fourth rows of \ref{Fig: smile_10angles}. Moreover, the artifacts that arise from the limited angles are less present in the third and the fourth rows. 
\paragraph{Case 2: Consider $n_a =30$ projection angles}
In this second case we consider the full number of angles in the dataset \verb|DataDynamic_128x30.mat|, i.e., $n_a = 30$ and we validate the performance of the methods by highlighting the importance of the number of the projection angles.
Here $\bA^{(t)} \in \R^{6,510\times 16,384}$ and the measured sinograms are obtained from 217 projections in 30 angles, that is, $\bd^{(t)} \in \R^{6510}$. Hence $\bF \in \R^{214,830 \times 540,672}$ and $\bd \in \R^{214,830}$. 
The reconstructions of the static problems \eqref{eq: static} are shown in the second row of Figure \ref{Fig: smile_30angles}. The third and the fourth rows represent the reconstruction by AnisoTV and Iso3DTV at time instances $t = 6,14,20, 26$ from left to right, respectively. The first remarks is that, similar to the case when we consider $n_a = 10$, by solving the dynamic inverse problem, the obtained reconstruction has enhanced quality with respect to the solution of the static inverse problem. In addition, we observe that increasing the number of projection angles from 10 to 30 helps in removing the background artifacts and better preserving the edges (jumps).  

Although we only report the performance of AnisoTV and Iso3DTV (for diversification of the methods in different numerical examples), we remark that other methods such as TVpluTikhonov, IsoTV, and Aniso3DTV produce reconstructions of similar quality to AnisoTV and Iso3DTV. In contrast to the other test problems that we presented above, where GS was one of the most accurate methods, in this example it is the less accurate method. This observation allows us to highlight one of the goals of this paper, that is to present a variety of regularization methods without advocating for one over the other as the methods that we describe are application dependent.

\begin{figure}[h!]
\centering
	\begin{minipage}{0.2\textwidth}
		\includegraphics[width=\textwidth]{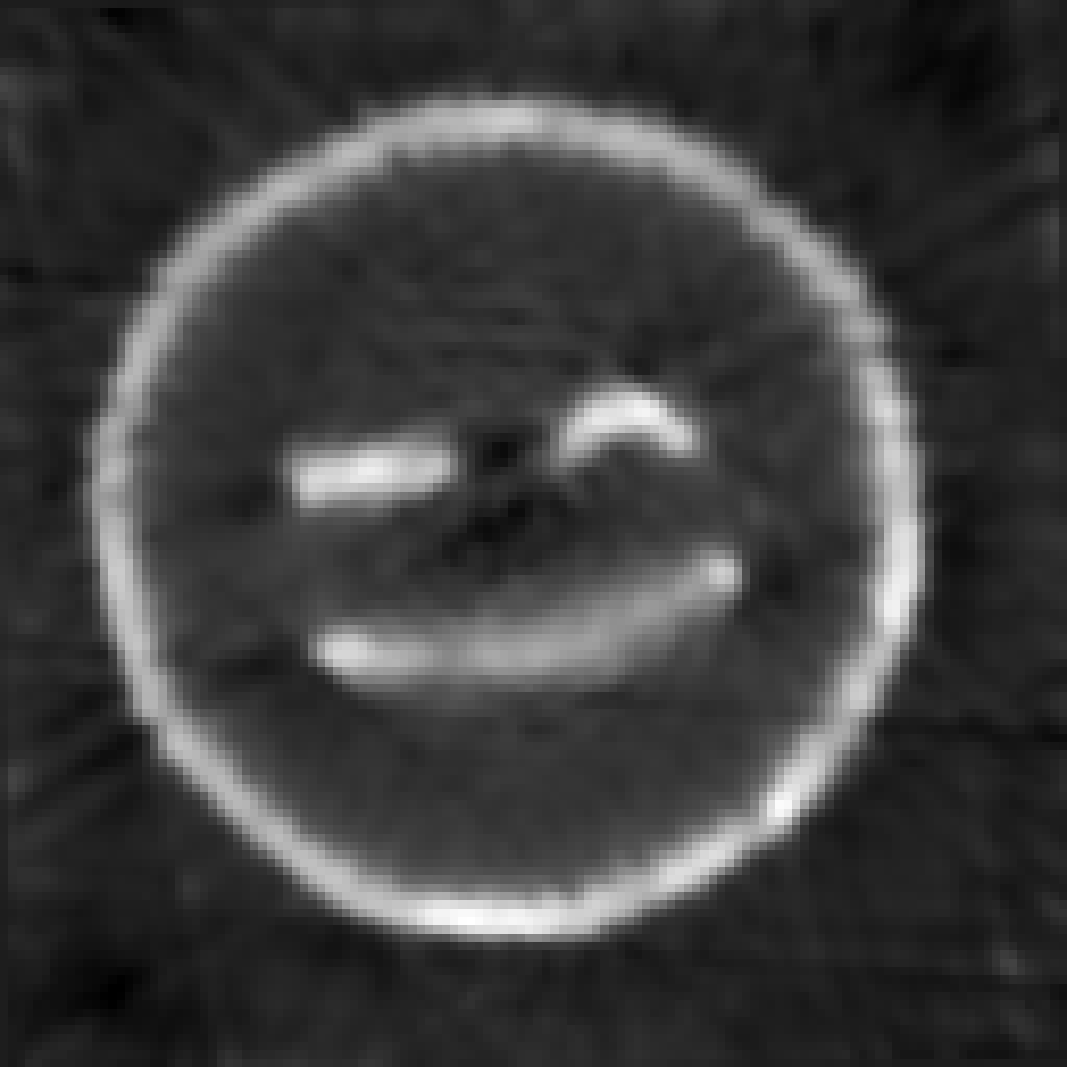}
	\end{minipage}
	\begin{minipage}{0.2\textwidth}
		\includegraphics[width=\textwidth]{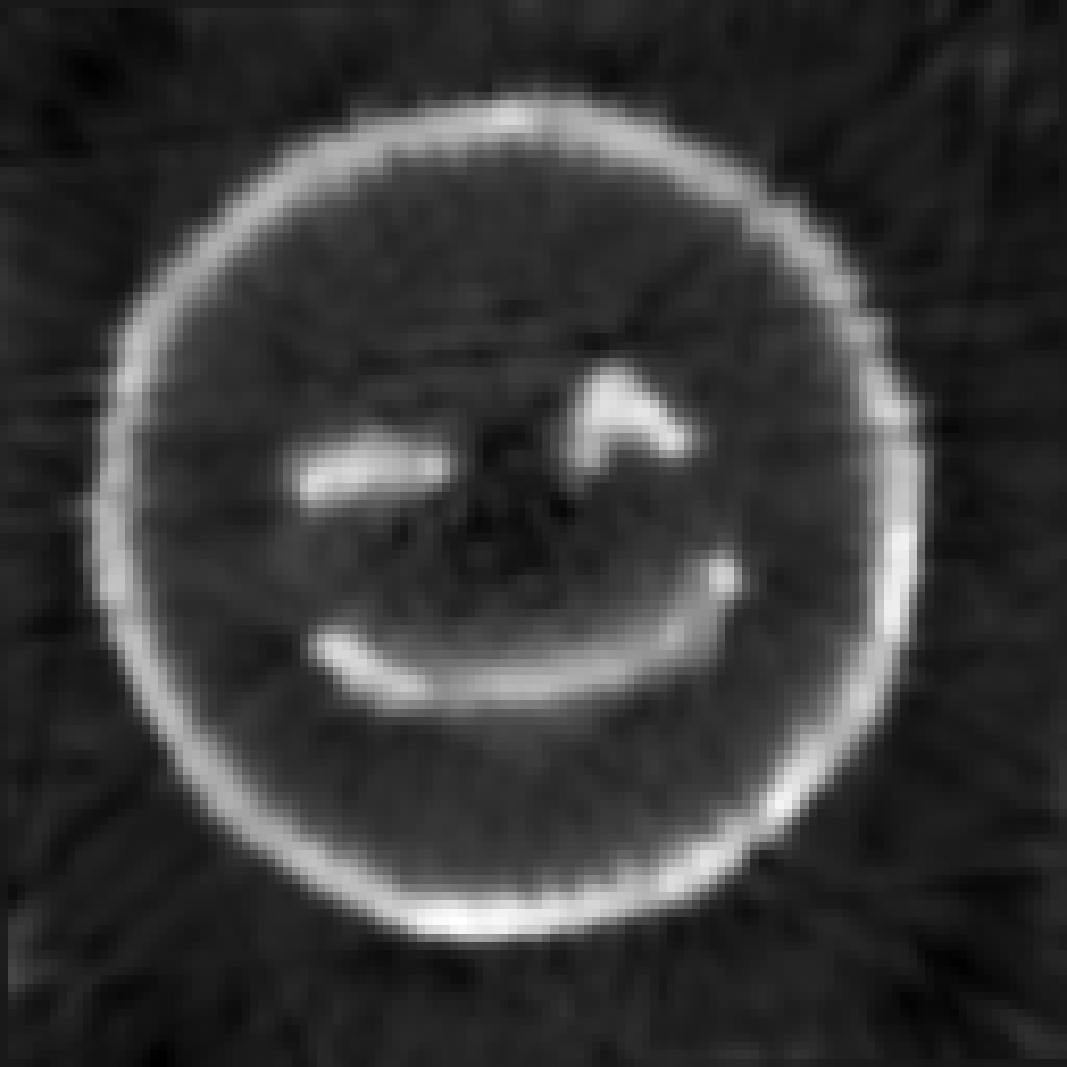}
	\end{minipage}
	\begin{minipage}{0.2\textwidth}
		\includegraphics[width=\textwidth]{Figures/Emoji_3DTV_var_rec_10_20.png}
	\end{minipage}
	\begin{minipage}{0.2\textwidth}
		\includegraphics[width=\textwidth]{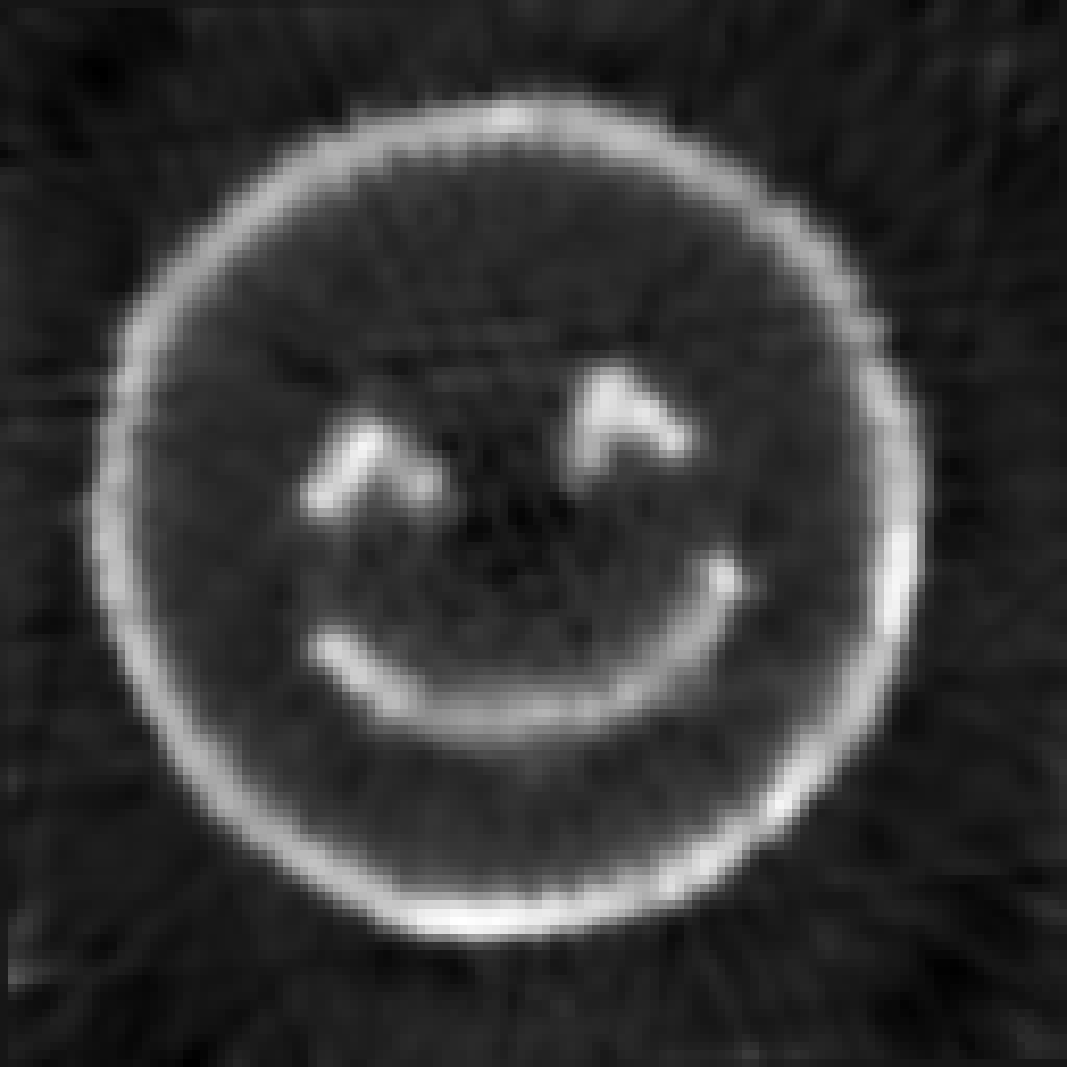}
	\end{minipage}

	\begin{minipage}{0.2\textwidth}
		\includegraphics[width=\textwidth]{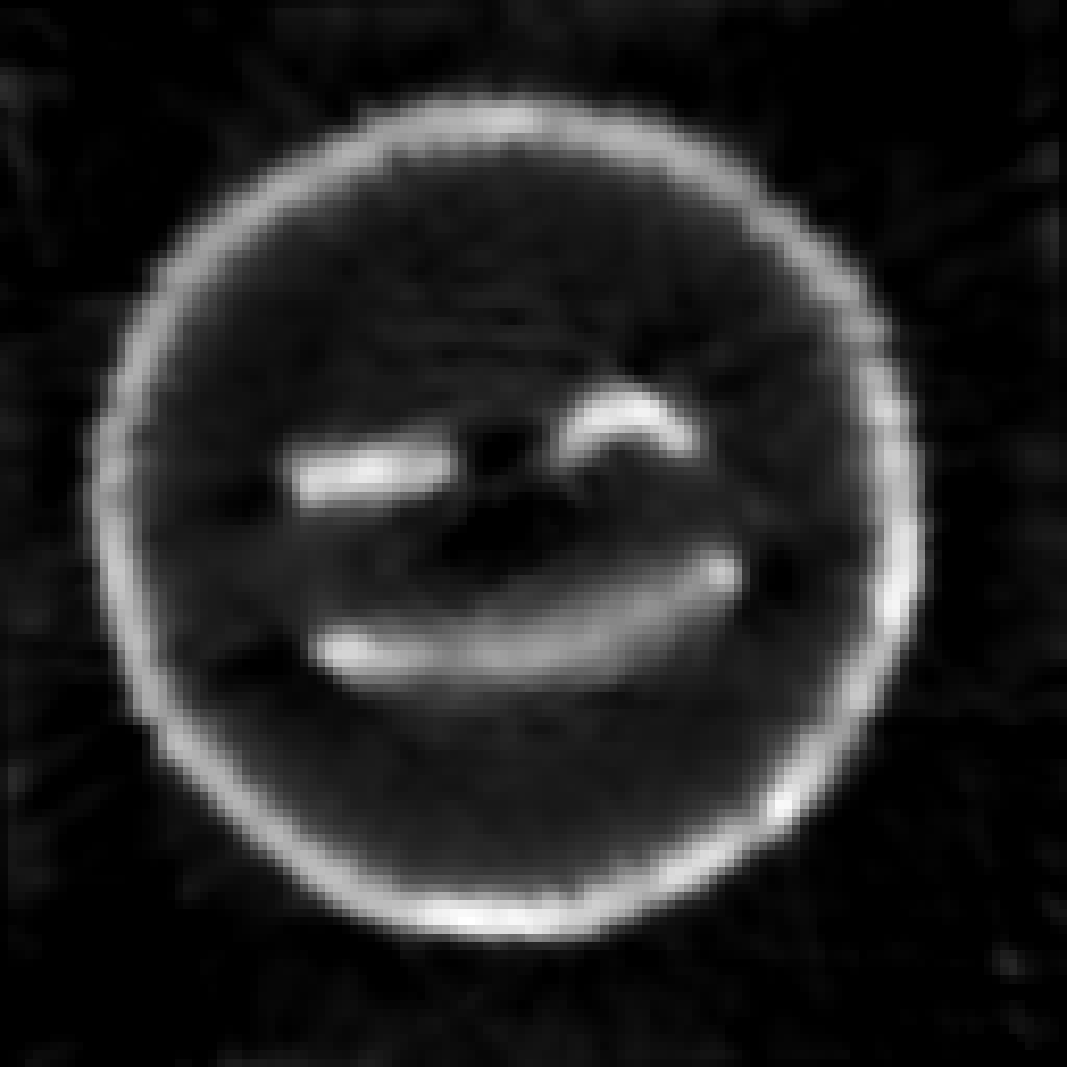}
	\end{minipage}
	\begin{minipage}{0.2\textwidth}
		\includegraphics[width=\textwidth]{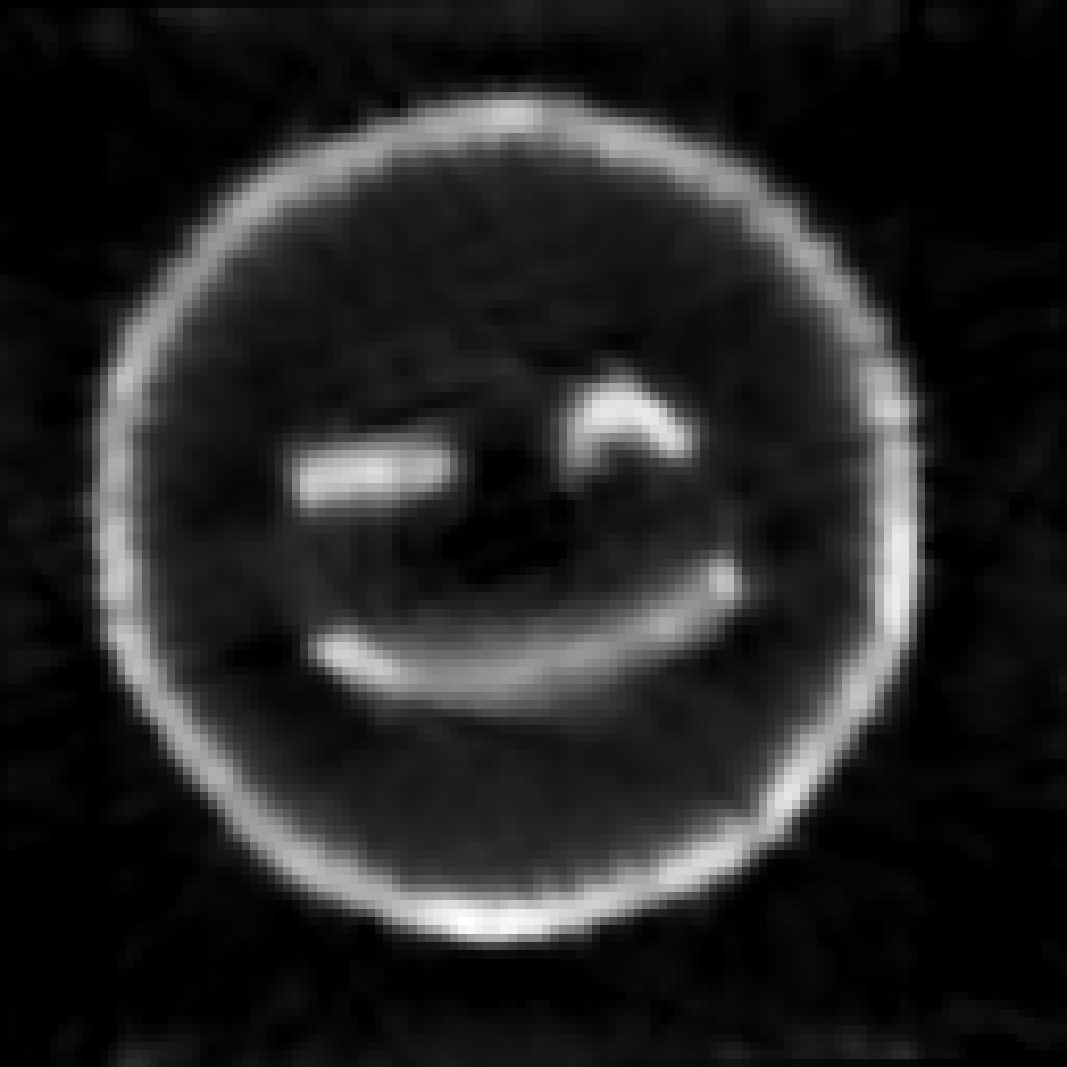}
	\end{minipage}
	\begin{minipage}{0.2\textwidth}
		\includegraphics[width=\textwidth]{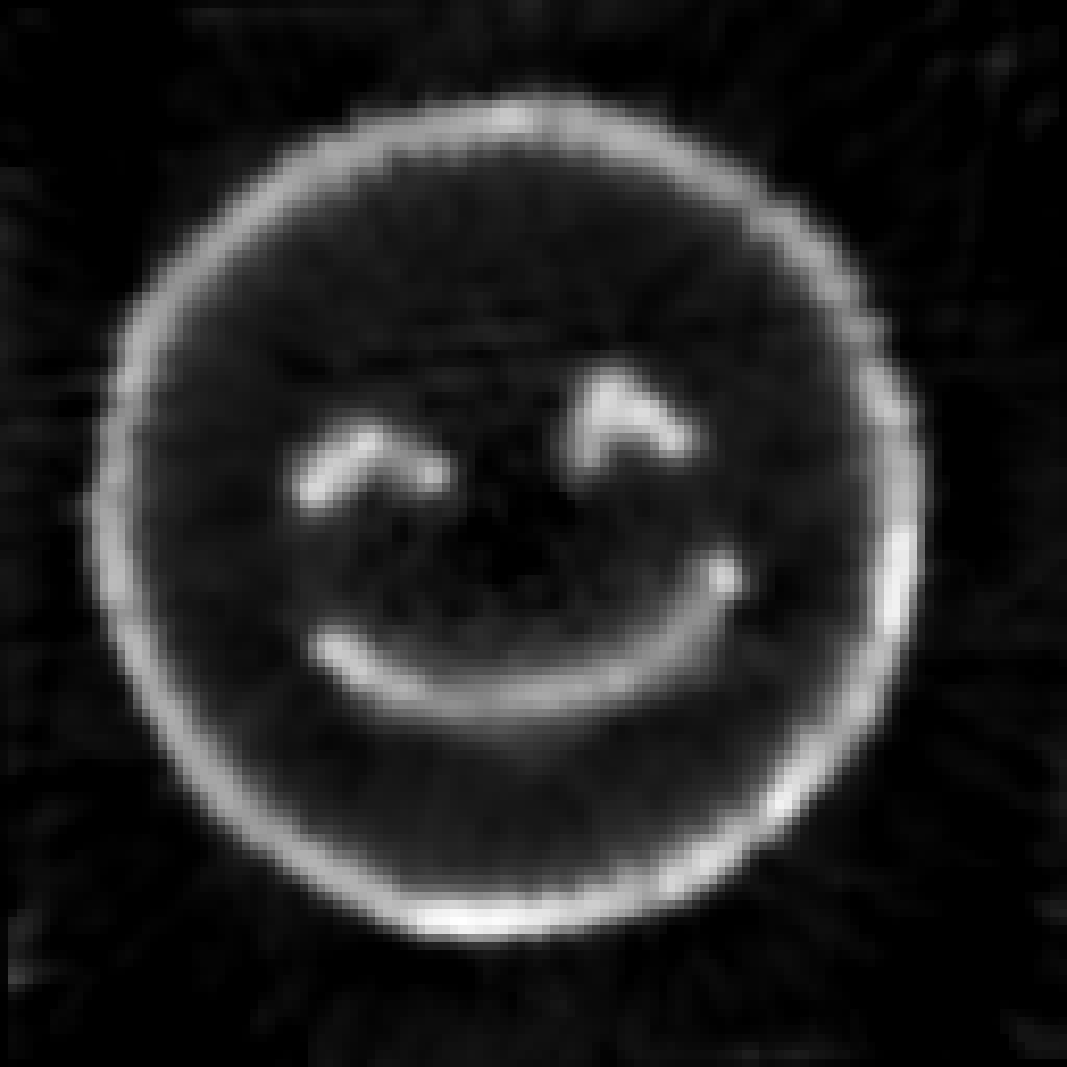}
	\end{minipage}
	\begin{minipage}{0.2\textwidth}
		\includegraphics[width=\textwidth]{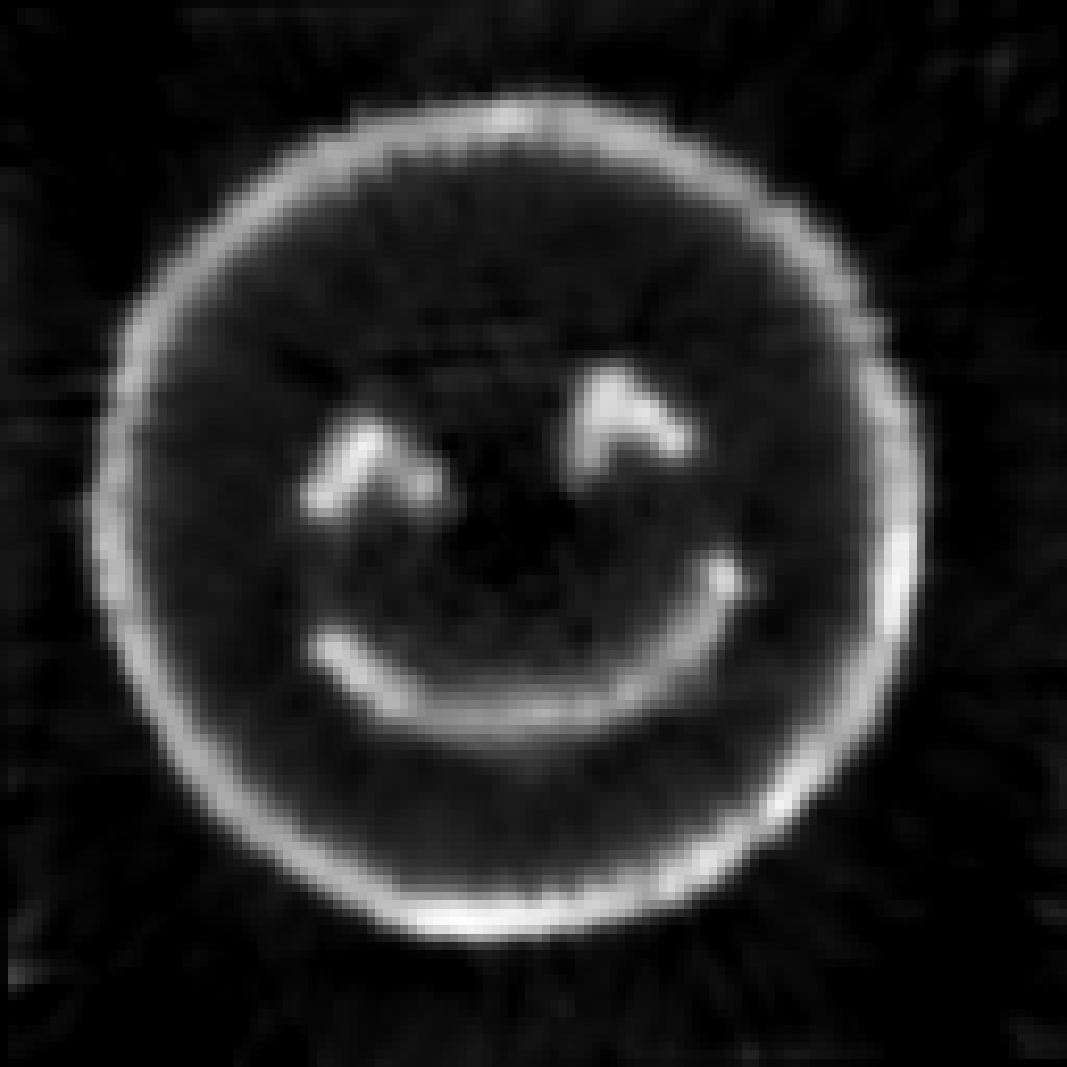}
	\end{minipage}

	\caption{Reconstruction results for the emoji test problem with $n_a = 10$. The firs row shows the reconstructions with Iso3DTV for the unconstrained problem and the second row shows the reconstructions with nonnegativity Iso3DTV at time steps %$n_t = 3,7,11,13$ 
	$n_t = 6,14,20, 26$
	respectively from left to right.}
	\label{Fig: 10angles_nonnegativity}
\end{figure}

\begin{table}[h!]
	\begin{center}
	\caption{ Dynamic X-Ray Tomography example: The number of iterations when the discrepancy principle is satisfied for the first time and the regularization parameters at those iterations for AnisoTV, and Iso3DTV.}
	\label{tab:ex3}
		\begin{tabular}{c|ccccc} 
			\hline 
			\multicolumn{1}{c}{$n_a$} & \multicolumn{1}{c}{-}&\multicolumn{1}{c}{AnisoTV}
			&\multicolumn{1}{c}{Iso3DTV} & 
			\\ \hline 
			\multirow{2}{4em}{$10$} &iter & $150^{*}$ & $150^{*}$\\
			&$\lambda$ &$0.26$& $0.56$&   \\ \hline
			\multirow{2}{4em}{$30$}&iter&$150^{*}$& $150^{*}$  \\ 
			&$\lambda$ &0.31& $0.61$  \\ \hline
		\end{tabular} 
	\end{center} %\vspace{.5cm}
\end{table}
\paragraph{Nonnegativity constraint}
In many application, such as medical imaging and  astronomical imaging, the pixels of the desired solution are nonnegative \cite{buccini2020linearized, buccini2020modulus, gazzola2017fast}, that is, the exact solution of 
\eqref{eq: l2-lq} is known to live in the closed and convex set
\[
\Omega_{0}=\{\bu\in\R^n:\bu_{\ell}\geq 0,~~\ell=1,2,\ldots,n\},
\]
where $\bu_{\ell}$ denotes the $\ell$-th entry of the vector $\bu$. In general, imposing nonnegativity helps mitigating the artifacts that arise from limited angles.
Here we consider the optimization problems \eqref{eqn:dynamic} subject to the constraint $\bu\in\Omega_0$. This is heuristically implemented by projecting the solution $\bu_{(k)}$ onto $\Omega_0$ at each iteration. 
We illustrate the effect of the nonnegativity constraint in Example 3 for \emph{Case 1}, where the number of projection angles $n_a = 10$ and 1\% Gaussian white noise was artificially added to the observations $\bd^{(t)}$, $t = 1,2,\dots,33$. 
The reconstructed images at time steps 
%3, 7, 11, and 13 
$t = 6,14,20, 26$
are shown in Figure \ref{Fig: 10angles_nonnegativity}. 
From visual inspection, the artifacts around the edges are less present when the nonnegativity constraint is applied. A more detailed analysis on the nonnegativity constraints is left as a future work.

\section{Conclusions and future direction}\label{sec: conclusions}
In this paper we proposed new methods for solving large scale dynamic inverse problems and providing solutions with edge preserving and sparsity promoting properties. The approaches that we propose here are grouped into isotropic TV methods (which include IsoTV and Iso3DTV), anisotropic TV (which include AnisoTV and an Aniso3DTV), and another group of methods based on the concept of group sparsity, GS. All the methods can be expressed in a unified framework using the majorization-minimization technique where the resulting least squares problem can be solved on a generalized Krylov subspace of relatively small dimension and the regularization parameter can be estimated efficiently. Several numerical examples, performed on both synthetic and real data, illustrate the performances of the proposed methods in terms of the quality of the reconstructed solutions.
Although we propose a unified and generic framework that can be used to solve a wide range of dynamic inverse problems, there are quite a few potential directions to investigate for future work. Some of them are listed within Section~\ref{sec: alternative}. Here we emphasise again that one direction of interest is to investigate more on the selection of domain dedicated regularization parameters, for instance, regularization parameters for temporal and spatial domain or adapted regularization parameters for different channels. Another direction includes alternative formulations along with their Bayesian interpretation and uncertainty quantification. Moreover, it is known that tensor formulations preserve the structure of the data, hence we are interested in investigating efficient tensor based regularization methods, \cite{semerci2014tensor}. Some applications of interest include video reconstruction, multi-channel x-ray spectral tomography, and moving object detection.

\section*{Acknowledgments}
This work was initiated as a part of the Statistical and
Applied Mathematical Sciences Institute (SAMSI) Program on Numerical Analysis in Data Science in 2020. Any opinions, findings, and conclusions or recommendations expressed
in this material are those of the authors and do not necessarily reflect the views of the National Science Foundation (NSF). 
The authors thank Professors Misha Kilmer and D.\ Andrew Brown for the constructive discussions and suggestions. MP gratefully acknowledges ASU Research Computing facilities for the computing resources used for testing purposes. AKS would like to acknowledge partial support from NSF through the awards DMS-1845406 and DMS-1720398. The work of SG is partially supported by EPSRC, under grant EP/T001593/1. The work by EdS was supported in part by NSF grant DMS 1720305.
 \bibliographystyle{abbrv}
 \bibliography{biblio}

\begin{thebibliography}{10}

\bibitem{achenbach2001detection}
S.~Achenbach, T.~Giesler, D.~Ropers, S.~Ulzheimer, H.~Derlien, C.~Schulte,
  E.~Wenkel, W.~Moshage, W.~Bautz, W.~G. Daniel, W.~A. Kalender, and U.~Baum.
\newblock Detection of coronary artery stenoses by contrast-enhanced,
  retrospectively electrocardiographically-gated, multislice spiral computed
  tomography.
\newblock {\em Circulation}, 103(21):2535--2538, 2001.

\bibitem{acton2009diffusion}
S.~T. Acton.
\newblock Diffusion partial differential equations for edge detection.
\newblock In {\em The Essential Guide to Image Processing}, pages 525--552.
  Elsevier, 2009.

\bibitem{ashwini2017sparse}
K.~Ashwini, R.~Amutha, and K.~Harini.
\newblock Sparse based image compression in wavelet domain.
\newblock In {\em 2017 International Conference on Signal Processing and
  Communication (ICSPC)}, pages 58--62. IEEE, 2017.

\bibitem{bach2008exploring}
F.~Bach.
\newblock Exploring large feature spaces with hierarchical multiple kernel
  learning.
\newblock {\em arXiv preprint arXiv:0809.1493}, 2008.

\bibitem{bach2012optimization}
F.~Bach, R.~Jenatton, J.~Mairal, and G.~Obozinski.
\newblock Optimization with sparsity-inducing penalties.
\newblock {\em Found. Trends Mach. Learn.}, 4(1):1--106, 2012.

\bibitem{blondel20043d}
C.~Blondel, R.~Vaillant, G.~Malandain, and N.~Ayache.
\newblock 3{D} tomographic reconstruction of coronary arteries using a
  precomputed 4{D} motion field.
\newblock {\em Physics in Medicine \& Biology}, 49(11):2197, 2004.

\bibitem{buccini2020linearized}
A.~Buccini, M.~Pasha, and L.~Reichel.
\newblock Linearized {K}rylov subspace {B}regman iteration with nonnegativity
  constraint.
\newblock {\em Numer. Algorithms}, pages 1--24, 2020.

\bibitem{buccini2020modulus}
A.~Buccini, M.~Pasha, and L.~Reichel.
\newblock Modulus-based iterative methods for constrained $\ell_p-\ell_q$
  minimization.
\newblock {\em Inverse Problems}, 36(8):084001, 2020.

\bibitem{buccini2}
A.~Buccini and L.~Reichel.
\newblock An $\ell_2-\ell_q$ regularization method for large discrete ill-posed
  problems.
\newblock {\em Journal of Scientific Computing}, 78:1526--1549, 2019.

\bibitem{buccini2021generalized}
A.~Buccini and L.~Reichel.
\newblock Generalized cross validation for $\ell_p-\ell_q$ minimization.
\newblock {\em Numerical Algorithms,
  https://doi.org/10.1007/s11075-021-01087-9}, 2021.

\bibitem{burger2017variational}
M.~Burger, H.~Dirks, L.~Frerking, A.~Hauptmann, T.~Helin, and S.~Siltanen.
\newblock A variational reconstruction method for undersampled dynamic x-ray
  tomography based on physical motion models.
\newblock {\em Inverse Problems}, 33(12):124008, 2017.

\bibitem{COS09}
J.-F. Cai, S.~Osher, and Z.~Shen.
\newblock Split {B}regman methods and frame based image restoration.
\newblock {\em Multiscale Modeling \& Simulation}, 8:337--369, 2009.

\bibitem{candes2006robust}
E.~J. Cand{\`e}s, J.~Romberg, and T.~Tao.
\newblock Robust uncertainty principles: {E}xact signal reconstruction from
  highly incomplete frequency information.
\newblock {\em IEEE Trans. Inform. Theory}, 52(2):489--509, 2006.

\bibitem{caselles2015total}
V.~Caselles, A.~Chambolle, and M.~Novaga.
\newblock Total variation in imaging.
\newblock {\em Handbook of mathematical methods in imaging}, 1:1455--1499,
  2015.

\bibitem{chambolle2005total}
A.~Chambolle.
\newblock Total variation minimization and a class of binary {MRF} models.
\newblock In {\em International Workshop on Energy Minimization Methods in
  Computer Vision and Pattern Recognition}, pages 136--152. Springer, 2005.

\bibitem{chartrand2007exact}
R.~Chartrand.
\newblock Exact reconstruction of sparse signals via nonconvex minimization.
\newblock {\em IEEE Signal Process. Lett.}, 14(10):707--710, 2007.

\bibitem{choksi2010anisotropic}
R.~Choksi, Y.~van Gennip, and A.~Oberman.
\newblock Anisotropic total variation regularized $\ell_1$-approximation and
  denoising/deblurring of 2{D} bar codes.
\newblock {\em arXiv preprint arXiv:1007.1035}, 2010.

\bibitem{chung2017motion}
J.~Chung and L.~Nguyen.
\newblock Motion estimation and correction in photoacoustic tomographic
  reconstruction.
\newblock {\em SIAM J. Imaging Sci.}, 10(1):216--242, 2017.

\bibitem{chung2018efficient}
J.~Chung, A.~K. Saibaba, M.~Brown, and E.~Westman.
\newblock Efficient generalized {G}olub--{K}ahan based methods for dynamic
  inverse problems.
\newblock {\em Inverse Problems}, 34(2):024005, 2018.

\bibitem{condat2017discrete}
L.~Condat.
\newblock Discrete total variation: {N}ew definition and minimization.
\newblock {\em SIAM Journal on Imaging Sciences}, 10(3):1258--1290, 2017.

\bibitem{cui2017classification}
B.~Cui, X.~Ma, X.~Xie, G.~Ren, and Y.~Ma.
\newblock Classification of visible and infrared hyperspectral images based on
  image segmentation and edge-preserving filtering.
\newblock {\em Infrared Physics \& Technology}, 81:79--88, 2017.

\bibitem{daubechies1992ten}
I.~Daubechies.
\newblock {\em Ten lectures on wavelets}.
\newblock SIAM, 1992.

\bibitem{dgm1986painless}
I.~Daubechies, A.~Grossmann, and Y.~Meyer.
\newblock Painless nonorthogonal expansions.
\newblock {\em J. Math. Phys}, 27:1271--1283, 1986.

\bibitem{engl1996regularization}
H.~W. Engl, M.~Hanke, and A.~Neubauer.
\newblock {\em Regularization of inverse problems}, volume 375.
\newblock Springer Science \& Business Media, 1996.

\bibitem{esedoglu2004decomposition}
S.~Esedo\={g}lu and S.~J. Osher.
\newblock Decomposition of images by the anisotropic {R}udin-{O}sher-{F}atemi
  model.
\newblock {\em Comm. Pure Appl. Math.}, 57(12):1609--1626, 2004.

\bibitem{gao2019iterated}
R.~Gao, F.~Tronarp, and S.~S{\"a}rkk{\"a}.
\newblock Iterated extended {K}alman smoother-based variable splitting for
  l$_1$-regularized state estimation.
\newblock {\em IEEE Trans. Signal Process.}, 67(19):5078--5092, 2019.

\bibitem{gazzola2019ir}
S.~Gazzola, P.~C. Hansen, and J.~G. Nagy.
\newblock {IR} tools: a {MATLAB} package of iterative regularization methods
  and large-scale test problems.
\newblock {\em Numer. Algorithms}, 81(3):773--811, 2019.

\bibitem{gazzola2017fast}
S.~Gazzola and Y.~Wiaux.
\newblock Fast nonnegative least squares through flexible {K}rylov subspaces.
\newblock {\em SIAM J. Sci. Comput.}, 39(2):A655--A679, 2017.

\bibitem{giraldo2013weighted}
E.~Giraldo, J.~Casta{\~n}o-Candamil, and G.~Castellanos-Dominguez.
\newblock A weighted dynamic inverse problem for electroencephalographic
  current density reconstruction.
\newblock In {\em 2013 6th International IEEE/EMBS Conference on Neural
  Engineering (NER)}, pages 521--524. IEEE, 2013.

\bibitem{gonzalez2017isotropic}
G.~Gonz{\'a}lez, V.~Kolehmainen, and A.~Sepp{\"a}nen.
\newblock Isotropic and anisotropic total variation regularization in
  electrical impedance tomography.
\newblock {\em Comput. Math. with Appl.}, 74(3):564--576, 2017.

\bibitem{guo2018edge}
F.~Guo, C.~Zhang, and M.~Zhang.
\newblock Edge-preserving image denoising.
\newblock {\em IET Image Process.}, 12(8):1394--1401, 2018.

\bibitem{hanke1993regularization}
M.~Hanke and P.~C. Hansen.
\newblock Regularization methods for large-scale problems.
\newblock {\em Surv. Math. Ind}, 3(4):253--315, 1993.

\bibitem{hansen1998rank}
P.~C. Hansen.
\newblock {\em Rank-deficient and discrete ill-posed problems: numerical
  aspects of linear inversion}.
\newblock SIAM, 1998.

\bibitem{he2016image}
K.~He, D.~Wang, and X.~Zheng.
\newblock Image segmentation on adaptive edge-preserving smoothing.
\newblock {\em J. Electron. Imaging}, 25(5):053022, 2016.

\bibitem{herzog2012directional}
R.~Herzog, G.~Stadler, and G.~Wachsmuth.
\newblock Directional sparsity in optimal control of partial differential
  equations.
\newblock {\em SIAM Journal on Control and Optimization}, 50(2):943--963, 2012.

\bibitem{4193460}
W.~S. {Hoge}, M.~E. {Kilmer}, C.~{Zacarias-Almarcha}, and D.~H. {Brooks}.
\newblock Fast regularized reconstruction of non-uniformly subsampled
  partial-{F}ourier parallel {MRI} data.
\newblock In {\em 2007 4th IEEE International Symposium on Biomedical Imaging:
  From Nano to Macro}, pages 1012--1015, 2007.

\bibitem{hu2016moving}
W.~Hu, Y.~Yang, W.~Zhang, and Y.~Xie.
\newblock Moving object detection using tensor-based low-rank and saliently
  fused-sparse decomposition.
\newblock {\em IEEE Trans. Image Process.}, 26(2):724--737, 2016.

\bibitem{huang2017majorization}
G.~Huang, A.~Lanza, S.~Morigi, L.~Reichel, and F.~Sgallari.
\newblock Majorization--minimization generalized {K}rylov subspace methods for
  $\ell_p-\ell_q$ optimization applied to image restoration.
\newblock {\em BIT}, 57(2):351--378, 2017.

\bibitem{hunter2004tutorial}
D.~R. Hunter and K.~Lange.
\newblock A tutorial on {MM} algorithms.
\newblock {\em The American Statistician}, 58(1):30--37, 2004.

\bibitem{kaipio2006statistical}
J.~Kaipio and E.~Somersalo.
\newblock {\em Statistical and computational inverse problems}, volume 160.
\newblock Springer Science \& Business Media, 2006.

\bibitem{kazantsev2018joint}
D.~Kazantsev, J.~S. J{\o}rgensen, M.~S. Andersen, W.~R. Lionheart, P.~D. Lee,
  and P.~J. Withers.
\newblock Joint image reconstruction method with correlative multi-channel
  prior for x-ray spectral computed tomography.
\newblock {\em Inverse Problems}, 34(6):064001, 2018.

\bibitem{kim2010tree}
S.~Kim and E.~P. Xing.
\newblock Tree-guided group lasso for multi-task regression with structured
  sparsity.
\newblock In {\em ICML}, 2010.

\bibitem{kolda2009tensor}
T.~G. Kolda and B.~W. Bader.
\newblock Tensor decompositions and applications.
\newblock {\em SIAM Rev.}, 51(3):455--500, 2009.

\bibitem{krishnan2009fast}
D.~Krishnan and R.~Fergus.
\newblock Fast image deconvolution using hyper-{L}aplacian priors.
\newblock {\em Adv Neural Inf Process Syst.}, 22:1033--1041, 2009.

\bibitem{lampe2012large}
J.~Lampe, L.~Reichel, and H.~Voss.
\newblock Large-scale {T}ikhonov regularization via reduction by orthogonal
  projection.
\newblock {\em Linear algebra and its applications}, 436(8):2845--2865, 2012.

\bibitem{lange2016mm}
K.~Lange.
\newblock {\em {MM} optimization algorithms}.
\newblock SIAM, 2016.

\bibitem{li2005motion}
T.~Li, E.~Schreibmann, Y.~Yang, and L.~Xing.
\newblock Motion correction for improved target localization with on-board
  cone-beam computed tomography.
\newblock {\em Phys. Med. Biol.}, 51(2):253, 2005.

\bibitem{li2010multichannel}
Y.~Li and R.~Verma.
\newblock Multichannel image registration by feature-based information fusion.
\newblock {\em IEEE Trans. on Medical Imaging}, 30(3):707--720, 2010.

\bibitem{lou2015weighted}
Y.~Lou, T.~Zeng, S.~Osher, and J.~Xin.
\newblock A weighted difference of anisotropic and isotropic total variation
  model for image processing.
\newblock {\em SIAM J. Imaging Sci.}, 8(3):1798--1823, 2015.

\bibitem{meaney2018tomographic}
A.~Meaney, Z.~Purisha, and S.~Siltanen.
\newblock Tomographic {X}-ray data of 3{D} emoji.
\newblock {\em arXiv preprint arXiv:1802.09397}, 2018.

\bibitem{meyer1992wavelets}
Y.~Meyer.
\newblock {\em Wavelets and Operators: {V}olume 1}.
\newblock Number~37. Cambridge university press, 1992.

\bibitem{mueller2012linear}
J.~L. Mueller and S.~Siltanen.
\newblock {\em Linear and nonlinear inverse problems with practical
  applications}.
\newblock SIAM, 2012.

\bibitem{mumford1985boundary}
D.~Mumford and J.~Shah.
\newblock Boundary detection by minimizing functionals.
\newblock In {\em Proc. IEEE Comput. Soc. Conf. Comput. Vis. Pattern
  Recognit.}, volume~17, pages 137--154. San Francisco, 1985.

\bibitem{pal2015brief}
C.~Pal, A.~Chakrabarti, and R.~Ghosh.
\newblock A brief survey of recent edge-preserving smoothing algorithms on
  digital images.
\newblock {\em arXiv preprint arXiv:1503.07297}, 2015.

\bibitem{patel2013sparse}
V.~M. Patel and R.~Chellappa.
\newblock {\em Sparse representations and compressive sensing for imaging and
  vision}.
\newblock Springer Science \& Business Media, 2013.

\bibitem{potts1952some}
R.~B. Potts.
\newblock Some generalized order-disorder transformations.
\newblock In {\em Mathematical proceedings of the Cambridge philosophical
  society}, volume~48, pages 106--109. Cambridge University Press, 1952.

\bibitem{rohde2003comprehensive}
G.~K. Rohde, S.~Pajevic, C.~Pierpaoli, and P.~J. Basser.
\newblock A comprehensive approach for multi-channel image registration.
\newblock In {\em International Workshop on Biomedical Image Registration},
  pages 214--223. Springer, 2003.

\bibitem{rudin1994total}
L.~I. Rudin and S.~Osher.
\newblock Total variation based image restoration with free local constraints.
\newblock In {\em Proceedings of 1st International Conference on Image
  Processing}, volume~1, pages 31--35. IEEE, 1994.

\bibitem{rudin1992nonlinear}
L.~I. Rudin, S.~Osher, and E.~Fatemi.
\newblock Nonlinear total variation based noise removal algorithms.
\newblock {\em Physica D: nonlinear phenomena}, 60(1-4):259--268, 1992.

\bibitem{saibaba2020efficient}
A.~K. Saibaba, J.~Chung, and K.~Petroske.
\newblock Efficient {K}rylov subspace methods for uncertainty quantification in
  large {B}ayesian linear inverse problems.
\newblock {\em Numer. Linear Algebra with Appl.}, 27(5):e2325, 2020.

\bibitem{schmitt2002efficient}
U.~Schmitt and A.~Louis.
\newblock Efficient algorithms for the regularization of dynamic inverse
  problems: {I}. {T}heory.
\newblock {\em Inverse Problems}, 18(3):645, 2002.

\bibitem{schmitt2002efficient2}
U.~Schmitt, A.~K. Louis, C.~Wolters, and M.~Vauhkonen.
\newblock Efficient algorithms for the regularization of dynamic inverse
  problems: {II}. {A}pplications.
\newblock {\em Inverse Problems}, 18(3):659, 2002.

\bibitem{semerci2014tensor}
O.~Semerci, N.~Hao, M.~E. Kilmer, and E.~L. Miller.
\newblock Tensor-based formulation and nuclear norm regularization for
  multienergy computed tomography.
\newblock {\em IEEE Trans. Image Process.}, 23(4):1678--1693, 2014.

\bibitem{sidky2008image}
E.~Y. Sidky and X.~Pan.
\newblock Image reconstruction in circular cone-beam computed tomography by
  constrained, total-variation minimization.
\newblock {\em Phys. Med. Biol.}, 53(17):4777, 2008.

\bibitem{stankovic2003haar}
R.~S. Stankovi{\'c} and B.~J. Falkowski.
\newblock The {H}aar wavelet transform: its status and achievements.
\newblock {\em Comput. Electr. Eng.}, 29(1):25--44, 2003.

\bibitem{stephane1999wavelet}
M.~Stephane.
\newblock A wavelet tour of signal processing, 1999.

\bibitem{stromberg1980weights}
J.-O. Str{\"o}mberg and A.~Torchinsky.
\newblock Weights, sharp maximal functions and {H}ardy spaces.
\newblock {\em Bulletin (New Series) of the American Mathematical Society},
  3(3):1053--1056, 1980.

\bibitem{sun2011image}
D.~Sun and M.~Ho.
\newblock Image segmentation via total variation and hypothesis testing
  methods.
\newblock 2011.

\bibitem{tan2018digital}
L.~Tan and J.~Jiang.
\newblock {\em Digital signal processing: fundamentals and applications}.
\newblock Academic Press, 2018.

\bibitem{thompson1985dynamic}
W.~B. Thompson, K.~M. Mutch, and V.~A. Berzins.
\newblock Dynamic occlusion analysis in optical flow fields.
\newblock {\em IEEE Trans. Pattern Anal. Mach. Intell.}, (4):374--383, 1985.

\bibitem{trabold2003estimation}
T.~Trabold, M.~Buchgeister, A.~K{\"u}ttner, M.~Heuschmid, A.~Kopp,
  S.~Schr{\"o}der, and C.~Claussen.
\newblock Estimation of radiation exposure in 16-detector row computed
  tomography of the heart with retrospective {ECG}-gating.
\newblock In {\em R{\"o}Fo-Fortschritte auf dem Gebiet der R{\"o}ntgenstrahlen
  und der bildgebenden Verfahren}, volume 175, pages 1051--1055. {\copyright}
  Georg Thieme Verlag Stuttgart{\textperiodcentered} New York, 2003.

\bibitem{vogel2002computational}
C.~R. Vogel.
\newblock {\em Computational methods for inverse problems}.
\newblock SIAM, 2002.

\bibitem{vogel1998fast}
C.~R. Vogel and M.~E. Oman.
\newblock Fast, robust total variation-based reconstruction of noisy, blurred
  images.
\newblock {\em IEEE Trans. Image Process.}, 7(6):813--824, 1998.

\bibitem{SSIM}
Z.~Wang, A.~C. Bovik, H.~R. Sheikh, and E.~P. Simoncelli.
\newblock Image quality assessment: {F}rom error visibility to structural
  similarity.
\newblock {\em IEEE Trans. Image Process.}, 13(4):600--612, 2004.

\bibitem{westman2012passive}
E.~Westman, K.~Luxbacher, and S.~Schafrik.
\newblock Passive seismic tomography for three-dimensional time-lapse imaging
  of mining-induced rock mass changes.
\newblock {\em The Leading Edge}, 31(3):338--345, 2012.

\bibitem{beck2012process}
R.~Williams and M.~S. Beck.
\newblock {\em Process tomography: principles, techniques and applications}.
\newblock Butterworth-Heinemann, 2012.

\bibitem{wohlberg2007tv}
B.~Wohlberg and P.~Rodriguez.
\newblock An iteratively reweighted norm algorithm for minimization of total
  variation functionals.
\newblock {\em IEEE Signal Process. Lett.}, 14(12):948--951, 2007.

\bibitem{xu2012l}
Z.~Xu, X.~Chang, F.~Xu, and H.~Zhang.
\newblock $ l\_ $\{$1/2$\}$ $ regularization: A thresholding representation
  theory and a fast solver.
\newblock {\em IEEE Trans. Neural Netw. Learn. Syst.}, 23(7):1013--1027, 2012.

\bibitem{yang2005adaptive}
X.~Yang, S.~Yao, K.~P. Lim, X.~Lin, S.~Rahardja, and F.~Pan.
\newblock An adaptive edge-preserving artifacts removal filter for video
  post-processing.
\newblock In {\em 2005 IEEE International Symposium on Circuits and Systems},
  pages 4939--4942. IEEE, 2005.

\bibitem{yin2007total}
W.~Yin, D.~Goldfarb, and S.~Osher.
\newblock The total variation regularized ${L}^1$ model for multiscale
  decomposition.
\newblock {\em Multiscale Model. Simul.}, 6(1):190--211, 2007.

\bibitem{zhang2009passive}
H.~Zhang, S.~Sarkar, M.~N. Toks{\"o}z, H.~S. Kuleli, and F.~Al-Kindy.
\newblock Passive seismic tomography using induced seismicity at a petroleum
  field in {O}man.
\newblock {\em Geophysics}, 74(6):WCB57--WCB69, 2009.

\end{thebibliography}
\end{document}